\documentclass[12pt]{amsart}
 \usepackage{amssymb}
 \usepackage{amsmath}
\usepackage{amsfonts}
\usepackage[monochrome]{xcolor}
\usepackage{comment}
 \usepackage{eucal}
      \makeatletter
     \def\section{\@startsection{section}{1}%
     \z@{.7\linespacing\@plus\linespacing}{.5\linespacing}%
     {\bfseries
     \centering
     }}
   \def\subsection{\@startsection{subsection}{1}%
     \z@{.1\linespacing\@plus\linespacing}{.5\linespacing}%
     {\bfseries
     \centering
     }}
    \def\@secnumfont{\bfseries}
     \makeatother
  \usepackage{a4wide}

\newtheorem{theorem}{Theorem}[section]
\newtheorem{lemma}[theorem]{Lemma}
\newtheorem{corollary}[theorem]{Corollary}
\newtheorem{proposition}[theorem]{Proposition}

\theoremstyle{definition}
\newtheorem{definition}[theorem]{Definition}
\newtheorem{problem}[theorem]{Problem}
\newtheorem{example}[theorem]{Example}
\newtheorem{remark}[theorem]{Remark}

\numberwithin{equation}{section}

\def \a{{\alpha}}
\def \b{{\beta}}
\def \D{{\Delta}}

\def \d{{\delta}}
\def \e{{\varepsilon}}

\def \g{{\gamma}}
\def \G{{\Gamma}}
\def \k{{\kappa}}
\def \l{{\lambda}}

\def \o{{\omega}}

\def \p{{\varphi}}
\def \t{{\vartheta}}

\def \m{{\mu}}
\def \s{{\sigma}}

\def \C{{\cal C}}

\def \E{{\bf E}\, }

\def \N{{\bf N}}

\def \P{{\bf P}}

\def \qq{{\qquad}}
\def \R{{\bf R}}

\def \Z{{\bf Z}}

\def \dd{{\rm d}}

\def \noi{{\noindent}}

\def\beq{\begin{equation}}
\def\eeq{\end{equation}}

\def\ben{\begin{eqnarray}}
\def\een{\end{eqnarray}}

\def\E{{\mathbb E \,}}

\def\P{{\mathbb P}}

\def\R{{\mathbb R}}

\def\Z{{\mathbb Z}}

\def\N{{\mathbb N}}

\def\C{{\mathbb C}}

 \scrollmode
 
   \font\sevenrm= cmr10 at 7 pt
   \font\ph=cmcsc10 at 14 pt
    \font\ph=cmcsc10 at 14 pt
\font\sevenit=cmti7
\font\sevenrm= cmr10 at 7,3 pt \font\eightit=cmti9
 \font\ninerm= cmr10 at 9 pt
 \font\eightrm= cmr10 at 8 pt
 
\font\mog=cmbx10 at 24 pt
\font\tf=cmbxsl10 at 11pt

\scrollmode

\def\ddate {\sevenrm \ifcase\month\or January\or
February\or March\or April\or May\or June\or July\or
August\or September\or October\or November\or December\fi\! {\the\day}, \!{\sevenrm\the\year}}

\scrollmode

\title[LOCAL LIMIT THEOREMS   AND ALMOST SURE VERSIONS]{\mog Classical and Almost Sure \\ \vskip 30 pt  Local Limit Theorems}

\begin{document}

\maketitle
\vskip 50 pt
\begin{center} {\ph  Zbigniew  Szewczak and Michel  Weber}
\end{center}
  \date{\today}
\vfill \vskip 300 pt
 {
{\it Key words and phrases}: almost sure limit theorems, local limit theorem,
lattice distribution
\vskip 1pt  2010 AMS  {\it Mathematics Subject Classification}: 60F05, 60F15, 60G50, 60J05, 11K50.}
\vfill\break
$${}$$
\vskip 30 pt$$\hbox{\bf Preface}$$
 \vskip 20pt
 {\leftskip = 2cm \rightskip= 2cm
  \noi \ninerm In this
monograph, we present and discuss the many results obtained concerning a famous limit theorem,
 the local limit theorem (LLT),  which has many interfaces, with Number Theory notably, and for which,
 in spite of considerable efforts,  the question concerning conditions of validity of the local limit theorem,
 has up to now no satisfactory solution. These results mostly concern sufficient conditions for the validity
 of the LLT and its interesting variant forms:
strong LLT, strong LLT with convergence in variation.
 Quite importantly are necessary conditions, and the results obtained are sparse,
essentially: Rozanov's necessary condition, Gamkrelidze's necessary condition, and, almost isolated among the
 flow of results, Mukhin's necessary and sufficient condition.
Extremely useful and instructive  are the counter-examples due to Azlarov and Gamkrelidze,
as well as necessary and sufficient conditions obtained for a class of random variables,
 such as Mitalauskas' characterization of the LLT in the strong form for random variables
   having stable limit distributions. The method of characteristic functions and the Bernoulli   part extraction
 method,  are presented and compared.
   The study of the LLT (old and recent results) consists of three parts:}
  \vskip 1 pt
    {  \leftskip = 2cm \rightskip= 2cm \ninerm The LLT for sums of independent, identically distributed random variables:  Gnedenko's Theorem, Ibragimov and Linnik
    characterization of the speed of convergence under moments conditions, stronger forms, Galstyan results, strong LLT's with convergence in variation, versions for
    densities,   the case of  weighted sums of i.i.d. random variables, local large deviations, Nagaev's result, Tkachuk and Doney's results,  Diophantine measures and LLT,
    Breuillard, Mukhin, Shepp, Stone, LLT and Edgeworth expansions, Breuillard, Feller, LLT's under arithmetical conditions.}
 \vskip 1 pt
  {  \leftskip = 2cm \rightskip= 2cm \ninerm The LLT for sums of independent  random variables: Prohorov's theorem, Richter LLT's  and large deviations, Maejima's LLT's with
  remainder term, LLT's with convergence in variation, Gamkrelidze's results, Rozanov's necessary condition and uniform asymptotic  distribution, Mitalauskas' LLT  for random
  variables having stable limit distribution,   Azlarov and  Gamkrelidze's counterexamples,  structural characteristics, Bernoulli part extraction, Dabrowski and McDonald's
  LLT, Giuliano and Weber's LLT with effective rate, Macht and Wolf's LLT using   H\"older-Continuity, R\"ollin  and Ross LLT's using Landau-Kolmogorov inequalities,
  Delbaen, Jacod, Kowalski and Nikeghbali recent LLT's under Mod-$ \phi$ convergence,  Dolgopyat's recent LLT    for sums of independent random vectors satisfying appropriate
  tightness assumptions,  Feller's LLT's and domains of attraction.}
\vskip 1 pt
 {  \leftskip = 2cm \rightskip= 2cm \ninerm The LLT for ergodic sums: essentially centered around interval expanding maps, results of Kac, Rousseau-Egele, Broise, Calderoni,
 Campanino and Capocaccia notably, and more recently Gou\"ezel, Szewczak.}

 \vskip 1 pt
 {  \leftskip = 2cm \rightskip= 2cm \ninerm An expanded and detailed  list of applications of the local limit theorem is provided.  The last part of the survey is devoted to
 the more recent study
 of the almost sure local limit theorem,   instilled by Denker and Koch.  The inherent second order study,
 which has its own interest,  is much more difficult than for establishing the almost sure central limit theorem.
 The   almost sure local limit theorems established already cover  the i.i.d. case, the stable case,
 Markov chains,   the model of the Dickman function,  and  the independent case, with almost
sure convergence of related series.}
\vskip 5 pt
{  \leftskip = 2cm \rightskip= 2cm \ninerm    Our aim in writing this
 monograph was  to survey and
   bring to knowledge many interesting results obtained since the sixties and up to now. Many of them were obtained during the sixties by the Lithuanian and Russian  Schools
   of Probability, and are  essentially written in Russian, and moreover often published in   Journals of   difficult  access.
In doing so, our feeling was  to   somehow help  with this whole coherent body of results and methods, researchers in the study of the local limit theorem,  at least it is
our hope.
}


\vskip 15 pt
{\leftskip 2cm \rightskip= 2cm  \eightit    Dedicated to the  Lithuanian and Russian  Schools of Probability for their important contributions in the study and the
applications of the local limit theorem,  in particular due to Azlarov,
 Gamkrelidze, Gnedenko, Ibragimov, Kolmogorov, Kubilius, Linnik, Matskyavichyus, Mitalauskas, Moskvin, Mukhin, Nagaev,  Petrov, Postnikov, Prokhorov,  
 Rva\v ceva, Richter, Rogozin, Rozanov, {\color{orange} Sirazhdinov,} {\color{orange} Statulyavichus,} Survila, Tkachuk. Let us also mention the contributions of
 {\color{orange} Cram\'er,} Dabrowski,
  Dolgopyat,
  Doney, Giuliano, Kac, Lasota, Yorke, Maejima, McDonald, Petit, Rousseau-Egele, Shepp, Stone, Szewczak, Weber, \ldots}

\vskip 35 pt
\noi {\bf Acknowledgments.} {\leftskip 2cm \rightskip= 2cm  \rm  We   thank Prs Nicko Gamkrelidze, Jan Rosinski and Jordan Stoyanov for their warm appreciations.}  

\vfill\break

 \tableofcontents




{\color{blue}

\section{The local limit theorem}
\label{1}

    A simple way to introduce to this topic
is to begin with citing McDonald \cite{M2} p.\,73: \lq\lq{\it The local limit theorem describes how the density of a sum of random variables follows the normal curve.
However the local limit theorem is often seen as a curiosity of no particular importance when compared with the central limit theorem.
Nevertheless the local limit theorem came first and is in fact associated with the foundation of probability theory by Blaise Pascal and Pierre de Fermat and was originally
formalized by Jakob Bernoulli, Abraham De Moivre and Pierre-Simon Laplace.}\rq\rq

\vskip  5 pt The central limit theorem  is
since long time a major piece of the Theory of Probability, with numerous applications, and  is  still actively studied.
The  local limit theorem  is comparatively much less investigated.
One reason is that the study of this  finer limit theorem  is intrinsically more complicated.
For instance, conditions on arithmetical properties of the support of a random variable are always present.

\vskip  3 pt On the other hand, it can be noted that although there  exists  a substantial literature   on the central limit theorem for sums of dependent random variables,
especially martingales and stationary sequences, there is a paucity of local versions of such works.

\vskip   3 pt Nevertheless, the field of application of the local limit theorem in combinatorics analysis   and number theory notably,  is considerable.
 {\color{green} The    interface between local limit theorems and structure theory of set addition  was much studied by  Freiman, Moskvin and Yudin notably.}
 The interaction with Number theory was  investigated by  the  Lithuanian and Russian  Schools of Probability, during the Soviet period, and later  by   Manstavichyus,
 Postnikov notably and Fomin, after earlier works of Gamkrelidze,  Mitalauskas, Mukhin and Raudelyunas, and many other contributors.

\vskip   9 pt
One can define simply  in the multi-dimensional case the local limit theorem. Let $S_n= \sum_{k=1}^n X_{nk}$, where $\{X_{nk}, k=1,\ldots,n\}$ are series of independent
random variables with values in   $\Z^s$, and such  that the integral limit theorem  holds:  there exist $A_n\in \R^s$ and real $B_n\to \infty$ such that the sequence of
distributions of $(S_n-A_n)/B_n$ {\color{green}weakly} converges   to an absolutely continuous distribution with density $g(x)$, which is uniformly continuous in $\R^s$.
 \vskip 3 pt Then the local limit theorem is valid if
\beq\label{ilt.llt.0}
\P\{S_n=m\}=B_n^{-s} g\Big(\frac{m-A_n}{B_n}\Big) + o(B_n^{-s}),
\eeq
uniformly in $m\in \Z^s$.
\vskip 5 pt

Mukhin  remarked in his 1991's paper \cite{Mu1}  that \lq\lq {\it In spite that the problem was posed a long time ago, the formulation is very simple, and there exists a
great number of investigations, the question concerning conditions of validity of the local limit theorem, has up to now no satisfactory solution. {\color{green} The
sufficient conditions for the validity of  the l.l.t., as well as necessary and sufficient conditions, are useful for the verification only in some particular cases, for
example such as the case  of identically distributed summands. And even in such situations they give no clear idea about the class of all sequences satisfying the l.l.t.} The
main obstacle for the l.l.t. to be valid is a concentration of a great amount of distributions of summands on a lattice whose span is greater than 1.}\rq\rq\
These remarks are  still relevant.

\vskip 8 pt
\vskip   3 pt
notably and Fomin, after earlier works of Gamkrelidze,  Mitalauskas, Mukhin and Raudelyunas, and many other contributors.

The first local limit theorem    was  established  nearly  three centuries ago, in the binomial case, which is also the fundamental case.\vskip 8 pt

\noi {\tf De Moivre-Laplace's Theorem.}
 It was proved by De Moivre (17}{\color{red}33}{\color{blue}) in {\it Approximatio ad Summam Terminorum Binomii $(a + b)^n$ in Seriem expansi} expanding on James Bernoulli's
 work} {\color{red}(see \cite{Arc})}{\color{blue}.
For the case $p = 1/2$ he proved:

\begin{theorem}[{\rm De Moivre--Laplace, 1730}]
\label{moivre}  Let
$0<p<1$,
$q=1-p$. Let $X$ be such that
${\mathbb P}\{X=1\}=p=1-{\mathbb P}\{X=0\}
$. Let
$X_1, X_2,\ldots$ be independent copies of
$X$ and let $S_n=X_1+\ldots +X_n$, $n\ge 2$. Let also  $0<\g<1$.
Then for all $k$ such that $|k-np|  \le \g npq$ and $    n \ge \max(p/q,q/p)$, letting
$ x= \frac{k-np}{\sqrt{  npq}}$,
we have
 \begin{equation*} {\mathbb P}\{S_n=k\}  \ =\
 \frac{e^{-  \frac{x^2}{  2  } }}{\sqrt{2\pi npq}} \   e^E ,
 \end{equation*}
with
$|E|\le   \frac{3|x|+ 2 |x|^3}{(1-\g)\sqrt{ npq}} +  \frac{1}{ 4n\min(p,q)(1 -  \g      )}$.

\vskip   3 pt  If $x= o(n^{1/6})$,  then
\begin{equation*}\Big| {\mathbb P}\{S_n=k\}- \frac{1 }{\sqrt{2\pi npq}}\,  e^{-  \frac{(k-np)^2}{  2 npq }}\Big|
\ \le  \   \frac{1 }{\sqrt{2\pi npq}}\ e^{-  \frac{(k-np)^2}{  2 npq }}\Big( \frac{c_1 |x|^3}{\sqrt{ n}}+ \frac{c_2|x|}{\sqrt n} +  \frac{c_3}{ n}\Big).
\end{equation*}
The constants $c_i$, $i=1,2,3$ are explicit and depend on $p$ only. \end{theorem}
}{\color{red} De Moivre's result was further improved by Khintchin \cite{Kh29}
(see also  \cite[Remark on p. 209]{Ri}).}
{\color{blue}We refer  to McDonald \cite{M2} for a short historical study of the local limit theorem,  to Todhunter \cite[Chap.\,IX]{To}  for a full analysis of De Moivre's
\lq\lq {\it Doctrine of Chance\,}\rq\rq, also to Bellhouse and  Genest \cite{BG}} {\color{red}and Bellhouse \cite{Be11}}{\color{blue} for a biography of De Moivre.

  \smallskip\par  De Moivre stated the local limit theorem for general $p$ but proved only the case $p = 1/2$. The proof for all $p$ was provided by Laplace  (1795) in {\it
  Th\'eorie analytique des probabilit\'es}.  This is why the above theorem is called the De Moivre--Laplace central limit theorem. As it concerns densities, it can be
  interpreted as a local limit theorem.

\begin{proof}[Proof of Theorem \ref{moivre}]
The proof  follows   the one given in Chow and Teicher \cite[p.\,46]{CT}.  Let $b(k,n,p)={\mathbb P}\{S_n=k\}$.  By using Stirling's formula,
\begin{equation} \label{stirling} n!= n^{n+\frac{1}{2}}e^{-n+\e_n}\sqrt{2\pi}  \qq {\rm with}\qq \frac{1}{12n +1}<\e_n <\frac{1}{12n } ,
\end{equation}
valid for all   $n\ge 2$, we have
\begin{eqnarray*}b(k,n,p)&=& p^kq^{n-k}\frac{ n!}{ k!(n-k)!}\ =\  \frac{p^kq^{n-k} n^{n+\frac1{2}}e^{\e_n-\e_k-\e_{n-k}} }
{\sqrt{2\pi}\, k^{k+\frac1{2}}(n-k)^{ n-k +\frac1{2}}  }\cr
&=&\frac{e^\e}{\sqrt{2\pi}}\Big(\frac{np}{k} \Big)^{k+ \frac1{2}}\Big(\frac{nq}{n-k} \Big)^{n-k+ \frac1{2}} (npq)^{-\frac1{2}},
\end{eqnarray*}
where $\e= \e_n-\e_k-\e_{n-k}$.

\vskip   3 pt
Whence (recalling that $x= \frac{k-np}{\sqrt{  npq}}$),
\begin{align*}
& \log \big(\sqrt{2\pi npq}\,b(k,n,p)\big)
\ =\
\e -(k+\frac1{2})\log \frac{k}{np} -  (n-k+ \frac1{2})\log \frac{n-k}{nq}
\cr =& \ \e-(np +x\sqrt{npq}+\frac1{2})\log \big(1+
\frac{x\sqrt q}{\sqrt{np}}\big)
 - (nq -x\sqrt{npq}+\frac1{2})\log
\big(1-\frac{x\sqrt p}{\sqrt{nq}}\big)
\cr =& \ \e-(np +x\sqrt{npq})\log \big(1+\frac{x\sqrt q}{\sqrt{np}}\big)
 - (nq -x\sqrt{npq})\log
\big(1-\frac{x\sqrt p}{\sqrt{nq}}\big) - \frac L{2},
\end{align*}
with $L=\log  \frac{1+x\sqrt{{q}/{np}}}{1-x\sqrt{{p}/{nq}}} $.

\vskip   3 pt
 Let $0<\g<1$.   Let $\p(y) = \log (1+y) -y +y^2/2$.   Then $\p'(y)= \frac{y^2}{1+y}$. Let $|y|\le \g$. By the mean value theorem,  $|\p(y)|=|\p(y)-\p(0)|\le
 |y|\frac{\xi^2}{1+\xi}$ for some $\xi \in ]-y,y[$, and so $|\p(y)|\le  \frac{|y|^3}{1-\g}$, namely
     \begin{eqnarray}\label{ineq1}\big|\log (1+y)-y+y^2/2 \big|\le  \frac{|y|^3}{ 1-\g}, \qq |y|\le \g.
  \end{eqnarray}
 Since $ \frac{|x|}{\sqrt n}
\max ( \sqrt{ {q}/{ p}}, \sqrt{ {p}/{ q}} ) \le   \frac{|x|}{\sqrt{ npq}  } \le \g$, by the assumption made,  we deduce that
 \begin{eqnarray*}\log \big(\sqrt{2\pi npq}\,b(k,n,p)\big)  &=& \e-(np +x\sqrt{npq})  \Big[ x\sqrt{\frac{q}{np}}-
\frac{x^2q}{2np} +A \Big]
\cr & &\quad + (nq -x\sqrt{npq})
\Big[ x\sqrt{\frac{p}{nq}} +\frac{x^2p}{2nq} +B\Big]- \frac L{2},
\end{eqnarray*}
with
$$|A|  \le \frac{1}{1-\g}\frac{|x|^3} { n^{\frac{3}{2}}}\big( \frac{q}{ p}\big)^{\frac{3}{2}}, \qq  |B|  \le
\frac{1}{1-\g}\frac{|x|^3} { n^{\frac{3}{2}}}\big( \frac{p}{ q}\big)^{\frac{3}{2}}.$$
Whence,
\begin{eqnarray*} \log \big(\sqrt{2\pi npq}\,b(k,n,p)\big) &=&  \e -\frac{x^2}{2}+ \frac{x^3q^{\frac{3}{2}}}{
2\sqrt{np}}-\frac{x^3p^{\frac{3}{2}}}{ 2\sqrt{nq}}- \frac L{2} +C ,
\end{eqnarray*}
where
\begin{eqnarray*} C&=& - (np +x\sqrt{npq}+\frac1{2})A +(nq -x\sqrt{npq}+\frac1{2}) B.
\end{eqnarray*}

  We have
$$np \big( \frac{q}{
p}\big)^{\frac{3}{2}}+nq\big( \frac{p}{ q}\big)^{\frac{3}{2}}\le \frac{n}{\sqrt{ pq}}, \qq \sqrt{ pq}\Big(\big( \frac{q}{
p}\big)^{\frac{3}{2}}+\big( \frac{p}{ q}\big)^{\frac{3}{2}}\Big)\le \frac{1}{pq}.$$
Thus \begin{eqnarray*} (1-\g) |C|
&\le & \frac{|x|^3}{n^{\frac{3}{2}}}\Big\{ \big(np +|x|\sqrt{npq}\big)\big( \frac{q}{
p}\big)^{\frac{3}{2}} +\big(nq +|x|\sqrt{npq}\big) \big( \frac{p}{ q}\big)^{\frac{3}{2}}\Big\}
\cr
&\le & \frac{|x|^3}{n^{\frac{3}{2}}}\Big\{\frac{n}{\sqrt{ pq}}+ \frac{|x|\sqrt n}{pq} \Big\}
\ =\  \frac{ |x|^3}{\sqrt{ npq}}+ \frac{|x|^4 }{npq}.
\end{eqnarray*}
Further
$  \frac{x^3q^{\frac{3}{2}}}{
2\sqrt{np}}-\frac{x^3p^{\frac{3}{2}}}{ 2\sqrt{nq}}
= \frac{x^3  (
q-p )}{ 2\sqrt{n pq}}
. $
We get
\begin{eqnarray}\label{estintermoivre}
\log \big(\sqrt{2\pi npq}\,b(k,n,p)\big)  &=&  \e - \frac{x^2}{2} -\frac{L}{2}+D ,
\end{eqnarray}
with
\begin{eqnarray}\label{estintermoivre1}|D|\le \frac{1}{(1-\g)}\Big\{ \frac{ |x|^3}{ \sqrt{ npq}}+ \frac{|x|^4 }{npq}\Big\}\ \le \  \frac{2}{(1-\g)} \frac{ |x|^3}{ \sqrt{
npq}},
\end{eqnarray}
since $|x| \le \g \sqrt{npq}$.

Now by \eqref{ineq1}, we have $\big|\log (1+y)\big|\le  \frac{3|y|}{ 1-\g}$ if $|y|\le \g$. Thus
$$|L|=\Big|\log  \frac{1+x\sqrt{{q}/{np}}}{1-x\sqrt{{p}/{nq}}} \Big|\le \frac{3|x|}{ (1-\g)\sqrt n}\big(\sqrt{q/p}+\sqrt{p/q}\big)=\frac{3|x|}{ (1-\g)\sqrt{ npq}}. $$

 By assumption $|k-np| \le \g npq$ and $    n \ge \max(p/q,q/p)$. Thus
$$  k\ge np(1 -\g q)  \ge np(1-\g) ,$$
 and   $   0\le \e_k \le \frac{1}{12k }\le \frac{1}{12np(1 -  \g      ) }$. Moreover,
 $$  n-k\ge nq-  \g npq \ge  nq(1-\g) ,$$
so that  $   0\le \e_{n-k} \le \frac{1}{12(n-k )}\le
\frac{1}{12nq(1 -  \g      ) }$. Recalling that $\e= \e_n-\e_k-\e_{n-k}$, we further have
 $$   \e \le  \frac{1}{ 6n\min(p,q) (1 -  \g      ) }+\frac{1}{12n }\le \frac{1}{ 4n\min(p,q)(1 -  \g      )}.  $$

We deduce from \eqref{estintermoivre1} that,
\begin{eqnarray*}\log \big(\sqrt{2\pi npq}\,b(k,n,p)\big)  &=&   - \frac{x^2}{2} +E ,
\end{eqnarray*}
with
$$|E|\le   \frac{1}{(1-\g)}\Big\{\frac{3|x|}{ 2\sqrt{ npq}} +\frac{2 |x|^3}{ \sqrt{ npq}}\Big\}+  \frac{1}{ 4n\min(p,q)(1 -  \g      )}\, .
$$ Consequently,
\begin{eqnarray*} {\mathbb P}\{S_n=k\}  &=&   \frac{e^{-  \frac{(k-np)^2}{  2 npq }} }{\sqrt{2\pi npq}}\   e^E .
\end{eqnarray*}

 Thus  if $x= o(n^{1/6})$,
\begin{eqnarray*}& &\Big| \P\{S_n=k\}-
\frac{e^{-
 \frac{(k-np)^2}{  2 npq }
 } }{\sqrt{2\pi npq}}  \Big|
\cr &\le &
\frac{
e^{-
\frac{(k-np)^2}{  2 npq }
 }
  }{\sqrt{2\pi npq}}
\Big\{ \, \frac{1}{(1-\g)}
\big(\frac{3|x|}{ 2\sqrt{ npq}} +\frac{2 |x|^3}{ \sqrt{ npq}}\big)+  \frac{1}{ 4n\min(p,q)(1 -  \g      )}\, \Big\}
\cr &:= &   \frac{e^{-
\frac{(k-np)^2}{  2 npq }
 } }
 {\sqrt{2\pi npq}}\Big( \frac{c_1 |x|^3}{\sqrt{ n}}+ \frac{c_2|x|}{\sqrt n} +  \frac{c_3}{ n}\Big).
\end{eqnarray*}
\end{proof}

 \rm The case of independent non-identically distributed binomial random variables was considered by Mamatov in \cite{Ma}.

 \smallskip\par
It is natural to compare the De Moivre--Laplace theorem with  nowadays classical local limit theorems.    The following precise result concerning  the case $p=q=1/2$, is
derived from a fine local limit theorem with asymptotic expansion \cite[Ch.\,7,\,Th.\,13]{P}.

\begin{theorem} \label{lltber}Let
$\mathcal B_n=\b_1+\ldots+\b_n$, $n=1,2,\ldots$  where
$
\b_i
$ are i.i.d.  Bernoulli r.v.'s (\,${\mathbb P}\{\b_i=0\}={\mathbb P}\{\b_i=1\}=1/2$).  There exists a numerical constant $C_0$ such that for all positive $n$
 \begin{eqnarray*}  \sup_{k}\, \Big|  {\mathbb P}\big\{\mathcal B_n=k\}
 -\sqrt{\frac{2}{\pi n}} e^{-{ (2k-n)^2\over
2 n}}\Big|\ \le \
  \  \frac{C_0}{n^{3/2}}  .
  \end{eqnarray*}
\end{theorem}

  One easily observes that    for
moderate deviation of type $x= \frac{|k-n/2|}{\sqrt{  n}/2}\sim n^{1/7}$,  this result is however  considerably less precise than the old    one of De Moivre (case
$p=q$).  Indeed,
\begin{eqnarray*} {\mathbb P}\{\mathcal B_n=k\}  -\sqrt{\frac{  2} {  \pi n }} e^{-  \frac{x^2}{  2  }}
&=&   \mathcal O\big(n^{-1/28}e^{-n^{1/7}}\big)\ \ll \ \mathcal O\big(n^{-3/2}\big).\end{eqnarray*}

 This  in particular   shows that
the usual formulation of the local limit theorem  does not   provide in fact the full information on the remainder term.}

\vskip 15 pt
{\color{blue} \noi{\tf Lattice-valued random variables.} We recall some classical facts  (\cite{P},\,p.\,10).
\vskip 3 pt
(I) A random variable $X$ has lattice distribution if and only if its characteristic function $f(t)=\E e^{it X}$ satisfies $f(t_0)=1$ for some $t_0\neq 0$.
\vskip 3 pt
Let $X$ be taking values in the lattice  $\mathcal L(v_{
0},D )$, namely defined by the
 sequence $v_{ k}=v_{ 0}+D k$, $k\in \Z$, where
 $v_{0} $ and $D >0$ are   real numbers.

 \vskip   3 pt  The span  $D$ is said to be {\it maximal},
 if   there are no  other real numbers
$v'_{0}
$ and
$D' >D$ for which
${\mathbb P}\{X
\in\mathcal L(v'_0,D')\}=1$.
\vskip 3 pt
(II) Let the span $D$ be maximal. Then,
$$ |f\big(\frac{2\pi}{D}\big)|=1, \qq {\rm and}\qq
  |f(t)|<1\ \ {\rm for}\   0<|t|<\frac{2\pi}{D}.$$
As a consequence, for every  $0<\e <\frac{2\pi}{D}$, there exists and $\theta>0$,  such that  $|f(t)|\le e^{-\theta}$, for $\e\le |t|\le \frac{2\pi}{D}-\e$.

\vskip   20 pt Let $ \{X_n , n\ge 1\}$  be a sequence of independent, square integrable
random variables taking values in a common lattice $\mathcal L(v_{
0},D )$, and let $S_n=\sum_{j=1}^nX_j$, for each $n$.
 Then  $S_n$  takes values in the lattice
$\mathcal L( v_{ 0}n,D )$.
   Put
\begin{equation}\label{not1} M_n= {\mathbb E\,} S_n
 ,  \qq B_n^2 ={\rm Var }(S_n).
\end{equation}

\begin{definition}\label{defllt} The sequence $ \{X_n , n\ge 1\}$ satisfies a local limit theorem if
 \begin{equation}\label{llt}  \D_n:=  \sup_{N=v_0n+Dk }\Big|B_n {\mathbb P}\{S_n=N\}-{D\over  \sqrt{ 2\pi } }e^{-
{(N-M_n)^2\over  2 B_n^2} }\Big| = o(1).
\end{equation}
\end{definition}
This is a fine limit theorem in Probability Theory, which also has
 connections with Number Theory.
  These two aspects of a same problem were  much studied in the past decades by the
  Lithuanian and
  Russian probabilists.

  Note that  the transformation
\begin{equation}\label{llt.transf.}
 X'_j= \frac{X_j-v_0}{D},
  \end{equation}
allows one to reduce  to the case  $v_0=0$, $D=1$.

{\color{green}\subsection{The i.i.d. case}} \label{1.2}

 Assume that the random variables $X_i$ are identically distributed, and let $\m={\mathbb E\,} X_1$, $\s^2={\rm Var}(
X_1)$. Then $M_n=n\m$, $B_n^2=n\s^2$, and  \eqref{llt} reduces to
 \begin{equation}\label{llt.iid}  \D_n =  \sup_{N=v_0n+Dk }\Big|  \s \sqrt{n}\, {\mathbb P}\{S_n=N\}-{D\over  \sqrt{ 2\pi } }e^{-
{(N-n\m)^2\over  2 n\s^2} }\Big| = o(1).
\end{equation}

\noi {\tf Gnedenko's Theorem.} The following  well-known result of Gnedenko \cite{G}    characterizes the local limit theorem in this case. It was proved in 1948 and can be
viewed as a generalization
of the De Moivre--Laplace theorem.

\begin{theorem}
\label{gnedenko}    Assume that  $ \{X_n , n\ge 1\}$ is an i.i.d.\  sequence and let $\m={\mathbb E\,} X_1$, $\s^2={\rm Var}(
X_1)$.  Then  \eqref{llt.iid} holds
if and only if the span $D$ is maximal.
\end{theorem}

{\color{green} The} proof was effected by {\color{green}means of} the
method of characteristic functions. See \cite{G1} or \cite{P} (Theorem 1 and proof of Theorem
2,  p.\ 193--195).

  \vskip 3 pt

  Remark
that   (\ref{llt.iid}) is significant only for the bounded domains of values
 \begin{equation}\label{lltrange}   {|N-n\m|} \le \s \sqrt{2n\log \frac{D }{ \e_n  }}   ,
 \end{equation}
where $\e_n\downarrow 0$ depends on the Landau symbol $o$.

   \vskip 3 pt
   It is {\color{green}also} worth observing that   (\ref{llt.iid}) cannot be deduced from a central limit
theorem with rate, even under   stronger moment assumption.  Suppose for instance that
$D=1$,
$X$ is centered and ${\mathbb E\,}| X|^3<\infty$. Using the Berry--Esseen  estimate   implies that
 $$\Big|\s  \sqrt{   n}  {\mathbb P}\{S_n=k\}-\s \sqrt n\int_{\frac{k}{\s \sqrt n}}^{\frac{k+1}{\s \sqrt
n}}e^{-t^2/2}\frac{{\rm d} t}{\sqrt{2\pi}}\Big|\le C\frac{ {\mathbb E\,}| X|^3}{\s^2} .$$ Further,
$$ \sup_{k+1\le \s \sqrt n}\big|\s \sqrt n\int_{\frac{k}{\s \sqrt n}}^{\frac{k+1}{\s \sqrt n}}e^{-t^2/2}\frac{{\rm d}
t}{\sqrt{2\pi}}-\frac{1}{\sqrt{2\pi}}e^{-\frac{k^2}{2\s^2   n}}\big|\le \frac{C}{\s \sqrt n} \to 0. $$
By substituting we  get,
  $$\Big|\s  \sqrt{   n}  {\mathbb P}\{S_n=k\}-\frac{1}{\sqrt{2\pi}}e^{-\frac{k^2}{2\s^2   n}}\Big|\le C\big(\frac{ {\mathbb E\,}| X|^3}{\s^2}+\frac{1}{\s \sqrt n}\big) .$$
 Letting $k=k_n\to \infty$, $k_n+1\le \s \sqrt n$, the  right-hand side is  bounded with $n$, whereas by  (\ref{llt.iid}), this one tends to zero with $n$.  Hence
 (\ref{llt.iid}) cannot follow from the Berry--Esseen   estimate. Note however that   if Cram\'er's condition is fulfilled, namely
 \beq\label{cramer.cond}\hbox{$\limsup_{|u|\to \infty}  |{\mathbb E\,} e^{iu X_1} |<1$},
 \eeq and higher moments exist, better rates of approximation in the Berry--Esseen theorem are available, see
 \cite[p.\,329]{BHW}, \cite[p.\,130]{P}.

 \bigskip

 Gnedenko's theorem is optimal. This was proved   by Matskyavichyus in \cite{Mat}.

 \begin{theorem}
 For any nonnegative function
$\p(n)\to 0$ as $n\to \infty$, there is a sequence $ \{X_n , n\ge 1\}$ of   i.i.d. integer valued random variables with ${\mathbb E\,} X_1=0$, ${\mathbb E\,} X_1^2<\infty$
such that for each $n\ge  n_0$,
 \begin{equation}\label{Matskyavichyus}  \D_n\ge \p(n).
 \end{equation}
\end{theorem}

 The common distribution of these random variables is a mixture of symmetrized    distributions,  namely with   characteristic function
 $$ f(t) = \sum_{k=1}^\infty \lambda_k \exp\big\{\nu_k(\cos t-1)\big\},$$
 where
\begin{eqnarray*}
\begin{cases}\sum_{k=1}^\infty \lambda_k=1, &\quad \cr
\nu_k>0, \lambda_k\ge 0\ \hbox{for each $k=1,2,\ldots$}
\end{cases}
\end{eqnarray*}
and  $\s_k,\lambda_k$ are chosen with respect to   $\p(n)$, see Lemma 2 in \cite{Mat}.  These random variables are therefore extremal for the local limit theorem.

\bigskip\par
 Stronger integrability properties yield finer
remainder terms. This is made precise in the following statement.

 \begin{theorem} \label{r}    Let
$F$ denote the distribution function of
$X_1$.

{\rm  (i) (\cite{IBLIN}, Theorem 4.5.3)}      In order that  the property
\begin{equation} \label{alfa}
 \sup_{N=an+Dk}\Big|
 {\color{green} \frac{\s \sqrt n}{D} }{\mathbb P}\{S_n=N\}-{{\color{green}1}\over  \sqrt{ 2\pi}\s}e^{-
{(N-n\m )^2\over  2 n \s^2} }\Big| ={\mathcal O}\big(n^{-\alpha{\color{green}/2}} \big) ,
  \end{equation}
 {\color{green} where $0<\a<1$},
 it is necessary and sufficient that the following conditions be satisfied:
 \begin{eqnarray*} (1) \   D \ \hbox{is maximal}, \ \qq\qq
(2)  \  \  \int_{|x|\ge u} x^2 F(dx) = \mathcal O(u^{-\a})\quad \hbox{as $u\to \infty$.}
\end{eqnarray*}

   {\rm (ii) (\cite{P} Theorem 6 p.\,197)} If ${\mathbb E\,} |X_1|^3<\infty$, then \eqref{alfa} holds with $\a =1/2$.

\end{theorem}
Although evident from the assumption made, the link with classical central limit theorem is not apparent in the statements.

\begin{remark}[Bernoulli case]
  Theorem \ref{lltber} provides  a finer remainder term  than the one   directly derived   from  Theorem \ref{r}.   Better formulations can   be established by using a  more
precise, but  less handy
 comparison term.
For instance, there exists    an absolute constant $C$ such that for all integers $n\ge 2$,
 \begin{eqnarray} \label{better} \sup_{k}\Big|{\mathbb P}\{ B_n  = k \} - {1\over  \pi}  \int_\R  e^{i(2k-n)v -  n
({v^2
\over   2}+{ v^4 \over 12} ) }\, d v\Big|\le C\,  {
\log^{7/2} n
\over n^{ 5/2}}.
\end{eqnarray}
\end{remark}

\vskip 7 pt
   Galstyan \cite{Gal} characterized in the theorem below a stronger form of the local limit theorem.
\begin{theorem}\label{galstyan}Let $\{X_k, k\ge 1\}$ be a sequence of independent, identically distributed {\color{green}random variables} taking values in a common lattice
$\mathcal L(v_{
0},D )$  with $\E X_1=0$, $\s^2= \E X_1^2<\infty$. Assume the   span $D$ is maximal.\vskip 2pt
Then in order that
\begin{equation} \sum_{n=1}^\infty n^{-1+\frac{1}{2\d}} \sup_{k\in \Z} \Big|\frac{\s \sqrt n}{D} \P\{ S_n= nv_0+ {\color{green}k}D\} -\frac{1}{\sqrt{2\pi}} e^{-\frac{(n\t_0+
{\color{green}k}D)^2}{2\s{\color{green}2} \sqrt n}}\Big| <\infty, \qq   0\le \d <1,
\end{equation}
it is necessary and sufficient that
\begin{eqnarray*}  & & \E |X_1|^{2+\d}<\infty,\qq\qq\qq \ \ {\it if }\ 0<\d<1,
\cr &\cr  & & \E |X_1|^{2}\log(1+|X_1|)<\infty,\qq {\it if }\ \d=0.
\end{eqnarray*}
\end{theorem}
This Theorem is an extension  for densities and for lattice valued variables of an earlier  result proved  by Heyde in \cite{hey}.

\vskip 11 pt
 Consider now the case  of sums of i.i.d. stable random variables.

\begin{definition}\label{lltdefstable}Let $G(x)$ be the stable distribution function for which $M(x)=0$, $N(x) = -\frac{1}{x^\a}$, $\s^2= 0$ and $\g(\tau)= \a
\tau^{1-\a}/(1-\a)$ in the L\'evy--Khintchin formula.
We say that the sequence $ \{X_n , n\ge 1\}$ satisfies  a local limit theorem if, letting $S_n= X_1+ \ldots+ X_n$ for each $n$, we have,
\begin{equation*}
B_n{\mathbb P}\{ S_n= m\} - g\big( \frac{m}{B_n}\big) \ \to \ 0
\end{equation*}
when $n\to \infty$, uniformly in $m$, $-\infty <m<\infty$, where $g(x)= G'(x)$.
\end{definition}

The local limit theorem
  for i.i.d. stable random variables {\color{blue} is due to Ibragimov and Linnik} and  states as follows.}

{\color{blue}
\begin{theorem}[\cite{IBLIN}, Th. 4.2.1]\label{th:G4}
Let $ \{X_n , n\ge 1\}$  be a sequence of i.i.d.   random variables taking values in a common lattice $\mathcal L(v_{
0},D )$, let $S_n= X_1+ \ldots+ X_n$ for each $n$.  In order that  for some choice of constants $A_n$ and $B_n$
$$\lim_{n \to \infty}\sup_{m\in\Z}\Big|\frac{B_n}{D}\,{\mathbb P}\{S_n=v_0n+m D\}-g\Big( \frac{v_0n+m D-A_n}{B_n}\Big)\Big|=0, $$
where $g$ is the density of some stable distribution $G$ with exponent $0< \alpha \leq 2$,
 it is necessary and sufficient that
$$ {\rm (i)}\ \
 \frac{S_n-A_n}{B_n}\ \buildrel{\mathcal D}\over{\Rightarrow}  \ G   \ \
\hbox{as $n \to \infty,$} \qq\qq
    {\rm (ii)}\ \   \hbox{$D $ is maximal}.$$
\end{theorem}
  We give the proof, although not with full details. It also serves
  as an illustration of  the well-known {\it method of characteristic functions}, which is the common approach of many   proofs of the local limit theorems investigated by
  the Russian and Lithuanian Schools of Probability. See \cite{IBLIN}, \cite{Pe} and  the references therein. See also Gamkrelidze's proofs, in particular hints of  proofs of
  Theorem \ref{gam.suff.cond.sllt}  and    Mitalauskas' Theorem \ref{mitalauskas.th1}, see also Shukri \cite{Shu}. This also makes alternate approaches using various
  characteristics of a random variable, or based on other arguments,  of particular  interest.\vskip 8 pt \noi
{\tf Method of   characteristic functions.}
There is no loss (see \eqref{llt.transf.}) in assuming for the proof of Theorem \ref{th:G4} that $v_0=0$, $D=1$.
Let $f$ denote the characteristic of $X_1$. By the Fourier inversion formula (letting in what follows $z=(m-A_n)/B_n$),
\beq
 \P\big\{S_n= m  \big\}\,=\,\frac{1}{2\pi} \int_{-\pi}^\pi e^{-i zt  -itA_n/B_n} f(t/B_n)^n\dd t.
\eeq

Let also  $v$ denote the characteristic function of the $G$. We further have
\beq
 g(z) \,=\, \int_{-\infty}^{\infty} e^{-i zt  } v(t)\dd t.
\eeq
We can write
\beq\label{basic.bound}
\Big|B_n\P\big\{S_n= m  \big\}-g\Big(\frac{m-A_n}{B_n}\Big)\Big|
\,\le \, I_1+I_2+I_3+I_4,
\eeq
where
\ben I_1 &=&  \int_{|t|\le A} \big| e^{-itA_n/B_n}f(t/B_n)^n-v(t)\big|\dd t,
\cr I_2 &=& \int_{|t|>A}  \big|  v(t)\big|\dd t,
\cr I_3 &=& \int_{A\le |t| \le \e B_n}  \big| f(t/B_n)  \big|^n\dd t,
\cr I_4 &=& \int_{\e B_n\le |t| \le \pi B_n}   \big| f(t/B_n)  \big|^n\dd t,
\een
and $A $, $\e$ are constants to be determined. That {\color{green}approach for} controlling the left-hand side  of \eqref{basic.bound} in bounding it with the sum of the
previous four integrals is   frequently used. See   proof  of   Theorem \ref{mitalauskas.th1} for a more delicate case.

\vskip 5 pt Integral $I_1$ can be made small in view of  assumption (i), and $I_2$ is small for $A$  large. Now   that $I_3$ is small for $n$ large follows from the fact that
for any $\d<\a$, there exists a positive constant $c(\d)$, such that for $t$ in some neighbourhood of $0$, and all $n$
$$  \big| f(t/B_n)  \big|^n\, \le\, e^{-c(\d) |t|^\d}.$$
Thus
$$ \int_{A\le |t| \le \e B_n}  \big| f(t/B_n)  \big|^n\dd t\le  \int_{  |t| \ge A}  e^{-c(\d) |t|^\d}\dd t\quad \to \, 0,$$
as $A\to\infty$.
\vskip 2 pt Finally as $1$ is maximal, for some positive constant $c$
$$ |f(t)| \,\le\, e^{-c},$$
if $\e\le |t|\le \pi$. Then
$$ \int_{\e B_n\le |t| \le \pi B_n}   \big| f(t/B_n)  \big|^n\dd t\,\le \,2\pi B_n\,e^{-cn}\quad \to \, 0,$$
since $B_n=o\big(e^{cn}\big)$ (see e.g. \cite{IBLIN}, Th. 2.1.1).

\vskip 10 pt

  Theorem \ref{th:G4} admits  a strenghtened version, that is, a strong local limit theorem with convergence in variation (Definition \ref{slltv}). More precisely, we have
  the following

  \begin{theorem}[\cite{IBLIN},   Th. 4.2.2]\label{th:G42}
 Under the assumptions of Theorem \ref{th:G4}, we have
 \beq \label{sllt1}
 \lim_{n\to\infty}
 \sum_{m\in \Z} \Big| \,{\mathbb P}\{S_n=v_0n+m D\}-\frac{D}{B_n} g\Big( \frac{v_0n+m D-A_n}{B_n}\Big)\Big|=0.
 \eeq
\end{theorem}

\vskip 20 pt Local limit theorems for densities in various norms (sup-norm, $L^p$-norms)
 were intensively  studied. We refer to \cite{IBLIN}, Ch.\,4 and to \cite{P} Ch.\,7. In the simplest case where $X_j$, $j\ge 1$ are i.i.d. centered   random variables with
 $\E X_j^2=\s^2<\infty$, $\s^2>0$,  assuming that for all $n$ large enough the distribution function of $Z_n= (\sum_{j=1}^n X_j)/\s\sqrt n$, is absolutely continuous with
  density function $p_n(x)$,   the following result holds.
\begin{theorem}[\cite{IBLIN}, Th. 4.5.1]\label{th:G43}
 In order  that
\beq \label{lltd1}
\sup_x\big|p_n(x)-\phi(x)\big| \,=\, \mathcal O \big( n^{-\d/2}\big),
 \eeq
with $0<\d<1$, where $\phi(x)$ is the standard normal density function, it is necessary and sufficient that
\beq \label{lltd2}
\int_{|u|\ge x}y^2 F(\dd y)  \,=\, \mathcal O\big( x^{-\d}\big),
 \eeq
 where $F$ is the distribution function of $X_1$, and that there exists $N$ such that
 \beq \label{lltd3}
\sup_x p_N(x) <\infty.
 \eeq\end{theorem}

There are   variants or extensions of this result,   to   independent non-identically distributed random variables, see  Theorem \ref{maej.th},    or to weighted sums of
i.i.d. random variables  as we shall see below. Some refinements in terms of asymptotical H\"older continuity are obtained in Macht and Wolf \cite{MW}, where a smoothing
inequality is used.

\vskip 7 pt
\noi {\tf Local limit theorems for weighted sums of i.i.d. random variables.}
  The following problem is considered. Let $\{\xi_k, k\in \Z\}$ be a sequence of i.i.d. random variables. Let $\{a_{k,n}, k\in\Z, n\ge 1\}$ be a matrix. Put for each $n\ge
  1$,
\beq\label{wsiid} S_n=\sum_{k\in \Z} a_{k,n} \,\xi_k.
\eeq
These sums are studied in Shukri \cite{Shu}, when the distribution function of $\xi_k$ belongs to the domain of attraction of some stable law with index $0<\a\le 2$.

Let $F(x)$ (resp.\,$F_n(x)$), $p(x)$ (resp.  $p_n(x)$ and $f(t)$ (resp. $f_n(t)$) successively denote the distribution function, the {\color{green} density} distribution
and characteristic function of the i.i.d. random variables $\xi_k$ (resp. $S_n$).
\vskip 3 pt
Prerequisites on the representation of stable laws are necessary. As $F(x)$ belongs to the domain of attraction of some stable law with index $0<\a\le 2$,   we have the
classical representation (\cite{IBLIN},  Ch. II, Section 6)
 \begin{equation*} \label{}
{\rm (I)} \qq\qq\qq\qq \begin{cases} F(x)=\frac{c_1+o(1)}{|x|^\a}\, h(|x|),\quad  & x\to-\infty,\ c_1\ge 0,\cr
1-F(x)=\frac{c_2+o(1)}{x^\a}\, h(x),\quad  & x\to \infty,\ c_2\ge 0, \ c_1 +c_2 >0,
\end{cases}
\end{equation*}
where $h(x)$ is a function  varying slowly as $x\to \infty$.
The case $\a=2$, ${\rm Var}( \xi_k)<\infty$
 was studied in \cite{DS} and \cite{ShSa}. If $\a=2$, ${\rm Var}( \xi_k)=\infty$, then {\color{red} we assume that} $F$ has the following representation
 \begin{equation*} \label{}
{\rm (I')} \qq\qq\qq
 \chi(x) \, =\, 1-F(x) + F(-x) \,= \, \frac{h(x)}{x^2}, \qq \quad x\to \infty,\qq\qq\qq
\end{equation*}
where  $h(x)$   is a function   that varies slowly as $x\to \infty$.
 The characteristic function $f(t)$ is representable in a neighborhood of 0 in the form
  \begin{equation*} \label{}
 {\rm (II)}\qq\qq \qq \qq \quad f(t) \,=\, \exp\{ -c|t|^\a e^{{\rm sign}(t) i\p}\}, \qq \qq \quad  c>0, \ 0<\a\le 2,
\end{equation*}
where $|\p|\le \pi\a /2|$ if $0<\a<1$, $|\p|\le (2-\a)\pi/2$ if $1<\a \le 2$, and
$\tilde h(t)$ is a function  varying slowly as $t\to 0$.
\vskip 2 pt
In particular, if $F(x)$ and $f(t)$ satisfy (I) or (I') and (II), $c$ can be chosen so that,
 \begin{equation*} \label{}
 \tilde h(t) \,\sim\, \begin{cases}
  h(1/|t|), \qq \quad   &\ 0<\a\le 2,\cr
  \int_0^{1/|t|} u^2\chi(\dd u) \qq \quad   &\ \a=2.
 \end{cases}
\end{equation*}

Conditions under which the series \eqref{wsiid} converges {\color{green} almost surely}
are given in the following proposition.
\begin{proposition}[\cite{Shu},\,Prop.\,1]\label{wslltiid0} Let $F(x) $ belong to the domain of attraction of some stable law with index $0<\a\le 2$. Assume that the
sequence $\{a_k, k\in \Z\}$ satisfies the condition
\beq
\sum_{k\in \Z}|a_{k}|^\a\, \tilde h( a_k)\, <\, \infty.
\eeq
Then $\sum_{k\in \Z}a_{k} \xi_k$  converges almost surely.
\end{proposition}

\vskip 10 pt The local limit theorems for $S_n$ are obtained under the following three conditions, the first arising naturally {\color{green} from}
 the previous Proposition.

 \vskip 3 pt
{\emph{Condition}} ${\rm A_1}$: $ \sum_{k\in \Z} |a_{k,n}|^\a \tilde h (a_{k,n})  \, \to \, 1$ as $n\to \infty$, where $\tilde h$ is the slowly varying function from the
representation  {\rm (II)}.
\vskip 2 pt
{\emph{Condition}} ${\rm A_2}$: $ \sup_{k\in \Z} |a_{k,n}|\, \to \,0$ as $n\to \infty$.
\vskip 2 pt
{\emph{Condition}} ${\rm A_3}$: For some $p>0$, $\int_\R |f(t)|^p \dd t <\infty.$
\vskip 3 pt
\noi Conditions ${\rm A_1}$, ${\rm A_2}$ are natural in regard of  regularity conditions of   matrix summation methods.
\begin{theorem} \label{wslltiid1}
Let $F(x)$ belong to the domain of attration of a stable law $G(x)$ with index $0<\a\le 2$, $\a\neq 1$, and let the matrix coefficients $\{a_{k,n}\}$ satisfy   Condition
${\rm A_2}$, and for some $0\le P\le 1$, the conditions
\ben\label{wslltiid2}
  {\rm (1)} & & \sum_{k\,:\, a_{k,n}>0} |a_{k,n}|^\a \tilde h (a_{k,n}) \ \to \ P,\quad n\to \infty,
\cr  {\rm (2)} & & \sum_{k\,:\, a_{k,n}<0} |a_{k,n}|^\a \tilde h (a_{k,n}) \ \to \ 1-P,\quad n\to \infty,
\een
where $\tilde h$ is the slowly varying function from the representation  {\rm (II)}. Then,
\beq\label{wsiid3} F_n(x) \ \to \ L(x)\quad n\to \infty,
\eeq
where
\beq\label{wsiid4} L(x)= G\big(xP^{-1/\a}\big)*(1-G\big(-x(1-P^{-1/\a})\big). \eeq
\end{theorem}

\begin{theorem} \label{wslltiid5}
Let the conditions of Theorem \ref{wslltiid1} be fulfilled. If moreover Condition ${\rm A_3}$ is satisfied , then the density $p_n(x)$ exists  for all $n$ larger than a
certain value, and
\beq\label{wsiid3a} \sup_x \big|p_n(x) -l(x)\big|\ {\color{green}\ \to \ 0, }
\quad n\to \infty,
\eeq
 where $l(x)$ is the density of the stable distribution $L(x)$ from Theorem \ref{wslltiid1}.
\end{theorem}

\vskip 7 pt
\noi {\tf Local large deviations.}
Let  $  \{ S_n, n\ge 0\}$  be an integer-valued random walk
such that $\frac{S_n}{a_n}$,  converges in distribution to a stable law of index $\a\in (0,1)$ as $n\to \infty$, then Gnedenko's local limit theorem provides a useful
estimate for ${\mathbb P}\{S_n=r\}$ for values of $r$ such that $r/a_n$  is bounded. However, under appropriate conditions,  another estimate  is valid when $r/a_n\to
\infty$,   establishing a large deviation local limit theorem.
\vskip 3 pt
This type of estimates proceeds from  Nagaev's classical and elegant result \cite{Na}, which we recall.
Let $ S= \{ S_n, n\ge 0\}$ be a random walk   whose increments $X_i$ are independent copies of $X$, where ${\mathbb E\,} X= 0$, and
\begin{equation}\label{nag0}{\mathbb P}\{X>x\} \sim x^{-\a} L(x),  \qq \hbox{as $x\to  \infty$},
 \end{equation}
where $\a\in (1, \infty[$ and $L$ is slowly varying at infinity. Then
for any $\e >0$ and uniformly in $x\ge \e n$,
\begin{equation}\label{nag}{\mathbb P}\{S_n >x\} \sim n\, {\mathbb P}\{X> x\}, \qq \hbox{as $n\to  \infty$}.
 \end{equation}
 \vskip 10 pt

 {\color{green}  Doney showed that for lattice-valued random walks, there is an analogous local limit theorem.
  \begin{theorem}[\cite{D1},\,Th.\,1]\label{doney.th}
    Suppose   that $\{ S_n, n\ge 0\}$ is an integer-valued random walk with ${\mathbb E\,} X= \m$ finite and such that
\begin{equation}\label{don0}{\mathbb P}\{X=m\} \sim m^{-(1+\a)} L(m),  \qq \hbox{as $m\to  \infty$},
 \end{equation}
where $\a\in (1, \infty[$ and $L$ is slowly varying at infinity.    Then for any $\e >0$ and uniformly in $m\ge (\m+\e) n$,
\begin{equation}\label{don}{\mathbb P}\{S_n =m\} \sim n\, {\mathbb P}\{X=[m-n\m]\}, \qq \hbox{as $n\to \infty$.}
 \end{equation}
\end{theorem}
}
 \vskip 10 pt Tkachuk \cite{Tka}, {\color{green} see also Doney \cite{D},} showed  that if $S$ is in the domain of attraction of a stable law with index $0<\a<1$ and
 positivity parameter $0<\rho \le 1$, then uniformly  for $n$ such that $x/a_n\to \infty$, where $a_n$ is a norming sequence for $S$,
 \begin{equation}\label{tkachuk}    \P\{S_n>x\}\sim n\P\{X>x\}\qq \hbox{as $x\to \infty$},
 \end{equation}

\vskip 10 pt

{\color{green} Doney  proved in \cite[Th.\,A]{D}   a local version of \eqref{tkachuk}.} Let $F(x)={\mathbb P}\{X>x\}=1-F(x)$. Introduce the tail ratio $\tau$ defined by,
 \begin{equation}\label{tailratio}\tau (r) = \frac{F(-r)}{F(r)}=\frac{\sum_{-\infty}^{-r} p_k}{\sum_{r+1}^\infty p_k}.
 \end{equation}

\begin{theorem}
Assume that the following assumptions are satisfied. For some slowly varying function $L$,
\begin{eqnarray}
{\rm (i)} & & p_r \sim \a r^{-(\a+1)}L(r),\qq \ \, \hbox{as $r\to \infty$,}
\cr
{\rm (ii)} & & \lim_{r\to \infty}\tau (r) =\rho^{-1} -1, \quad\ \hbox{where  $0<\rho\le 1$.}
\end{eqnarray}
Then, uniformly in $n$ such that $r/a_n\to \infty$,
 \begin{equation}\label{moderate.llt}{\mathbb P}\{S_n =r\}\sim n\,{\mathbb P}\{X=r\}\qq \hbox{as $r\to \infty$}.
 \end{equation}
\end{theorem}

Results of  this kind can be used for establishing   local versions of generalized renewal theorems. Assume that $X$ takes non-negative values only. Thus $S$ is a discrete
renewal process. Consider the renewal mass function,
\begin{equation}\label{renewmassfunc} u_r=\sum_{n=0}^\infty {\mathbb P}\{S_n=r\}.
 \end{equation}

 Here we have  {\color{green}by using Theorem 8.7.3 in Bingham, Goldies and Teugels \cite{BGT},} that
 \begin{equation}\label{equiv1} F(n)=\sum_{k=n+1}^\infty p_k \sim n^{-\a} L(n), \qq \hbox{as $n\to \infty$},
 \end{equation}
 is equivalent to
 \begin{equation}\label{equiv1a}  \sum_{k=0}^n u_k \sim n^{\a} \frac{\Gamma(1-\a)}{L(n)\Gamma(1+\a)}, \qq \hbox{as $n\to \infty$},
 \end{equation}

\begin{theorem} Assume that $X$ takes non-negative values only, \eqref{equiv1} holds and further
 \begin{equation}\label{equiv3}  \sup_{n\ge 0}\, \frac{np_n}{F(n)}<\infty.
 \end{equation}
Then \begin{equation*}\lim_{n\to \infty}  n^{1-\a} L(n)u_n={\Gamma(1-\a)}\big/{\Gamma(1+\a)}.
 \end{equation*}
\end{theorem}
The two preceding theorems are respectively Theorems A and B in Doney \cite{D}.
}

\vskip 15 pt
{\color{blue} \noi{\tf Diophantine measures and local limit theorem.}
  Breuillard studied in \cite{Br} the local limit theorem in $\R $, also in $\R^d$, in the i.i.d. case for Diophantine measures. Let $\{X_i, i\ge 1\}$ be a sequence of i.i.d.
  $\R $-valued, centered random variables of common law $\m$. In order to study the asymptotic behavior as $n$ tends to infinity of the means
  \beq \E f(S_n)= \int  f \dd \m^n
  \eeq
for $f$ defined on $\R $, the author introduces the notion of Diophantine measure, {\color{green} in analogy with the {\it   type}  of a real number. Let $\eta$ be a positive
real number or infinity. An irrational number $\a$ is said to be of type {\it   type $\eta$} if $\eta$ is the supremum of all $\g$ for which
$$\liminf_{q\to \infty}q^\g \{q\a\}=0 ,$$
 where  $q$ runs through the positive integers. See Kuipers and Niederreiter \cite {K.N}, p.\,121.}  So that if $\eta$ is finite and $\eta'>\eta$, then for some $C>0$
$$  \{q\a\}\ge \frac{C}{q^{ \eta' }},$$
 for any $q$ large enough. A measure $\m$ on $\R $ is $l$-Diophantine, for some real $l\ge 0$, if there exists  $C>0$ such that for any $x\in \R$ large enough in absolute
 value,
\beq\label{mu.l.dioph}
\inf_{y\in \R} \int \{xa+y\}^2 \dd \m(a) \ge \frac{C}{|x|^l}.
\eeq
Here $\{b\}$ denotes the distance from $b$ to the nearest integer. One says  that $\m$ is Diophantine if $\m$ is $l$-Diophantine, for some real $l\ge 0$.

\vskip 5 pt

This notion shares strong similarities with Mukhin's {\color{green}structural characteristic} $ \mathcal D(X,x)$\footnote{\ {\color{blue}not quoted in \cite{Br}}}  defined in
\eqref{ssprime2}, since the quantity to be estimated, {\color{green}namely the right hand side of \eqref{mu.l.dioph},} is the same. Indeed,
 \begin{eqnarray*}  \mathcal D(X,x)=\inf_{v\in \R}{\mathbb E\,} \{ (X-v)x\}^2=\inf_{y\in \R}\int  \{ xa+y\}^2 \dd \m (a).
 \end{eqnarray*}
limit theorems with remainder terms.
Note that the question considered was also studied in 1964 by Shepp in a  seminal paper \cite{Sh64}\footnote{\ {\color{blue}not quoted in \cite{Br}}}, see also Stone
\cite{Sto}, \cite{Sto1}.
\vskip 3pt

However   condition  \eqref{mu.l.dioph}  is specific, to our knowledge, to Breuillard's work, and is not present in Mukhin's papers. That condition seems not easy to check in
general. When $\m$ is a Dirac measure at point $a$, the $l$-Diophantinity of $\m$ coincides with the $l$-type of irregularity of $a$, expressing that $a$ is $l$-badly
approximable, see \cite{K.N}, Def. 3.4. This clearly justifies the terminology introduced.
\vskip 2pt
Condition  \eqref{mu.l.dioph}   translates to the Fourier transform of $\m$. Let $l\ge 0$. Then $\m$ is $l$-Diophantine if and only if there exists a real $C>0$ such that for
any $x$ large enough in absolute value,
\beq\label{fourier.cond.BR}
|\widehat \m(x)|\,\le\, 1-\frac{C}{|x|^l}.
\eeq
Other equivalences are formulated in term of the symmetrized measure $\m*\m^{-1}$. See \cite{Br}, Proposition 3.1. The proof is short and based on the elementary inequality
($c_1$, $c_2$ being positive constants)
$$  c_1\{x\}^2\le 1-\cos (2\pi x)\le c_2 \{x\}^2.$$
\vskip 3pt

The author generalizes to Diophantine measures some results of Feller \cite{F}, XVI Th. 2 (Edgeworth expansion  of densities under the assumption that $\widehat \m$ is
integrable). For instance (\cite{Br}, Th.\,3.2),
\begin{theorem}\label{th3.2.BR} Let $r$ be some non-negative integer, and a real $l\ge 0$. Assume that $\m$ is centered and has moment of order $r+2$ finite. If moreover $\m$
is $l$-Diophantine, then for any function $f$ of class $C^k$ with $k>l(r+1)/2+1$ and such that
\beq C^k(f):= \max_{0\le j \le k} \big\| f^{(j)}\|_1<\infty,
\eeq
one has the following asymptotic Edgeworth expansion
\beq\E f(S_n) \,=\, \sum_{p=0}^r \frac{1}{n^{p/2}} \int_\R f(x\sqrt n) P_p(x) g(x) \dd x + C^k(f)\,.\, o\Big(\frac{1}{n^{(r+1)/2}}\Big)
\eeq
where the constants hidden in the symbol $o$ depend on $r$, $\m$ only, $g$ is the Gaussian density associated to $\m$ through the central limit theorem, and $P_p$ are
polynomials of degree $\le 3p$, whose coefficients only depend of the moments of $\m$ of order less or equal to $p+2$. \end{theorem}

\vskip 3 pt
We say that $\m$ is aperiodic if
\beq\label{aperiodic} |\widehat \m(x)|<1, \qq \quad (\forall x\in \R\backslash \{0\})
\eeq
The author also examines for $\m$ centered, aperiodic on $\R$, with variance $\s_2$ and admitting a moment of order $r\ge 2$, when the limit (where $I_s=[-s,s]$)
\beq \label{} \lim_{n\to \infty} \frac{\m^n(I_s+x)}{\nu^n(I_s+x)}=1
\eeq
  is uniform in $x$ and $s$, {\color{green}$\nu$ being the standard Gaussian distribution}. For instance the limit is uniform when $|x|+s \le \sqrt{c\s_2n\log n}$. See
  \cite{Br}, Th.\,4.2.
  }

{\color{blue}

\vskip 15 pt
\noi {\tf Local limit theorems with arithmetical   conditions.}
The first approaches use arithmetical conditions of type: For all $q\ge 2$
\begin{equation}
  \max_{0\le r\le q-1}{\mathbb P}\{X_j\equiv r \ {\rm (mod\, {\it q})}\}\le 1-\a_j,
  \end{equation}
for all $j$, where $\a_j$ is some specific sequence of reals decreasing to $0$. Conditions of this sort already appeared in Mitalauskas 1960's paper \cite{Mit1}, and 1964's
paper by Raudelyunas \cite{Rau}, later in Fomin's paper \cite{Fo}. Roughly speaking, one requires the random variables to do not overly much concentrate in a particular
residue class $r$ (mod $q$)  of $\Z$.

 \smallskip\par
A condition of this kind is used in the following multidimensional result by Fomin \cite{Fo}, which however seems not of an easy  application.

\begin{theorem}\label{fomin.th} Let $X^{(n)}_\ell= (X^{(n)}_{\ell1},\ldots, X^{(n)}_{\ell d})$, $\ell=1,\ldots, n$, $n=1,2,\ldots$ be an array of independent i.i.d.
$\Z^d$-valued, centered random vectors with  components having third moments, and assume  that   the covariance matrix $R^{(n)}$ of $X^{(n)}$ is positive definite.   Let
$\s^2_{ni}= {\mathbb E\,} (X^{(n)}_{ni})^2$,  $B_{ni}^2=n\s_{ni}^2$ and $\b^2_{ni}={\mathbb E\,} |X^{(n)}_{ni}|^3$,  $i=1,\ldots, d$. Let also $S_n=\sum_{k=1}^n
X^{(n)}_k=(S_{n1}, \ldots, S_{nd})$, $\overline{S}_n=( \frac{S_{n1}}{B_{n1}}, \ldots,  \frac{S_{nd}}{B_{nd}})$. Assume that
 $$ \b_{ni}= o\big( \s_{ni}^3\D_n^2\sqrt n\big), \qq i=1, \ldots, d.$$
where ${\color{green}\D_n}=\det(R^{(n)})$.
  Further assume that for each $q\ge 2$, each $\overline{a} =(a_1, \ldots, a_d)\in \Z^d$, with $\gcd(a_1, \ldots, a_d, q)=1$,
\begin{equation}\label{fomin.cond}
  \max_{0\le r\le q-1}{\mathbb P}\{\overline{a}.X^{(n)}_1\equiv r \ {\rm (mod\, {\it q})}\}\le 1-\a_n,
  \end{equation}
   where
  $$ \a_n= K \max_{1\le i\le d}\max\Big(\frac{\b_{ni}^2}{\s^4_{ni}{\color{green}\D_n^2}n}, \frac{\s_{ni}}{\sqrt n}\Big)\log n, $$
 and   $K$ is some constant.  Then,
\begin{equation}\label{fomin.llt}
B_{n1}\ldots B_{nd} \, {\mathbb P}\big\{\overline{S}_n= \overline{z} \big\}-\phi_{R^{(n)}}(\overline{z})\ \to \ 0 ,
  \end{equation}
where  $ \phi_{R^{(n)}}(\overline{z})$ is the Gaussian density with covariance matrix $R^{(n)}$.
  \end{theorem}

The following theorem provides a fine local limit theorem for triangular arrays of i.i.d. sums. It was proved by Freiman, Moskvin and Yudin \cite{FMY}.

\begin{theorem}\label{freiman.moskvin.yudin.th}
For each $n=1,2,\ldots$, let $\tilde X_n= (X_{1n}, \ldots, X_{nn})$  be composed with i.i.d. $\Z$-valued  random variables, having third moments, and let ${\mathbb E\,}
X_{1n}=a_n$, $\E(X_{1n}-a_n)^2= \s_n^2$.
Assume that
\begin{eqnarray*} {\rm (1)}& & \hbox{$\tilde X_n$ satisfies the central limit theorem as $n\to \infty$,}
\cr {\rm (2)}& & \hbox{$\s_n^2= \mathcal O(n^\rho)$ where $\rho < \frac{\log (2+c)}{\log 2} -1$ and $c<1$.}
\cr {\rm (3)}& & \hbox{For $q=2, 3, \ldots$ and $\o>0$ sufficiently large,}
\cr & &\qq \max_{1\le r \le q} {\mathbb P}\big\{X_{1n}\equiv r \ {\rm (mod}\, q)\big\}< 1-\o \max\Big(\frac{\rho_n^2}{n\s_n^4},\frac{1}{n^{1-\m}}\Big)\, \log n,
\cr & & \hbox{where $\rho_n={\mathbb E\,} |X_{1n}-a_n|^3$ and $\m=\frac{(1+\rho)\log 2}{\log (2+c)}$.} \end{eqnarray*}
Then the local limit theorem for the sequence $(X_{1n}, \ldots, X_{nn})$, $n\ge 1$ holds.
\end{theorem}

In both statements, the third moment condition is restrictive. However, in some important applications, the random variables are bounded and this condition is automatically
satisfied, see   \cite{FMY}.
}

{\color{green}\subsection{The independent case}}  \label{1.4}

 {\color{blue}
 Let $X=\{X_i , i\ge 1\}$ be independent,  square integrable
random variables taking values in a common lattice $\mathcal L(v_{
0},D )$.
 Let
$S_n=X_1+\ldots +X_n$, $n\ge 1$.  Let also    $a_n= {\mathbb E\,} S_n$, $\s_n^2={\rm Var}(S_n)\to \infty$.

 \vskip 5 pt
 According to  Definition \ref{defllt},   the sequence $  X$ satisfies a local limit theorem if
 \begin{equation}\label{def.llt.indep}    \sup_{N=v_0n+Dk }\Big|\s_n {\mathbb P}\{S_n=N\}-{D\over  \sqrt{ 2\pi } }e^{-
{(N-a_n)^2\over  2 s_n^2} }\Big| = o(1), \qq \quad n\to\infty.
\end{equation}

}

{\color{blue}Let us fi{\color{green}r}st consider  the following  example   due to Prohorov \cite{Pr}.
\begin{example} Assume that $X_1$ takes the values 1 and 0 with probabilities $p$ and $q$, respectively, and let $X_j$, $j\ge 2$  be equal to $2$ with probability $\frac12$
and $0$ with the same probability. Then the  sequence $X$  satisfies the local limit theorem iff $p = q = \frac12$ and the local limit theorem does not hold for the sequence
$ X'= \{ X_j , j\ge 2\}$.
\end{example} Thus in this example, the fulfilment of the local limit theorem  depends on the behavior of the first members of $X$. Hence  it is reasonable to introduce the
following definition.
  \begin{definition}[\cite{Pr}] \label{deflltsf}Let $ X= \{ X_j , j\ge 1\}$ be a
sequence of independent random variables taking only integral values. A local limit theorem   in  the {\it strong form}
(or {\it in a strengthened form})
 is said to be
applicable to $  X$, if a local limit theorem in the usual form (Definition \ref{defllt}) is applicable to any subsequence extracted from $  X$, which differs from
$  X$ only in a finite number of members.
  \end{definition}
{\color{blue}
\begin{remark}The   terminology \lq  in  the strong form\rq, also \lq in strong sense\rq, is used by various authors, for instance in Gamkrelidze \cite{Gam3}, \cite{Gam},
Mitalauskas \cite{Mit}, Mitalauskas-Stepanauskas  \cite{MiSt}, Mukhin \cite{Mu}, \ldots  whereas the other is used in Petrov \cite{P}, Prohorov \cite{Pr}, Rozanov \cite{Ro},
notably. Different equivalent forms of the above definition are considered in   \cite{Gam3}.
\end{remark}}
\vskip 5 pt This definition can be made a bit more convenient (\cite{Gam3}).
Let
\ben  & &S_{k,n}=\xi_{k+1}+ \ldots +\xi_{k+n},\qq A_{k,n}=\E S_{k,n}, \qq B^2_{k,n}={\rm Var} (S_{k,n}).
\een
The local limit theorem   in  the  strong form  holds if and only if
\beq\label{lltsf.ref}
\P\Big\{ \sum_{j=k+1}^n X_j=m\Big\}= {D\over B_{k,n} \sqrt{ 2\pi } }\ e^{-
{(m-A_{k,n})^2\over  2 B_{k,n}^2} }+o\Big({1\over B_{k,n}  }\Big),
\eeq
uniformly in $m$ and every finite $k$, $k=0,1,2, \ldots$, as $n\to \infty$ and $B_{k,n}\to \infty$.

\vskip 7 pt

\noi {\tf Prohorov's Theorem.} Gnedenko's Theorem \ref{gnedenko} was generalized by Prohorov   in 1954, to sequences of   independent uniformly bounded random variables.

\begin{theorem}[\cite{Pr}] \label{Prohorov} Let $X= \{ X_j , j\ge 1\}$ be independent, uniformly bounded integral-valued random variables with partial sums
$S_n= X_1+\ldots +X_n$. Let $A_n= {\mathbb E\,} S_n$, $B_n^2= {\rm Var}(S_n)\to \infty$, and suppose that
$${\mathbb P}\{X_n=0\}= \max_{j\ge 1}\, {\mathbb P}\{X_n=j\}.$$
Then $  X$ satisfies a local limit theorem in the strong form if and only if
1 is the highest common factor of the values of $j$ with $\sum_{n=1}^\infty {\mathbb P}\{X_n=j\}=\infty$.\end{theorem}
\vskip 10 pt

Local limit theorems for asymptotically stable lattice distribution functions have been obtained by   Rva\v ceva \cite{Rauk} for the multi-dimensional case, generalizing
those of Gnedenko for the one-dimensional case. Stone \cite{Sto} obtained results of this sort for non-lattice distribution functions.
}
\vskip 7 pt

{\color{blue} Maejima  investigated in   \cite{maej2}  (continuing   \cite{maej1}),  the question of estimating  the remainder term in the local limit theorem for independent
random variables.  He notably deduced from a more general result he proved in this paper, an extension to the non-identically distributed case of  a result due to
Ibragimov-Linnik.

\vskip 3 pt
 Let $  \{X_k, k\ge 1\}$   be a sequence of independent
random variables with $\E X_k=0$, $\E X_k^2=\s_k^2<\infty$, ($\s_k\ge 0$ and not all null) and with distribution $F_k(x)$. Suppose for instance $\s_1>0$, and let $s_n^2=
\sum_{k=1}^n\s_k^2$, $Z_{ n}=s_n^{-1} \sum_{k=1}^n X_k$, $f_k(t)= \E e^{it X_k}$, $\bar f_n(t)= \E e^{itZ_n}$, $R_k(x) = \int_{|u|>x}u^2\,F_k(\dd u)$, $Q_k(x) =
\big|\int_{|u|\le x}u^2\,F_k(\dd u)\big|$. Further let $\bar p_n(x)$  be the density function of $Z_n$ and $\phi(x)$ be the standard normal density function.
 \vskip 3 pt

Consider the following two simple classes of functions:
\begin{align*}G&\,=\, \big\{g:\R\to \R: \hbox{$g$ is even, and on $]0,\infty)$, $g(x)>0$, $\frac{x}{g(x)}$ is non-decreasing} \big\}
\cr G_0&\,=\, \big\{g\in G: \hbox{  $\frac{x^\a}{g(x)} $ is non-decreasing on $]0,\infty)$  for some $0<\a<1$}\big\}.
\end{align*}
Set for $g\in G$,
$$\l_k (g)=\sup_{x>0} g(x) R_k(x), \quad \m_k (g)=\sup_{x>0} \frac{g(x)}{x} Q_k(x), \qq \rho_k (g)=\l_k (g)+\m_k (g).$$
Note that if $\E X_k^2\, g(X_k)<\infty$, then of course $\rho_k (g)<\infty$.
\begin{theorem}\label{maej.th} Let $g\in  G$. Suppose $\rho_k(g)<\infty$ for $1\le k\le n$. Further assume that
\begin{eqnarray*}{\rm (a)} & & \hbox{$s_n^2< Kn$, for some positive constant $K$},
\cr {\rm (b)} & & \hbox{$\sup_{x} p_k(x)<M$, for some  $M>0$ independent of $k$.}
\end{eqnarray*}
Then we have
\begin{eqnarray*}
\sup_x \big|\bar p_n(x) -\phi (x) \big|\,\le \, C \, \frac {\sum_{k=1}^n\rho_k (g)}{s_n^2g(s_n)} .\end{eqnarray*}
In particular if $\E X_k^2\, g(X_k)<\infty$ for $1\le k\le n$, then
\begin{eqnarray*}
\sup_x \big|\bar p_n(x) -\phi (x) \big|\,\le \, C \, \frac {\sum_{k=1}^n\E X_k^2\, g(X_k)}{s_n^2g(s_n)} .\end{eqnarray*}
\end{theorem}
 Applying  Theorem \ref{maej.th} with $g(x)=x^\d$, gives  the following extension of the local limit theorem by Ibragimov-Linnik (Theorem  \ref{th:G43}) to the case of
 non-identically distributed random variables.
\begin{corollary}\label{maej.cor}  Assume in addition to conditions {\rm (a)} and {\rm (b)}, that $\liminf_{n\to \infty} s_n^2/n>0$. Let $0<\d\le 1$. Then in order that
\begin{eqnarray*}
\sup_x \big|\bar p_n(x) -\phi (x) \big|\,= \,  \mathcal O\big( n^{-\d/2}\big),\end{eqnarray*}
it is sufficient for $0<{\color{green}\d}<1$, that
\beq \label{maej.cor.Rk} \frac{1}{n} \sum_{k=1}^n \sup_{x>0}R_k(x) =\mathcal O(1), \quad {\text as}\qq n\to \infty,
\eeq
 and for $\d=1$, that \eqref{maej.cor.Rk} with $\d=1$ and
\beq \label{maej.cor.Qk} \frac{1}{n} \sum_{k=1}^n \sup_{x>0}Q_k(x) =\mathcal O(1), \quad {\text as}\qq n\to \infty,
\eeq
hold.\end{corollary}

\vskip 3 pt
One key argument used in the course of the proof of Theorem \ref{maej.th} is Survila's inequality:

\begin{proposition}[\cite{Sur}] \label{survila} Under conditions {\rm (a)} and {\rm (b)}, we have for $1\le k_1\neq k_2\le n$,
\begin{eqnarray*}
|\bar f_n(t)|\ \le \
\begin{cases}\big|f_{k_1}(t/s_n)f_{k_2}(t/s_n)\big|\,\exp\{-cn\},\quad &{\text for} \ |t|\ge \frac{\pi}{\sqrt 2}\sqrt n,\cr
 \exp\{-ct^2\},\quad &{\text for} \ |t|< \frac{\pi}{\sqrt 2}\sqrt n.
 \end{cases}
\end{eqnarray*}
\end{proposition}
\vskip 25 pt
We consider also the strong strong local limit theorem with convergence in variation (\cite{IBLIN}, \cite{Gam3}).
Let $X=\{X_n, n\ge 1\}$   be a sequence of independent, square integrable, integer-valued
random variables, and let $S_{k,n}=X_{k+1}+ \ldots +X_{k+n}$ for each $n$.

 \begin{definition}\label{slltv}
   The sequence of partial sums $\{S_n, n\ge 1\}$ is said to satisfy the strong local limit theorem with convergence in variation if there are constants $A_{k,n}$ and
   $B_{k,n}$  such that for each $k=0,1,2,\ldots$,
 \beq \label{sllt1.d}
 \sum_{m\in \Z} \Big|\P\{ S_{k,n}= m\}-\frac{1}{B_{k,n}\sqrt{2\pi}}\exp\Big\{-\frac{(m-A_{k,n})^2}{2B^2_{k,n}}\Big\} \Big|\quad \to \ 0,
 \eeq
 as $n\to\infty$ and $B_{k,n}\to \infty$.

 \end{definition}

  A necessary condition for the sequence $X$ to satisfy the strong local limit theorem with convergence in variation, is given in    Gamkrelidze \cite{Gam3}. More precisely,
\begin{proposition}\label{gam.nec.cond.slltv} Let the sequence $X$ satisfy the strong local limit theorem with convergence in variation. Then for every finite $k$, there
exist a positive integer $n_0 = n_0(k)$ and a positive number $\l = \l(k)$ such that
\beq  \d_{S_{k,n_0}} =\sum_{m\in\Z}\big|\P\big\{S_{k,n_0}=m\big\}-\P\big\{S_{k,n_0}=m-1\big\}\big|<\l.
\eeq
\end{proposition}

The smoothness characteristic $\d_{S_{k,n_0}} $  is defined for general integral-valued random variables in \eqref{delta} and it
  was introduced and investigated by Gamkrelidze in  \cite{G2} notably,  see   Proposition \ref{deltaX.prop} for some properties.
\vskip 4 pt
{\color{blue}  In \cite{Gam3} Gamkrelidze  made the following instructive remark: \lq\lq {\it  Until now very general conditions have been obtained under which a local limit
theorem follows from an integral limit theorem in some form. However, as a rule, these conditions are hard to verify in the case of summands that are not identically
distributed. In \cite{G1},
sufficient conditions for the local limit theorem are expressed in very simple terms, but they do not have great generality. The purpose is to give sufficiently general
conditions for the local limit theorem to hold by means of simple and intuitive characteristics of the summands.
}

 \vskip 4 pt In the same paper \cite{Gam3}, he {\color{green}further} showed that his necessary condition is sharp, and proved the following result improving upon Theorem
 \ref{gamkrelidze}.
\begin{theorem}\label{gam.suff.cond.sllt}
Assume that the sequence $X$ satisfies the following conditions:
\vskip 2 pt {\rm (1)} There exists a positive integer $l$, a positive $\l$ such that $\d_{S_{k,l}} \le \l<\sqrt 2$  for some $k$,
\vskip 2 pt {\rm (2)} For every finite $k$, $\frac{ S_{k,n}-A_{k,n} }{B_{k,n}}$ converges in distribution to $\mathcal N(0,1)$ as $n\to \infty$.
\vskip 2 pt {\rm (3)}
For every finite $k$, $B_{k,n}$ tends to infinity not rapidly than the power function, i.e. $B_{k,n}=\mathcal O( n^t)$  as $n\to \infty$.
\vskip 2 pt
Then $X$ verifies the local limit theorem   in  the  strong form, namely
\beq \label{sllt1a}
 \sup_{m\in \Z} \Big|B_{k,n}\P\{ S_{k,n}= m\}-\frac{1}{\sqrt{2\pi}}\exp\Big\{-\frac{(m-A_{k,n})^2}{2B^2_{k,n}}\Big\} \Big|\quad \to \ 0,
 \eeq
 as $n\to\infty$.
 \end{theorem}
We provide some hints of proof.  The proof    follows the usual scheme of proof of local limits theorems using characteristic functions.

 \vskip 8 pt
 Let $f(t, S_{k,n})$ denote the characteristic function of $ S_{k,n}$. By the Fourier inversion formula
\begin{equation}
2\pi \Big( B_{k,n}\P\{ S_{k,n}= m\}-\frac{1}{\sqrt{2\pi}}e^{-\frac{z^2}{2}}\Big)=I_1+I_2+{\color{green}I_3+I_4},
\end{equation}
where $z=z_{n,m}(m-A_{k,n})/B_{k,n}$ and
\begin{eqnarray}
I_1&=& \int_{|t|\le A} e^{-izt}\Big[f\Big( \frac{t}{B_{k,n}},  S_{k,n} \Big) -e^{-\frac{t^2}{2}}\Big]\dd t\cr
I_2&=& \int_{|t|\ge A}e^{izt-\frac{t^2}{2}}\dd t\cr
I_3 &=& \int_{A\le |t|\le \e B_{k,n}} e^{-izt} f\Big( \frac{t}{B_{k,n}},  S_{k,n} \Big)  \dd t
\cr I_4 &=& \int_{\e B_{k,n}\le |t|\le \pi B_{k,n}} e^{-izt} f\Big( \frac{t}{B_{k,n}},  S_{k,n} \Big)  \dd t. \end{eqnarray}
Note that $I_1, I_2$ can be made small, the first because of condition (2) and the second by choosing $A$ {\color{green}large}. Now $I_3$ is small for $ |t|\le \e B_{k,n}$,
since because of Prohorov's bound in \cite{Pr}, $|f (  {t}/{B_{k,n}},  S_{k,n}  )|\le e^{-t^2/3}$. Finally $| I_4|$ is made small by using Lemma \ref{L2G1}.
The proof is quite similar to the one given in \eqref{I3.G}.

}\vskip 7 pt

\noi{\tf A general necessary condition.} {\color{blue} Before stating Rozanov's   Theorem, let us first introduce a useful definition.
 \begin{definition}\label{aud}
 Let $\{X_n, n\ge 1\}$   be a sequence of independent, integer-valued
random variables, and {\color{green}let} $S_n=\sum_{k=1}^n X_k$, for each $n$. The sequence of partial sums $\{S_n, n\ge 1\}$ is said to be asymptotically uniformly
distributed with respect
to lattices of span $h$, in short a.u.d.($h$), if for $m = 0,1,\ldots,h-1$, we
have
 \beq \label{aud1}\lim_{n\to \infty}\ \P\{S_n \equiv m \,{\rm (mod}\,h)\}=\frac1h.
 \eeq
 We say that the sequence $\{S_n, n\ge 1\}$ is  asymptotically uniformly distributed, in short a.u.d., if \eqref{aud1}  holds true for any $h\ge 2$ and $m = 0,1,\ldots,h-1$.
 \end{definition}
\vskip 5 pt
Dvoretzky and Wolfowitz \cite{DW} proved the following characterization. Let $\{X_k,k\ge 1\}$ be a sequence of independent random variables taking on only the values
$$ 0, 1,\ldots, h-1.$$
Let   $S_n=\sum_{k=1}^n X_k$, for each $n$. In order that the partial sums $\{S_n, n\ge 1\}$  be a.u.d.($h$), it is necessary and sufficient that
\beq \label{aud.dw.ns} \prod_{k=1}^\infty\bigg( \sum_{m=0}^{h-1}\P\{X_k=m\}\,e^{\frac{2i\pi }{h}rm}\bigg) \,=\, 0, \qq \quad (r=1,\ldots, h-1).
 \eeq
Equivalently,
\beq \label{aud.dw.ns.e} \prod_{k=1}^\infty{\color{red} \E}\big(  e^{\frac{2i\pi }{h}rX_k}\big) \,=\, 0, \qq \quad (r=1,\ldots, h-1).
 \eeq

 The condition \beq \label{aud.dw.S} \sum_{k=1}^\infty\bigg( \min_{m=0}^{h-1}\P\{X_k=m\} \bigg) \,=\, \infty ,
 \eeq
is sufficient for $\{S_n, n\ge 1\}$  be a.u.d.($h$).}

\vskip 15 pt   Rozanov's   necessary condition   was found in  1957 and states as follows.

\begin{theorem}[\cite{Ro},\,Th.\,I]\label{rozanov.I} Let $  X= \{ X_j , j\ge 1\}$ be a sequence of independent, square integrable random variables taking only integral
values. Let $b_k^2= {\rm Var}(X_k)$, $B_n^2 =b_1^2+\ldots+ b_n^2$. Assume that
\begin{equation}\label{Ro.A}B_n \to \infty\qq\qq {\it as} \ n\to \infty.
\end{equation}
The following condition is necessary for the applicability of a local limit theorem  in the strong form to the sequence $X${\color{green},}
\begin{equation}\label{rozn}\prod_{k=1}^\infty
\big[ \max_{0\le m< h} {\mathbb P}\big\{X_k\equiv m\, {\rm (mod {\it \, h })} \big\}\big]=0\qq {\it for\ any\ }  h\ge 2 .
\end{equation}
\end{theorem}

This is naturally an important result in the theory of the local limit theorem. The proof is based on the following lemma, {{\color{blue}which provides a necessary condition
for the validity of the local limit theorem}.
{\color{blue}
 \begin{lemma}  \label{roz.cond.equi}Suppose that the   {\color{green}local limit theorem}
   is applicable to the sequence $X$. Then the partial sums
$S_n$ are {a.u.d.}
\end{lemma}
 {\color{green}This is Lemma 1 in \cite{Ro},
 we provide  a detailed proof.}
\footnote{\color{blue}\,  In Petrov \cite{P}, Lemma 1,\,p.\,194, also in   Rozanov's \cite{Ro} Lemma 1,\,p.\,261,  Lemma \ref{roz.cond.equi} is   stated  under the
assumption that a  local limit theorem  in the strong form holds, which is not necessary.}}
\begin{proof}
By assumption,
\beq\label{lltsf.ref.a}
\P\big\{ S_n=m\big\}= {D\over B_n \sqrt{ 2\pi } }\ e^{-
{(m-M{\color{ green}_n})^2\over  2 B_n^2} }+o\Big({1\over B_n }\Big),
\eeq
uniformly in $m$, as $n\to \infty$.
  Let {\color{green}$\e>0$  be fixed} and  note that  ${\mathbb P}\{   |S_n -M_n |> B_n/\sqrt \e\}    \le      \e$.
Then
\begin{eqnarray*}
{\mathbb P}\big\{ S_n\equiv\, m\ \hbox{\rm (mod $h$)}\big\}&=&  \sum_{|k-M_n|\le B_n/\sqrt \e
\atop k\equiv m\, (h)} {\mathbb P}\big\{S_n=k\}+ \mathcal
O(\e)
 \cr &=& {D\over  \sqrt{ 2\pi }B_n }    \sum_{|k-M_n|\le B_n/\sqrt \e
\atop k \equiv m\, (h)} \Big( e^{-
{(k-M_n)^2\over  2 B^2_n} }+o\Big({1\over B_n }\Big)\Big)+ \,\mathcal O(\e)
\cr &=& {D\over  \sqrt{ 2\pi }B_n }    \sum_{|k-M_n|\le B_n/\sqrt \e
\atop k \equiv m\, (h)}   e^{-
{(k-M_n)^2\over  2 B^2_n} } +{D\over  \sqrt{ 2\pi   \e} }  \,o(1)+ \,\mathcal O(\e),
\end{eqnarray*}
where the above $o(1)$ term tends to 0 as $n$ tends to infinity.
{\color{green}Now as \begin{eqnarray*} \frac{1}{B_n}\,\sum_{|k-M_n|> B_n/\sqrt \e } e^{-
{(k-M_n)^2\over  2 B^2_n} }&= &    o(\e)   ,
\end{eqnarray*}
 we get
\begin{eqnarray*}
{\mathbb P}\big\{ S_n\equiv\, m\ \hbox{\rm (mod $h$)}\big\}
 &=& {D\over  \sqrt{ 2\pi }B_n }    \sum_{ k \equiv m\, (h)}  e^{-
{(k-M_n)^2\over  2 B^2_n} }+   \,\mathcal O(\e)
.
\end{eqnarray*}
 Noticing that
 \begin{equation}\label{est2}    {  1\over \sqrt{2 \pi}B_n}  \sum_{ k \equiv m\, (h)}   e^{-
{(k-M_n)^2\over  2 B^2_n} } =\frac1h + \mathcal O\Big( {1\over    B_n   }\Big) .
\end{equation}
we obtain\begin{eqnarray*}
{\mathbb P}\big\{ S_n\equiv\, m\ \hbox{\rm (mod $h$)}\big\}
  &=&  \frac{D}{ h} +\mathcal O(\e),
\end{eqnarray*}
as $n$ tends to infinity. But $\e$ could be chosen arbitrarily small. And so the claimed result is proved.}
\end{proof}

   \vskip 3 pt

We now deduce from  Lemma \ref{roz.cond.equi},  that if  the   local limit theorem    in the strong form  is applicable to the sequence $X$, then
\begin{equation} \label{roz.cond}
\sum_{k=1}^\infty \  {\mathbb P}\big\{ X_k\not\equiv 0\, ({\rm mod} \ h) \big\}= \infty, \qq {\it for\ any} \ h\ge 2.
\end{equation}

Indeed, otherwise if $\sum_{k=1}^\infty  {\mathbb P}\{ X_k\not\equiv 0\, ({\rm mod} \ h) \}< \infty$ for some $h\ge 2$, then by the  Borel--Cantelli lemma,   on a set of
measure greater than $3/4$, $X_k\equiv 0\, ({\rm mod} \ h)$ for all $k\ge k_0$, say. The new sequence $X'$ defined by $X'_k=0$ if $k< k_0$, $X'_k=X_k$ unless, with partial
sums $S'_n$, verifies ${\mathbb P}\{ S'_n\equiv 0\, ({\rm mod} \ h) \}>3/4$ for all $n$ large enough, thereby contradicting the fact that ${\mathbb P}\{ S'_n\equiv 0\, ({\rm
mod} \ h) \}$ should converge to $1/h$.
 This  establishes \eqref{roz.cond}.
\vskip 5 pt
\vskip 5 pt

 Now the proof of Theorem \ref{rozanov.I} is achieved as follows. Since the  local limit theorem   is applicable to the sequence $X'$, it is also applicable to the \lq\lq
 shifted\rq\rq\,sequence $X''$ where $X''_k=X'_k-m_k$ and $m_k$ are arbitrary integers. This is easy to check. Let $m_k$ be defined by
 $$\max_{0\le m< h} {\mathbb P}\big\{X_k\equiv m\, {\rm mod({\it h})} \big\}= {\mathbb P}\big\{X_k\equiv m_k\, {\rm mod({\it h})} \big\}=  {\mathbb P}\big\{X''_k\equiv 0\,
 {\rm mod({\it h})} \big\}. $$
By \eqref{roz.cond},
\begin{equation*}
\sum_{k=1}^\infty \  {\mathbb P}\big\{ X_k\not\equiv m_k\, {\rm mod} \ h) \big\}=\sum_{k=1}^\infty \ {\mathbb P}\big\{X''_k\not\equiv 0\, {\rm mod({\it h})} \big\}= \infty,
\end{equation*}
 for any $h\ge 2$. Thus
 \begin{equation*}\prod_{k=1}^\infty
\big[ \max_{0\le m< h} {\mathbb P}\big\{X_k\equiv m\, {\rm (mod({\it h})} \big\}\big]=
\prod_{k=1}^\infty \Big(1-  {\mathbb P}\big\{ X_k \not\equiv m_k\, {\rm mod({\it h})} \big\}\Big)=0.
\end{equation*}
 Thus   \eqref{rozn} is obtained, and thereby Theorem \ref{rozanov.I} is proved.

{\color{blue}\begin{remark}  A sequence $X$ and  any subsequence extracted from it, which differs from
it only in a finite number of members, is a.u.d. if and only if    \eqref{rozn} holds.
 \end{remark}}

\begin{remark} \rm
It follows from Gnedenko's theorem   or   Prohorov's theorem  that the necessary condition \eqref{rozn} is also sufficient  in the case
when the
$X_j$ are identically distributed or bounded.
\end{remark}
\vskip 8 pt
Less restrictive conditions on the random variables $X_k$ are  also given in Rozanov \cite{Ro}, under which condition (\ref{rozn}) is sufficient. More precisely, we have the
following result.

\begin{theorem}[{\rm \cite{Ro}, Theorem II}] \label{rozanov.thII} Let $X = \{X_j, j\ge 1\}$ be as in Theorem \ref{rozanov.I}, and let  $a_k={\mathbb E\,} X_k$, $ k\ge 1$.
Assume that condition \eqref{Ro.A} is satisfied and further that
\begin{equation}\label{Ro.B} \frac{1}{b_k^2}\sum_{|j-a_k|\le N} (j-a_k)^2\, {\mathbb P}\{X_k =j\}\ \to \ 1\qq\  \text{\rm as } \ N\to \infty,
\end{equation}
uniformly in $k$.

Then condition \eqref{rozn} is sufficient for the applicability of the local limit theorem    in the strong form to the sequence $X$.
\end{theorem}

 Condition \eqref{Ro.B} together with condition \eqref{Ro.A} imply that  Lindeberg  condition is fulfilled. The proof  consists with a careful application of the method of
 characteristic functions.

\vskip 10 pt

 {\color{blue}}

{\color{blue} In \cite{VMT}, Vo An'\,Zung, Mukhin   and To An'\,Zung, proved   the following interesting theorem generalizing    Theorems \ref{Prohorov} and
\ref{rozanov.thII} of Prohorov and Rozanov respectively.

\vskip 2 pt

 Let $\{X_n, n\ge 1\}$   be a sequence of independent, integer-valued
random variables, and let $\{X'_n, n\ge 1\}$ be the corresponding sequence of
symmetrized random variables. Let  $\tilde p_{km} = \P\{\tilde X_k=m\}$ for each $k$ and $m$,
and $S_n=\sum_{k=1}^n X_k$ for each $n$.
 \vskip 2 pt
 Set for $n\ge 1$ and $M\ge 1$,
 $$ B_n^2(M)=\sum_{k=1}^n\sum_{|m|\le M} m^2 \tilde p_{km}.$$
\begin{theorem} \label{VMT}
Assume that
\vskip 3 pt {\rm (1)} For each $k$, ${\rm Var }(X_k)<\infty$ and $b_n^2 =\sum_{k=1}^n {\rm Var }(X_k)\to \infty$ as $n\to \infty$. Further there exist $M$ and $\g>0$ such
that $B_n^2(M)\ge b_n^2$.

\vskip 3 pt {\rm (2)} There exists a sequence $\{a_n,n\ge 1\}$ such that
$$ \P\{S_n-a_n<x b_n\}\,\to \, \Phi(x) = (2\pi)^{-1/2} \int_{-\infty}^x \exp(-u^2/2)\dd u,\qq  \hbox{as $n\to \infty$}.$$
\vskip 3 pt {\rm (3)} For each $d$, $2\le d\le M$, the sequence $\{S_n, n\ge 1\}$ is   a.u.d.(d).
\vskip 3 pt
\noi Then,
$$\P\{S_n =m\}= \frac{1}{\sqrt{2\pi}\,b_n} e^{-\frac{(m-\E S_n)^2}{2b_n^2}}+ o\big(\frac{1}{b_n}\big).$$
\end{theorem}

\vskip 15 pt

}

{\color{blue} Another important necessary condition {\color{green}for} the validity of the local limit theorem,  due to Gamkrelidze,    takes the form of a remarkable
explicit lower bound of the essential integral majorant $I_4$ (see \eqref{basic.bound}).
\vskip 10 pt
Concerning   results of this {\color{green}kind}, Gamkrelidze remarked in \cite{Gam80}, p.\,275,  (citing also Gnedenko and Kolmogorov \cite{GK}) that: \lq\lq {\it ... A good
deal of probability theory consists of the study of limit theorems, because  in reality the epistemological value of the theory of probability is revealed only by limit
theorems.
Important part of this area consists of the upper estimation of the rate of convergence in the   limit theorem. It is quite reasonable turn to the construction of the lower
estimates. Unfortunately, many mathematicians  working in the field of the theory of  limit theorems pay less attention to such kind of problems.\rq\rq}

\vskip 10 pt
\noi{\tf Gamkrelidze's lower bound.}
Let $\{X_k\}$ be a sequence of independent square integrable random variables taking only integer values and let $S_n=\sum_{k=1}^n X_k,$ $B_n^2={\rm Var}(S_n),$
$M_n={\E}(S_n),$   $\phi_k(t)={\E}(e^{itX_k})$. Recall  (Definition \ref{defllt}) that
$$  \D_n=  \sup_{N=v_0n+Dk }\Big|B_n {\mathbb P}\{S_n=N\}-{D\over  \sqrt{ 2\pi } }e^{-
{(N-M_n)^2\over  2 B_n^2} }\Big| .$$

{\color{red}
\begin{theorem}[{\rm \cite{Gam80} Theorem 1}]
\label{gamt}
Suppose $S_n=X_1+\ldots+X_n$ takes only integer values. For $k\geq 1$
\begin{equation}
\label{egam90}
{{1}\over{4\pi}}
\int_{{{2\pi}\over{2k+1}}\leq |t|\leq \pi}|\phi_n(t)|^2 dt \leq
{{1}\over{2\sqrt{\pi}}B_n}(1-e^{-{{k^2}\over{4B_n^2}}}) + {{2\Lambda_n}\over{B_n}},
\end{equation}
where $\Lambda_n=2.01(\Delta_n + {{1}\over{2\sqrt{\pi}}} e^{-\pi^2B_n^2} )$.
\end{theorem}
\begin{proof}
{\color{blue}We follow the proof given in \cite{Gam80}.
}
Let $\tilde{X}_k$ be an independent copy of $X_k$ and set $\hat{S}_n = \sum_{k=1}^n X_k-\tilde{X}_k$.
Suppose $v_0=0,$ $D=1,$ $B_n\geq 1,$  fix $k$ and  write
\begin{eqnarray}
\label{egam81}
\lefteqn{\P(\hat{S}_n=k) = \sum_{\nu\in{\mathbb Z}} \P(S_n=k+\nu)\P(S_n=\nu)}\nonumber\\
& & =  \sum_{\nu\in{\mathbb Z}} \P(S_n=k+\nu)(\P(S_n=\nu) - Q(x_{n\nu})) + \sum_{\nu\in{\mathbb Z}} (\P(S_n=k+\nu) - Q(x_{n,k+\nu}))Q(x_{n\nu})\nonumber\\
& & + \sum_{\nu\in{\mathbb Z}} Q(x_{n,k+\nu})Q(x_{n\nu})
\end{eqnarray}
where $Q(x_{n\nu})={{1}\over{\sqrt{2\pi}B_n}}\exp\{-{{x_{n\nu}^2}\over{2}}\},$ $x_{n\nu}={{(\nu-M_n)}\over{B_n}}$.
Define
\[
S_1 =  \sum_{\nu\in{\mathbb Z}} Q(x_{n\nu}) \qquad {\rm and} \qquad S_2 = \sum_{\nu\in{\mathbb Z}} Q(x_{n,k+\nu})Q(x_{n\nu}).
\]
Recall formula (5.12) on p.633 in \cite{F}
\begin{equation}
\label{f512}
{{1}\over{\sqrt{2\pi t}}}\sum_{\nu\in{\mathbb Z}}\exp\{-{{1}\over{2t}}(s+2\nu\pi)^2\} =
{{1}\over{2\pi}}\sum_{m\in{\mathbb Z}} e^{-{{m^2t}\over{2}}}\cos{(ms)},
\end{equation}
where $t>0$ and $s\in{\mathbb R}$. Put $s=-2\pi M_n$ and $t=4\pi^2B_n^2$ in (\ref{f512}), thus
\[ S_1= {{1}\over{\sqrt{2\pi}B_n}} \sum_{\nu\in{\mathbb Z}} \exp\{-{{1}\over{2B_n^2}}(\nu-M_n)^2\}
= 1+ 2\sum_{m\geq 1} e^{-2m^2\pi^2 B_n^2}\cos{(m\cdot2\pi M_n)}. \]
It follows from
\[ 2 e^{-2\pi^2 B_n^2}\sum_{m\geq 1} e^{-2(m-1)^2\pi^2 B_n^2} \leq  2e^{-2\pi^2 B_n^2}\left(1+\sum_{m\geq 2} (e^{-2\pi^2 B_n^2})^m\right)
\leq 2e^{-2\pi^2 B_n^2}\left( 1+{{e^{-4\pi^2}}\over{1-e^{-2\pi^2}}}\right)
\]
that
\begin{equation}
\label{egam82}
S_1 = 1+\theta_1\cdot 2.01\cdot e^{-2\pi^2 B_n^2}, \qquad\qquad |\theta_1|<1.
\end{equation}
Further,
\[
S_2=   {{1}\over{2\pi B_n^2}} \sum_{\nu\in{\mathbb Z}} \exp\{-{{1}\over{2B_n^2}}((k+\nu-M_n)^2-(\nu-M_n)^2)\}
\]
and by
\[
(k+\nu-M_n)^2 +(\nu-M_n)^2 = 2(\nu-(M_n -0.5k))^2 +0.5k^2
\]
we have
\begin{equation}
\label{egam83}
S_2 = {{1}\over{2\pi B_n^2}} e^{{-{k^2}\over{4B_n^2}}}\sum_{\nu\in{\mathbb Z}} \exp\{-{{1}\over{B_n^2}}(\nu-(M_n-0.5k))^2).\}
\end{equation}
Now, by (\ref{f512}) with $s=-2\pi(M_n-0.5k)$ and $t=2\pi^2B_n^2$ we get
\begin{equation}
\label{egam84}
S_2 =  {{1}\over{2\sqrt{\pi}B_n}} e^{{-{k^2}\over{4B_n^2}}}
\left(1+ 2\sum_{m\geq 1} e^{-m^2\pi^2 B_n^2}\cos{(m\cdot2\pi (M_n-0.5k))}\right).
\end{equation}
Consequently,
\begin{equation}
\label{egam85}
S_2 = {{1}\over{2\sqrt{\pi}B_n}} e^{{-{k^2}\over{4B_n^2}}}(1+\theta_2\cdot 2.01\cdot e^{-\pi^2 B_n^2}), \qquad\qquad |\theta_2|<1.
\end{equation}
Whence by (\ref{egam81}) we get
\begin{equation}
\label{egam86}
|\P\{\hat{S}_n=k\} - {{1}\over{2\sqrt{\pi}B_n}}e^{-{{k^2}\over{4B_n^2}}}|
\leq B_n^{-1}\Delta_n + B_n^{-1}\Delta_n\cdot S_1 + (S_2 - {{1}\over{2\sqrt{\pi}B_n}} e^{{-{k^2}\over{4B_n^2}}} ).
\end{equation}
Therefore by (\ref{egam82}) and  (\ref{egam85}) we obtain
\begin{equation}
\label{egam87}
|\P\{\hat{S}_n=k\} - {{1}\over{2\sqrt{\pi}B_n}}e^{-{{k^2}\over{4B_n^2}}}|
\leq 2.01(\Delta_n + {{1}\over{2\sqrt{\pi}}} e^{-\pi^2B_n^2} )B_n^{-1}:= \Lambda_n B_n^{-1}.
\end{equation}
On the other hand if we set $\varphi_n(t)=E(e^{itS_n})$ then
\[ \vartheta_k:=\P\{\hat{S}_n=0\} - {{1}\over{2k+1}}\sum_{\nu=-k}^k \P(\hat{S}_n=\nu)
= {{1}\over{2\pi}}\int_{-\pi}^{\pi} (1- {{1}\over{2k+1}}\sum_{\nu=-k}^k e^{-it\nu})|\varphi_n(t)|^2dt.\]
From
\[
\sum_{\nu=-k}^k (e^{it})^{\nu} = 1 +2\sum_{\nu=1}^k \cos{t\nu} = {{\sin{(2k+1)}}\over{\sin{(0.5t)}}}
\]
we get
\begin{equation}
\label{egam88}
\vartheta_k = {{1}\over{2\pi}}\int_{-\pi}^{\pi}\left(1- {{\sin{(2k+1)}}\over{(2k+1)\sin{(0.5t)}}} \right)|\varphi_n(t)|^2dt.
\end{equation}
Now, observe that by (\ref{egam87}) for $|\nu|\leq k$
\[ \P\{\hat{S}_n=0\}\leq  {{1}\over{2\sqrt{\pi}}B_n} + \Lambda_n B_n^{-1}\quad{\rm and}\quad
\P\{\hat{S}_n=\nu\}\geq  {{1}\over{2\sqrt{\pi}}B_n}e^{-{{k^2}\over{4B_n^2}}} - \Lambda_n B_n^{-1}
\]
whence
\begin{equation}
\label{egam89}
\vartheta_k \leq {{1}\over{2\sqrt{\pi}}B_n}(1-e^{-{{k^2}\over{4B_n^2}}}) + 2\Lambda_n B_n^{-1}.
\end{equation}
Finally, note that for $|t|\leq \pi$ we have $|\sin{0.5t}|\geq {{|t|}\over{\pi}}$
thus for $|t|\geq {{2\pi}\over{2k+1}}$ we get
\[  \left|{{\sin{(2k+1)}}\over{(2k+1)\sin{(0.5t)}}}\right| \leq   {{1}\over{(2k+1)|\sin{(0.5t)}|}}
\leq   {{\pi}\over{(2k+1)|t|}}\leq 0.5.
\]
By this, (\ref{egam88}) and (\ref{egam89}) we obtain claimed lower bound.
\end{proof}

}

Gamkrelidze  \cite{Gam2} (see also \cite{P}, p.\,215, supplement 12) proved the following inequality
\begin{equation}
\label{game1}
\left({{1}\over{8\sqrt{\pi}}}+2\right)\Big({{2(1+\sqrt{2{\color{red}\pi}}{\color{red}\Delta_n})}\over{\pi B_n}}+2\D_n\Big) \, \geq \,{{B_n}\over{4\pi}} \int_{\pi\geq |t|\geq
{{2{\color{red} \pi}}\over{2K_n+1}}}\prod_{k=1}^n|\varphi_k(t)|^2\, \dd t,
\end{equation}
where{\footnote{\color{red}{In  \cite{Gam2} in formula (3) $\lambda_n$ is missing after $\sqrt{2\pi}$ (the same misprint is in the definition
of $c_1$ on p. 280 in \cite{Gam80}) and also
in formula (4) should be ${{2\pi}\over{2K+1}}$ under
intergral.}}}
 $$K_n=\Big\lfloor \Big({{2(1+\sqrt{2{\color{red} \pi}}{\color{red}\Delta_n})}\over{\pi B_n}}+2\D_n\Big)^{1/2}  B_n \Big\rfloor. $$
 \vskip 3 pt Assume that $B_n\to \infty$ with $n$.  A consequence is that the   condition
\beq\label{gam.nec.cond}
B_n\int_{\e_n\le t\le 2\pi}\prod_{k=1}^n|\varphi_k(t)|^2 \dd t\ \to \ 0
\eeq
for any positive $\e_n$ which tends to 0 as $n$ tends  to infinity,  is necessary for the local limit theorem to hold. {\color{blue} An application   is given in Gamkrelidze
\cite{Gam4}, p.\,149}.

\vskip 15 pt
From (\ref{game1}) it follows that for i.i.d. sequences
there exists $t_0$ such that
\[   {{2(1+\sqrt{2{\color{red} \pi}}{\color{red}\Delta_n}})\over{\pi B_n}}+2\D_n  \geq {\rm Const.}\, |\phi(t_0)|^{2n+1}. \]
\vskip 5 pt
{\color{blue}
An analog statement to Theorem \ref{gamt} is also proved for multidimensional characteristic functions in  \cite{G2a}.}

{\color{blue}

\vskip 15 pt 
\noi  {\tf Stable limit distributions.}
Let $
X=\{X_i , i\ge 1\}$ be independent  random variables taking only integer values. We assume that  $X_i$ has  distribution function $F_i(x)$ defined by
\begin{equation}\label{mitalauskas.das} F_i(x) =
\begin{cases}
0 \quad & \quad \hbox{for $x\le 0$}
\cr
1-\big[c_i+\a_i(x)\big]\frac{1}{x^\a} & \quad \hbox{for $x>0$,}
\end{cases}
\end{equation}
where  $|\a_i(x)|\le \a(x)$, $\a(x)\to  0$ as $x\to \infty$, $0<c'<c_i<c''<\infty$, $i=1,2,\ldots $ and $0<\a<1$.
\vskip 2 pt
 Let also
\begin{equation*}
B_n= \Big( \sum_{i=1}^n c_i\Big)^{1/\a},\qq S_n= \sum_{i=1}^n X_i, \qq p_{kj}= {\mathbb P}\{X_k= j\}.
\end{equation*}

{\color{blue}Let $G(x)$ denote the stable distribution function for which $M(x)=0$, $N(x)= -1/x^\a$, $\s^2=0$ and $\g(\tau)= \a\tau^{1-\a}/(1-\a)$, in the L\'evy-Khintchin
formula.
\vskip 2 pt
   Rogozin  \cite{Rog} proved that, as $n$ tends to infinity,
\begin{equation}\label{Rog.lim.th}
\P\big\{ \frac{S_n}{B_n}\,<\, x\big\}\,\to\, G(x).
\end{equation}
 This result was later generalized by Banys \cite{Ba}.} \vskip 2 pt
A local limit theorem in the strong form holds (Definition \ref{lltdefstable}) for the sequence $  X$, if
$$B_n\P\{S_n=m\}-g\big( \frac{m}{B_n}\big) \to \ 0, $$
as $m\to \infty$ uniformly in $m$, $-\infty<m<\infty$, where
$$g\big( \frac{m}{B_n}\big)=G'\big( \frac{m}{B_n}\big). $$
}
 \vskip 2 pt
The following characterization was proved by Mitalauskas in \cite{Mit}.
\begin{theorem}\label{mitalauskas.th1}
In order that the sequence $ X$ satisfies a local limit theorem in the strong form, it is necessary and sufficient that Rozanov's condition (see \eqref{rozn}) be fulfilled,
namely that
\begin{equation*}
\sum_{k=1}^\infty \ \min_{0\le m<q} {\mathbb P}\big\{ X_k\not\equiv m\, ({\rm mod} \ q) \big\}= \infty,\qq \hbox{for all integers $q\ge 2$.}
\end{equation*}
\end{theorem}
{\color{blue} We give some hints. Let $\p_n(t)$ be the characteristic function of $S_n$. By the Fourier inversion formula
\begin{equation*}
 2\pi \P\{ S_n=m\}= \int_{-\pi}^\pi e^{-itm} \p_n(t)\dd t=\frac{1}{B_n} \int_{-\pi B_n}^{\pi B_n} e^{-itm/B_n} \p_n(t)\dd t.
 \end{equation*}
\begin{equation*}
 2\pi g\big(\frac{m}{B_n}\big)=    \int_{-\infty}^{\infty} e^{-itm/B_n} f(t)\dd t,
 \end{equation*}
 where
 $$ f(t)=\exp\Big\{-\G(1-\a) \cos \frac{\pi \a}{2}|t|^\a \Big[1- i\frac{t}{|t|}\Big]\Big\}.$$
Thus
\begin{eqnarray} B_n\P\{S_n=m\}-g\big(\frac{m}{B_n}\big) &=&\int_{-\pi B_n}^{\pi B_n} e^{-itm/B_n} \p_n(t)\dd t- \int_{-\infty}^{\infty} e^{-itm/B_n} f(t)\dd t
\cr &=& I_1+I_2+I_3+I_4,
\end{eqnarray}
where
\begin{eqnarray}
I_1 &=& \int_{|t|<\theta}  e^{-itm/B_n} \big(\p_n\big(\frac{t}{B_n}\big)-f(t)\big)\dd t
\cr I_2 &=& \int_{|t|\ge\theta}  e^{-itm/B_n}  f(t) \dd t
\cr I_3 &=& \int_{\theta\le |t|<\e B_n}  e^{-itm/B_n}  \p_n\big(\frac{t}{B_n}\big)\dd t
\cr I_4 &=& \int_{  \e B_n\le |t|\le \pi B_n}  e^{-itm/B_n}\p_n\big(\frac{t}{B_n}\big)\dd t ,
\end{eqnarray}
and $\theta$ is some large positive real. The last integral requires a delicate analysis which starts with
\begin{eqnarray}
 |I_4| &\le & \int_{  \e B_n\le |t|\le \pi B_n} \big|\p_n\big(\frac{t}{B_n}\big)\big|\dd t
 \cr &\le & \int_{  \e/2\pi  \le |u|\le 1/2} \big|\exp\Big\{-p\sum_{k=1}^n\sum_j p_{kj}\sin^2 \pi u j\Big\}\dd u .
\end{eqnarray}
Next by a standard diophantine approximation result (Dirichlet theorem),  one can write
$$ u= \frac{a}{q}+t$$
 where $a,q$ are integers, $(a,q)=1$, $0<q<\tau$, $|t|\le 1/q\tau$.
By reporting, $ |I_4|$ will be bounded from above by a sum over $\frac{a}{q}$ of three integrals, of which the estimation achieves the proof.

 \vskip 4 pt In the same paper the author stated a version of the previous theorem for independent, symmetric,   lattice valued random variables.
Let $
X=\{X_i , i\ge 1\}$ be independent  random variables taking only integer values, and  assume that  $X_i$ has  distribution function $F_i(x)$ defined by
\begin{equation}\label{mitalauskas.das.1} F_i(x) =
\begin{cases}
\big[c_i+\a_i(x)\big]\frac{1}{x^\a}  \quad & \quad \hbox{for $x\le 0$}
\cr
1-\big[c_i+\a_i(x)\big]\frac{1}{x^\a} & \quad \hbox{for $x>0$,}
\end{cases}
\end{equation}
where  $|\a_i(x)|\le \a(x)$, $\a(x)\to  0$ as $x\to \infty$, $0<c'<c_i<c''<\infty$, $i=1,2,\ldots $ and $0<\a<2$.
\vskip 2 pt
\begin{theorem}\label{mitalauskas.th1.tilde}
In order that the sequence $X$ satisfies a local limit theorem in the strong form, it is necessary and sufficient that Rozanov's condition be fulfilled.
\end{theorem}

\vskip 15 pt

\vskip 15 pt

\vskip 2 pt \noi {\tf Gamkrelidze's counter-examples.}
  The first counter-example provides a negative answer to a question raised by Prohorov. It is known that the local limit theorem always implies the integral limit theorem.
  It follows from the intermediate step of the proof of Rozanov's theorem, namely Lemma \ref{roz.cond.equi}, (see also \cite{GK}),  that a necessary condition for the
  validity of the local limit theorem is that the partial sums $S_n$ are {\it a.u.d.}\,. It was also proved that under  relatively general supplementary conditions in
  \cite{Pr}, \cite{Ro}, \cite{Mit} notably, see also Theorems \ref{rozanov.thII}, \ref{VMT},
   \ref{mitalauskas.th1}, \ref{mitalauskas.th1.tilde}, from the integral limit theorem and {\it a.u.d.} property for partial sums follows the local limit theorem. An
   hypothesis was stated by Prohorov that the  integral limit theorem plus the {\it a.u.d.} property for partial sums  is equivalent to the local limit theorem. This
   hypothesis was contradicted by Gamkrelidze \cite{Gam} who constructed  the following remarkable counter-example.

\begin{example} Let $\xi=\{\xi_j,j\ge 1\}$ be a sequence of independent random variables defined as follows:
\vskip 2pt --- The variables $\xi_j$ with odd subscripts $j=2k-1$ are distributed according to a symmetrized Poisson law with characteristic functions
\beq f(\xi_{2k-1},t)=\exp\big\{\Lambda_k\big(\cos(th_k)-1\big)\big\},
\eeq
where
$$\Lambda_k=k^{-1/2},\qq\qq h_k =2\big\lfloor e^{k^2/2}\,.\, k^{-1/4}\big\rfloor.
$$
\vskip 2pt --- The variables $\xi_j$ with even subscripts $j=2k $ take the values
$$ -k, -k+1, -\ldots, 1,2,\ldots, k$$
with probability $1/2k$ each, and so are uniformly distributed mod $k$.
Let $S_n = \sum_{j=1}^n \xi_j$ and set
$$ A_n=\E S_n, \qq B_n^2 = {\rm Var}(S_n),\qq P_n(m)= \P\{S_n =m\}.$$
\vskip 8pt
\vskip 2pt
(i) We have $A_n=0$. Further $S_n/B_n$ is asymptotically normal,
\beq \label{ex.sn.normal} \P\big\{S_n<xB_n\big\}\ \to \ \Phi(x),
\eeq
so that the variables $\xi_j$ satisfy the central limit theorem. In fact, let $n'=2k'-1$ be the largest odd number less than $n$. Then
$$\Lambda_{k'}h_{k'}^2\,\sim\, 4e^{{k'}^2}, $$
and
$${\rm Var}(\xi_{2k'-1}) <B_n^2\le \sum_{k\le k'}{\rm Var}(\xi_{2k})+\sum_{k< k'}{\rm Var}(\xi_{2k-1})+{\rm Var}(\xi_{2k'-1}),$$
or
$$ \Lambda_{k'}h^2_{k'}<B_h^2\le {k'}^3+4k'e^{(k'-1)^2}+ \Lambda_{k'}h^2_{k'}.$$
Thus
$$ B_h^2\sim \Lambda_{k'}h^2_{k'}\sim 4e^{{k'}^2}.$$
Further
$$S_n= \sum_{2k\le n}\xi_{2k}+\sum_{ k<k'}\xi_{2k-1}+\xi_{2k'-1}:= r_n+ \rho_n+\xi_{2k'-1}.$$
One checks that the random variables $r_n$, $\rho_n$ are tending to 0 in probability, and the random variable
$$\xi_{2k'-1}\big(\Lambda_{k'}h^2_{k'}\big)^{-1/2} ,$$
is asymptotically normal. This thus implies \eqref{ex.sn.normal}.
\vskip 2pt
\vskip 2pt
(ii) The sums $S_n$ are {\it a.u.d.}, this follows from the fact that if a random variable $\tau$ is {\it a.u.d.} mod $h$, and $\eta$ is an independent integral-valued random
variable, then the sum $\tau +\eta$ is also {\it a.u.d.} mod $h$. Since the random variables $\xi_{2k}$ are {\it a.u.d.} mod ($k$), it follows that $S_n$ is {\it a.u.d.} mod
$h$, for any $0<h<\lfloor n/2\rfloor$.

\vskip 2pt
\vskip 2pt
(iii) The sequence $\xi$ does not satisfy the local limit theorem. Assume that the contrary is true.
{\color{blue} We use the following   Lemma.
\begin{lemma}[stability Lemma]\label{stability.llt} Assume that the local limit theorem    in the strong form is applicable to the sequence $X$. Then
$$ \P\{ S_n=m'\}\,\sim \, \P\{ S_n=m''\}
$$
for $m'= o(B_n)$, $m''= o(B_n)$.
\end{lemma}
}
Thus by the stability Lemma, we must have
\beq \label{stab.llt}\P\{ S_n=m'\}\,\sim \, \P\{ S_n=m''\}
\eeq
for $m'= o(B_n)$, $m''= o(B_n)$. Choose $m''=0$, $m'=\frac12 h_{k'}= o(B_n)$. Then
\beq \P\{ S_n=m'\}\,= \, \P\{ r_n+ \rho_n=m''-\xi_{2k'-1}\}\,\le  \,2 \P\{ r_n+ \rho_n\ge \frac12 h_{k'}\}\,\le  \,2 \P\{   \rho_n\ge \frac12 h_{k'}-{k'}^2\}.
\eeq
But for $k'\ge k_0$, $\frac14 h_{k'}\ge {k'}^2$, so that
\beq\label{exponent.example.} \P\{ S_n=m'\} \,\le  \,2 \P\{   \rho_n\ge \frac14 h_{k'}\}.
\eeq
Let $\frac14 h_{k'}=T$. By Chebyshev's inequality,
\beq \P\{   \rho_n\ge T\}\,\le \,e^{-uT}\E e^{u\rho_n}=e^{-uT+\sum_{k<k'} \Lambda_k(\cosh uh_k-1)} ,
\eeq
for any $u>0$. Choose $u$ such that the equation
\beq  e^{uh_{k'}-1}\,=\,\frac{h_{k'}}{8\Lambda_{k'-1}h_{k'-1}},
\eeq
is realized, that is
$$u=\frac{1}{h_{k'-1}} \,\log \frac{h_{k'}}{8\Lambda_{k'-1}h_{k'-1}}.$$
For such an $u$,
\ben \label{exponent.example} -uT+\sum_{k<k'} \Lambda_k(\cosh uh_k-1) &\le &-uT + k'\Lambda_{k'-1}\, e^{uh_{k'-1}}
\cr &=&-\frac14 \,\frac{h_{k'}}{h_{k'-1}} \,\log \frac{h_{k'}}{8\Lambda_{k'-1}h_{k'-1}}+ \frac18\,\,\frac{h_{k'}}{h_{k'-1}}\,k'.
\een
Now observe that
$$\frac{h_{k'}}{h_{k'-1}}\,\sim\,
e^{-1/2}e^{k'}, \qq \quad \log \frac{h_{k'}}{8\Lambda_{k'-1}h_{k'-1}}\,\sim \, k'.$$
So that the right-hand side of  \eqref{exponent.example} is equivalent to
$$-\frac18 \,\frac{h_{k'}}{h_{k'-1}}\,k' \,\sim \, -\frac18 \,e^{-\frac12} \,{k'}\,e^{k'}$$
and so, does not exceed
\beq\label{exponent.example1}-\frac1{16} \, {k'}\,e^{k'},
\eeq
for $k\ge k_1$, say. It follows from \eqref{exponent.example.}, \eqref{exponent.example1} that for $n$ sufficiently large
\beq\label{exponent.example2}\P\{ S_n=m'\}\le 2\exp\big\{-\frac1{16} \, {k'}\,e^{k'}\big\}.
\eeq
But at the same time
\beq\label{exponent.example3}\P\{ S_n=0\}\,\sim\,\frac{1}{\sqrt{2\pi} B_n}\,\sim\,\frac{1}{2\sqrt{2\pi} }\,e^{-{k'}^2/2},
\eeq
which produces a contradiction with \eqref{stab.llt}. The local limit theorem is therefore inapplicable to the sequence $\xi$.
\end{example}

\vskip 15 pt

We now pass a second  remarkable  counter-example {\color{green}due to Azlarov \cite{Az}}.
{\color{red} Sometimes the local limit theorem (LLT) is equivalent to the central limit theorem (CLT) (see e.g.
\cite{Az} (also Chapter VII, \S4, 14 on p.216 in \cite{P}) and  \cite{Pr63}).
Let $\xi_1,\xi_2,\ldots,\xi_n\ldots$ be such that
\[\P\{\xi_k=m\}={{1}\over{2N_k+1}},\quad m\in\{-N_k, -N_k+1,\ldots,-1,0,1\ldots,N_k\},\quad N_k\in\N\cup\{0\}. \]
Set
\[ \sigma_k^2={\rm Var}(\xi_k)= {{N_k(N_k+1)}\over{3}},
\]
\[ B_n^2= \sum_{k=1}^n \sigma_k^2 = {{1}\over{3}}\sum_{k=1}^n N_k(N_k+1),
\]
\[ \lambda_n = B_n\min_{1\leq k\leq n} N_k^{-1}.
\]
\begin{theorem}[{\rm \cite{Az}}]
\label{azth}
The condition $\lim_n\lambda_n=\infty$ is necessary and sufficient for $\{\xi_k\}$
to satisfy CLT and LLT. Moreover, for every $n\geq 2$ such that $\lambda_n\geq 2$  we have $\Delta_n\leq C\lambda_n^{-1}$, where
$C>0$ is an absolute constant.
\end{theorem}
\begin{proof}
We have
\begin{eqnarray*}
{{1}\over{4}}\lambda_n^{-2} & < & B_n^{-2}\E^2(\max_{1\leq k\leq n}|\xi_k|)
\leq  B_n^{-2}\E (\max_{1\leq k\leq n}\xi_k^2)\\
&\leq &  B_n^{-2}\max_{1\leq k\leq n}N_k\E(\max_{1\leq k\leq n}|\xi_k|) \leq 4\lambda_n^{-1}.
\end{eqnarray*}
Thus $\lambda_n\to\infty$ is equivalent to $B_n^{-2}\E(\max_{1\leq\k\leq n}\xi_k^2)\to 0$ which in turn is equivalent to
the Lindeberg condition. Therefore in view of Theorem 22 on pp.99--100 in \cite{P}
it is enough to prove the second statement of Theorem \ref{azth}.
For this note that
\[ f_k(t)=\E(e^{it\xi_k})= {{1}\over{2N_k+1}}{{\sin{{{2N_k+1}\over{2}}}t}\over{\sin{{{t}\over{2}}}}}.\]
Suppose that we proved that for $t\in[0,{{2}\over{2N_k+1}}]$
\begin{equation}
\label{aze1}
|f_k(t)| \leq e^{-Ct^2\sigma_k^2}
\end{equation}
and for $t\in[{{1}\over{2N_k+1}},\pi)$
\begin{equation}
\label{aze2}
|f_k(t)|\leq e^{-C},
\end{equation}
where $C$ are some constants (possibly different). Let us decompose
\[ 2\pi B_n\P\{S_n=m\} = \int_{|t| < {{1}\over{4L_n}}} e^{-it{{m}\over{B_n}}}\prod_{k=1}^n f_k({{t}\over{B_n}})dt
+ \int_{{{1}\over{4L_n}} \leq |t|\leq \pi B_n} e^{-it{{m}\over{B_n}}}\prod_{k=1}^n f_k({{t}\over{B_n}})dt = I_1 + I_2,
\]
where $L_n=B_n^{-3}\sum_{k=1}^n \E(|\xi_k|^3)$.
We have
\[ I_1 = \int_{|t|< {{1}\over{4L_n}}} e^{-it{{m}\over{B_n}}}
\left( \prod_{k=1}^n f_k({{t}\over{B_n}}) - e^{-{{t^2}\over{2}}} \right)dt
+ \sqrt{2\pi}e^{-it{{m^2}\over{B_n^2}}}
- \int_{|t|\geq {{1}\over{4L_n}}} e^{-it{{m}\over{B_n}}}  e^{-{{t^2}\over{2}}} dt.
\]
By Lemma 1 on p.109 in \cite{P} the absolute value of the first integral is bounded by $CL_n$
while
\[ \left|  \int_{|t|\geq {{1}\over{4L_n}}} e^{-it{{m}\over{B_n}}}  e^{-{{t^2}\over{2}}} dt \right|
\leq 8L_n \int_{t\geq {{1}\over{4L_n}}} t e^{-{{t^2}\over{2}}}dt = 8L_n e^{-{{1}\over{32L_n}}}.
\]
Note that $L_n\leq B_n^{-3}\sum_{k-1}^n N_k \E(\xi_k^2)\leq B_n^{-1}\max_{1\leq k\leq n}N_k=\lambda_n^{-1}$.
In view of this
\begin{equation}
\label{aze3}
|I_1| \leq \sqrt{2\pi} e ^{-{{1}\over{2}} {{m^2}\over{B_n^2}}} + C L_n + C L_n e^{-{{1}\over{32L_n^2}}}
\leq C\lambda_n^{-1}.
\end{equation}
Now, we fix $n$ (satisfying $\lambda_n\geq 2$) and rearrange
$\xi_1,\ldots,\xi_n$ such that $N_1 \geq N_2\geq\cdots\geq N_n$. Whence  $\sigma_1^2 \geq \sigma_2^2\geq\cdots\geq \sigma_n^2$
and $\lambda_n= N_1^{-1}B_n$. So that
\begin{equation}
\label{aze4}
|I_2|\leq 2B_n\int_{{{1}\over{4N_1}}\leq t\leq \pi}\prod_{k=1}^n |f_k(t)|dt = \sum_{i=0}^n A_i,
\end{equation}
where  for $i=1,2,\ldots,n-1$
\[A_i = 2B_n \int_{{{1}\over{N_i +1/2}}\leq t < {{1}\over{N_{i+1}+1/2}}} \prod_{k=1}^n|f_k(t)|dt \]
and
\[
A_0 = 2B_n \int_{{{1}\over{4N_1}}\leq t < {{1}\over{N_{1}+1/2}}} \prod_{k=1}^n|f_k(t)|dt,\quad
A_n = 2B_n \int_{{{1}\over{N_n +1/2}}\leq t < \pi} \prod_{k=1}^n|f_k(t)|dt.
\]
For $i=1,2,\ldots,n$ we have $0\leq t \leq {{1}\over{N_1+1/2}}\leq {{1}\over{N_i+1/2}}$, therefore
by (\ref{aze1}) for some $C$
\begin{equation}
\label{aze5}
A_0 \leq 2\int_{t\geq {{\lambda_n}\over{4}}} e^{-Ct^2}dt \leq {{C}\over{\lambda_n}}
\end{equation}
and
\[
A_1\leq 2B_n \int_{{{1}\over{N_1+1/2}}\leq t< {{1}\over{N_2+1/2}}} \exp\{-Ct^2(B_n^2-\sigma_1^2)\}dt.
\]
In view of
\[
B_n^2 - \sigma_1^2 = B_n^2 - {{N_1(N_1+1)}\over{3}}\geq B_n^2 - {{2}\over{3}}{{B_n^2}\over{\lambda_n^2}}\geq {{5B_n^2}\over{6}}
\]
we get
\begin{equation}
\label{aze6}
A_1\leq 2B_n \int_{{{1}\over{N_1+1/2}}\leq t< {{1}\over{N_2+1/2}}} \exp\{-Ct^2 B_n^2)\}dt
\leq 2 \int_{t\geq {{1}\over{N_1+1/2}}} e^{-Ct^2} dt \leq {{C}\over{\lambda_n}}.
\end{equation}
By the inequality (\ref{aze2}) for $k\geq 2$ and $t\in[{{1}\over{N_k+1/2}},\pi)$ we have  (since $N_k\leq N_2$)
\[
|f_k(t)|\leq \exp\{-C {{(N_k+1/2)^2}\over{(N_2+1/2)^2}}\}.
\]
Therefore
\[
A_k\leq 2B_n\exp\{-C {{B_n^2-\sigma_1^2-\sigma_2^2}\over{(N_2+1/2)^2}}\} \int_{{{1}\over{N_k+1/2}}\leq t< {{1}\over{N_{k+1}+1/2}}} |f_1(t)||f_2(t)|dt
\]
and
\[
\sum_{k=2}^{n-1} A_k\leq 2B_n\exp\{-C {{B_n^2-\sigma_1^2-\sigma_2^2}\over{(N_2+1/2)^2}}\} \int_{{{1}\over{N_2+1/2}}\leq t< {{1}\over{N_n+1/2}}} |f_1(t)||f_2(t)|dt.
\]
Taking into account inequalities
\[
B_n^2-\sigma_1^2 -\sigma_2^2 \geq B_n^2 - 2\sigma_1^2\geq B_n^2\left(1-{{4}\over{3\lambda_n^2}}\right)\geq {{2}\over{3B_n^2}}
\]
and using inequality $\pi\sin{x} \geq 2x,$ $x\in[0, \pi/2],$ we get
\[
\left(\int_{{1}\over{N_2+1/2}}^{{1}\over{N_n+1/2}}\! |f_1(t)||f_2(t)|dt\right)^2
\leq {{\pi^4}\over{16}}{{1}\over{(N_1+1/2)^2(N_2+1/2)^2}}\left(\int_{{1{}\over{N_2+1/2}}}^\infty\!\! {{dt}\over{t^2}}\right)^2
\leq {{\pi^4}\over{16}}{{1}\over{(N_1+1/2)^2}}.
\]
Whence
\begin{equation}
\label{aze7}
\sum_{k=2}^{n-1} A_k \leq 2B_n{{\pi^2}\over{4}}{{1}\over{N_1+1/2}}e^{-C{{B_n^2}\over{(N_1+1/2)^2}}} \leq C\lambda_n e^{-C\lambda_n^2}
\leq {{C}\over{\lambda_n}}.
\end{equation}
Analogously
\begin{equation}
\label{aze8}
A_n\leq  2B_n\exp\{-C {{B_n^2-\sigma_1^2-\sigma_2^2}\over{(N_2+1/2)^2}}\} \int_{{{1}\over{N_n+1/2}}}^\pi |f_1(t)||f_2(t)|dt
\leq {{\pi^2}\over{4}}{{2B_n}\over{N_1+1/2}}e^{-C{{B_n^2}\over{(N_2+1/2)^2}}}
\leq {{C}\over{\lambda_n}}.
\end{equation}
Now, in view of (\ref{aze3}), taking into account (\ref{aze5}),(\ref{aze7}) and (\ref{aze8}) we obtain the second statement
of Theorem \ref{azth}.

It remains to proof (\ref{aze1}) and  (\ref{aze2}). As for the first one observe that for $C={{1}\over{12}}$
and $x\in[0,1]$ we have
\begin{equation}
\label{aze9}
{{\sin{x}}\over{x}}\leq e^{-Cx^2}\quad {\rm and}\quad {{\sin{{{x}\over{2}}}}\over{{{x}\over{2}}}}\geq e^{-6C(x/2)^2}.
\end{equation}
By (\ref{aze9}) for $t\in [0,{{1}\over{N_k+1/2}})$ we obtain
\[
|f_k(t)|\leq {{t/2}\over{\sin{(t/2)}}}e^{-Ct^2(N_k+1/2)^2}
\leq {{t/2}\over{\sin{(t/2)}}}e^{-Ct^2/4}e^{-Ct^2(N_k^2+N_k)}
\leq {{t/2}\over{\sin{(t/2)}}}e^{-Ct^2/4}e^{-Ct^2\sigma_k^2}.
\]
By $3\sigma_k^2=N_k(N_k+1)\geq 2$ for any k we have
\[
|f_k(t)| \leq {{t/2}\over{\sin{(t/2)}}}e^{-Ct^2/4}e^{-4/3Ct^2}e^{-Ct^2\sigma_k^2}
\leq {{t/2}\over{\sin{(t/2)}}}e^{-6Ct^2/4} e^{-Ct^2\sigma_k^2}\leq e^{-Ct^2\sigma_k^2}.
\]
Finally we turn to the proof of (\ref{aze2}). For any $t\in(0, \pi{{N_k-1/2}\over{N_k+1/2}}]$ we have
\[
|f_k(t+{{\pi}\over{N_k+1/2}})|\leq {{1}\over{2N_k+1}}\left| {{\sin{\left[(N_k+1/2)t+\pi\right]}}\over{\sin{(t/2+ {{\pi}\over{2N_k+1}})}}}\right|
\leq {{1}\over{2N_k+1}} {{|\sin{((N_k+1/2)t)}|}\over{\sin{(t/2)}}} =|f_k(t)|
\]
hence for $1\leq (N_k+1/2)t\leq \pi(N_k+1/2)$ the following inequality is true
\[
|f_k(t)| \leq \max_{s\in\left[{{1}\over{N_k+1/2}}, {{\pi+1}\over{N_k+1/2}} \right]}|f_k(s)|.
\]
One can prove that $|f_k(s)|$ is decreasing for $s\in\left[{{1}\over{N_k+1/2}}, {{\pi}\over{N_k+1/2}} \right]$.
Moreover, if $s\in\left[{{1}\over{N_k=1/2}}, {{\pi+1}\over{N_k+1/2}} \right]$ then
\[
|f_k(s)|\leq {{s/2}\over{\sin{(s/2)}}}\left| {{\sin{((N_k+1/2)s)}}\over{(N_k+1/2)s}}\right|\leq
 {{1}\over{\pi}}{{\pi/2}\over{\sin{(\pi/2)}}}=1/2.
\]
Consequently, we obtain
\[
|f_k(t)|\leq \max\{1/2, |f_k({{1}\over{N_k+1/2}})|\}.
\]
Thus by (\ref{aze1}) it follows
\[
|f_k({{1}\over{N_k+1/2}})| \leq \exp\{-C{{\sigma_k^2}\over{(N_k+1/2)^2}}\}\leq e^{-C/6}<1.
\]
Theorem is proved.
\end{proof}
\begin{remark}
\label{azr}
From the proof of Theorem \ref{azth} it follows that the conclusions are true if i) $n-s(n)\to\infty$, where $s(n)$ is the
number of random variables in $\{\xi_1,\xi_2,\ldots,\xi_n\}$ for which $N_k=0$; ii) we drop  the assumption of symetricity around the
origin.  Also $N_k$ can dependent on $n$ in such a way that still $\lim_{n\to\infty}\lambda_n=\infty$.
\end{remark}
}
\vskip 8 pt This equivalence {\color{red} in  Theorem \ref{azth}}
is however not true in general, as shown by  Gamkrelidze, who built   in \cite{Gam4} an unexpected counter-example   based on the properties of the Fibonacci sequence.
\vskip 3 pt

\begin{example}
Let $[1,1,\ldots,1\ldots]$ be a continued fraction representation of the number $\p= (1+\sqrt 5)/2$. Let $P_j/Q_j$ denote the convergents of the continued fraction of $\p$.
The sequence $\{P_j, j\ge 0\}$ is the Fibonacci sequence and $P_{j-1}=Q_j$ for $j\ge 1$.
 \vskip 3 pt
 Consider now the sequence of independent integer-valued random variables defined as follows
 \ben \label{F1} 1. & & \xi_1,\ldots, \xi_{n_1},\cr
  2. & &  \xi_{n_1+1},\ldots,  \xi_{n_1+n_2},\cr
   & & \ldots \cr
 j. & & \xi_{n_1+n_2+n_{j-1}+1},\ldots, \xi_{n_1+n_2+\ldots+n_{j}},\cr
 & &  \ldots \cr
 \een
Any  random variable $\xi_r$ of the line $j$   takes values $0, Q_j, P_j$ with respective probabilities $(P_j-2)/P_j$, $1/P_j$, $1/P_j$.
\vskip 8 pt For this sequence, the  u.a.d. property  is satisfied, as well as the central limit theorem. Further it is  uniformly asymptotically negligible, namely
$$\sup_{1\le j\le n}\P\{|\xi_j|>x\}\to 0\qq {\rm as } \ n\to \infty,\  {\rm for\ all} \   x>0. $$
But the local limit theorem fails to hold.
\vskip 8 pt
Let us analyze the characteristic properties of  this example.
    At first,
   \begin{eqnarray*}\E \xi_j \,=\,\frac{P_j+Q_j}{P_j}\,=\,\frac{P_{j+1}}{P_j}
   \cr
{\rm Var}(\xi_j) \,=\,\frac{P_j^2+Q_j^2}{P_j}-\frac{(P_j+Q_j)^2}{P_j}
 \end{eqnarray*}
 We note
 that
 $${\rm Var}(\xi_j) \,\ge \,(1-1/P_j)\frac{P_j^2+Q_j^2}{P_j}\,>\,\frac13\, \frac{P_j^2+Q_j ^2}{P_j} $$   {\color{green} The}  characteristic function verifies
 $$ f(t,\xi_j)=\E e^{it \xi_r}= \frac{P_j-2+e^{itQ_j}+e^{itP_j}}{P_j},$$
 $$ \big|f(t,\xi_r)\big|^2
 = \frac{(P_j-2)^2+2 }{P_j^2}+ \frac{2}{P_j^2}
 \cos \, t(P_j-Q_j)+ \frac{2(P_j-2)}{P_j^2}\big( \cos \, t Q_j+ \cos \, tP_j \big).$$

Let
$$n_j = \big\lfloor P_j^{3/2} \big\rfloor .$$
Let $N_k = n_1+\ldots+n_k$ and
 $$ B_{N_k}^2  \,=\,{\rm Var}(S_{N_k})\,=\,\sum_{j=1}^k \Big( \big\lfloor P_j^{3/2} \big\rfloor+1\Big)\, {\rm Var}(\xi_{N_k})  \,=\, \mathcal O (P^{5/2}_k)
.$$
\vskip 3 pt (i) The sequence has the u.a.n. property. Let $n$ be arbitrary and let $k$ be such that $N_{k-1}<n\le N_k$. Then
$$ \max_{1\le j\le n} \big|\xi_j - \E \xi_j \big|\,\le\, P_k, \qq \quad
B_{N_k}^2\,\ge \, 3^{-1}\sum_{j=1}k(P_j^2+Q_j^2) \frac{n_j}{P_j}.
$$
Whence
$$\frac{ \max_{1\le j\le n} \big|\xi_j - \E \xi_j \big|}{B_n}\,\le\, C\, P_{k-1}^{-1/4}\ \to\ 0,$$
as $n$ tends to infinity.
\vskip 3 pt (ii)  The Liapunov's condition
$$\frac{1}{B_n}\, \Big(\sum_{1\le j\le n} \big|\xi_j - \E \xi_j \big|^{2+\d}\Big)^{1/(2+\d)}\ \to \ 0
$$
as $n$ tends to infinity, holds for some $\d>0$.

\vskip 3 pt (iii) The a.u.d. condition is fulfilled.

\vskip 3 pt (iv) The necessary condition \eqref{gam.nec.cond} for the validity of the local limit theorem fails to hold.

\end{example}

\vskip 15 pt
\noi{\tf Integral limit theorem and local limit theorem.} We shall consider in this subsection a wider setting  and  first introduce a  definition.
\begin{definition}\label{ilt} Let $\{X_{nk}, k=1,\ldots,k_n\}$ be series of independent random variables with values in   $\Z^s$ and set
\beq\label{}S_n= \sum_{k=1}^{k_n} X_{nk}, \qq \quad n=1,\ldots
\eeq  We say that the integral limit theorem ({\rm i.l.t.}) holds if there exist $A_n\in \R^s$ and real $B_n\to \infty$ such that the sequence of distributions of
$(S_n-A_n)/B_n$ converges weakly to an absolutely continuous distribution with density $g(x)$, which is uniformly continuous in $\R^s$.
\end{definition}

\vskip 3 pt

}
This preliminary raises the question of finding all  the limit distributions for these sums. This  has been considerably investigated and we refer to Petrov \cite{P},
Ch.\,IV, \S\,1, where these questions are   well exposed. Some natural restrictions are necessary.

\begin{definition}\label{infinite.smallness.def}
We say that the series of random variables satisfies a {\it condition of infinite smallness} if
\beq\label{infinite.smallness}
\max_{1\le k\le k_n} \P\{|X_{nk}|\ge \e\}\ \to \ 0,
\eeq
for every fixed $\e>0$.
\end{definition}

\vskip 3 pt
We have the following fundamental result.

\begin{theorem}[\cite{P},\, Th.\,1,\,p.\,72]\label{infinite.smallness.th}
The set of distribution functions that are limits  (in the sense of weak convergence) of the distributions of sums  $\sum_{k=1}^{k_n} X_{nk}$ of independent random variables
satisfying the condition of infinite smallness, coincides with the set of infinitely divisible distributions functions. \end{theorem}

 \vskip 7 pt
 We denote by the class $L$ the set of all distribution functions which are limit of distribution of sums
 \beq \frac{1}{a_n} \sum_{k=1}^n X_k -b_k
 \eeq
where $\{X_n, n\ge 1\}$ is a sequence of independent random variables,
  $\{a_n, n\ge 1\}$, $\{b_n, n\ge 1\}$ are sequences of constants with $a_n>0$, and the following condition is satisfied,
  \beq \max_{1\le k\le n}\P\{|X_k|\ge \e a_n\}\ \to \ 0, \qq \hbox{for every fixed $\e>0$.}
\eeq
\vskip 7 pt By the previous theorem the class $L$ is a subset of the set of infinitely divisible distributions functions. We have the following characterization.
\begin{theorem}[\cite{P},\, Th.\,9,\,p.\,85]\label{classL.th}
An infinitely divisible distribution function $F(x)$ belongs to the class $L$ if and only if the corresponding spectral function $L(x)$ is continuous at every point $x\neq 0$
and has left and right continuous derivatives and the function $xL'(x)$ is non-increasing {\rm (}here $L'(x)$ denotes either the left or the right derivative{\rm )}.
\end{theorem}

}
{\color{blue}
Let $S_n= \sum_{k=1}^{k_n} X_{nk}$, where $\{X_{nk}, k=1,\ldots,k_n\}$ are series of independent random variables with values in   $\Z^s$, and satisfying the {\rm i.l.t.},
namely (Definition \ref{ilt})  there exist $A_n\in \R^s$ and real $B_n\to \infty$ such that the sequence of distributions of $(S_n-A_n)/B_n$ converges weakly to an absolutely
continuous distribution with density $g(x)$, which is uniformly continuous in $\R^s$.
\begin{definition}\label{LLTgeneral.} The  local limit theorem is valid if
\beq\label{ilt.llt}
\P\{S_n=m\}=B_n^{-s} g\Big(\frac{m-A_n}{B_n}\Big) + o(B_n^{-s}),
\eeq
uniformly in $m\in \Z^s$.
\end{definition}
As already observed,   the local limit theorem always implies the integral limit theorem.
\vskip 3 pt

{\color{red} By Theorem in Sato \cite{Sa82} if $\mu$ is a {\color{green} genuinly}} $s$-dimensional distribution of class $L$, then $\mu$ is absolutely continuous (with
respect to Lebesgue measure). In this case its characteristic function $\hat{\mu}(t),$ $t\in\R^s$, is uniquely represented in the form
\beq\label{Levy-Khintchine}
\hat{\mu}(t) = \exp\left\{ iat - {\bf A}(t) + \int_{\R^s} \left( e^{itu} -1 - {{itu}\over{1+|u|^2}}\right)\nu(du)\right\},
\eeq
where $a\in\R^s$, ${\bf A}(t)$ is a non-negative quadratic form, and $\nu$ is a measure (called {\it L\'evy measure}) on $\R^s$ such that
$\nu(\{0\})=0$ and
\beq
\int_{\R^s} {{|u|^2}\over{1+|u|^2}} \nu(du) < \infty.
\eeq
}

\vskip 15 pt
\noi {\tf Mukhin's necessary and sufficient condition.}     Let $\{S_n,n\ge 1\}$ be a sequence of {\color{red}integer-valued} random variables  such  that an integral limit
theorem ({i.l.t.}) holds:  there exist $a_n\in \R$ and real $b_n\to \infty$ such that the sequence of distributions of $(S_n-a_n)/b_n$ converges weakly to an absolutely
continuous distribution $G$ with density $g(x)$, which is uniformly continuous in $\R$. Let
\beq \label{en}\e_n:=\sup_{x\in \R}\Big|\P\Big\{\frac{S_n-a_n}{b_n}<x\Big\}-G(x)\Big|\ \to \ 0.
\eeq\vskip 3 pt
Mukhin \cite[Th.\,1]{Mu2} has shown in 1984 the following important result   relating the {i.l.t.} to  {l.l.t.}\,.

\begin{theorem}\label{Mukhin[NSC]} The following assertions are equivalent.
\vskip 3 pt {\rm (A)} \quad There exists a sequence of integers $v_n=o(b_n)$ such that
\beq \sup_{m}\Big|\P\{S_n=m+v_n\big\}-\P\{S_n=m \big\}\Big|\,=\,\,o\Big(\frac{1}{b_n}\Big)
\eeq
\vskip 3 pt
{\rm (B)} \quad
\beq\P\{S_n=m\}\,=\, \frac{1}{b_n} \,g\Big(\frac{m-a_n}{b_n}\Big) +  \,o\Big(\frac{1}{b_n}\Big),
\eeq
 \end{theorem}

 \vskip 10 pt  However from the proof available   in  \cite{Mu2},  we could only  obtain the following slightly weaker result, see Weber \cite{W3}.

\begin{theorem}\label{Mukhin[NSC]weaker}
Let $v_n$ be a  sequence of positive  integers such that $v_n=o(b_n)$.  The following assertions are equivalent.
\vskip 3 pt
 {\rm (A')} \quad  \beq \sup_{m, k\in\Z\atop |m-k|\le v_n}\Big|\P\{S_n=m\big\}-\P\{S_n=k \big\}\Big|\,=\,\,o\Big(\frac{1}{b_n}\Big)+ {\color{blue} \frac{1}{v_n}\,\mathcal
 O(\e_n)},\qq \quad n\to \infty.
\eeq

\vskip 3 pt
{\rm (B')} \quad
\beq\sup_{m}\Big|\P\{S_n=m\} -\frac{1}{b_n} \,g\Big(\frac{m-a_n}{b_n}\Big)\Big| \,=\,   \,o\Big(\frac{1}{b_n}\Big)+ {\color{blue} \frac{1}{v_n}\,\mathcal O(\e_n)},\qq \quad
n\to \infty.\eeq
 \end{theorem}
\vskip 3 pt
Choosing
$v_n= \max\{1, [\sqrt \e_n b_n]\}$
 we get the following

\begin{corollary}\label{Mukhin[NSC]w}  A necessary and sufficient condition for the local limit theorem in the usual form to hold is
\beq \sup_{m, k\in\Z\atop |m-k|\le \max\{1, [\sqrt \e_n b_n]\}}\Big|\P\{S_n=m\big\}-\P\{S_n=k \big\}\Big|\,=\,\,o\Big(\frac{1}{b_n}\Big).
\eeq
 \end{corollary}
\vskip 2 pt
 \begin{remark}
By Theorem 2 in \cite{Mu2}, for arbitrary $\D>0$, the relation
\beq \label{(3)}\P\{S_n\in[x,x+\D]\}= \frac{\D}{b_n}\,g\Big(\frac{m-a_n}{b_n}\Big) +  \,o\Big(\frac{1}{b_n}\Big),\qq \quad n\to \infty,
\eeq
holds uniformly in $x\in\R$ if and only if for any $\l>0$, $v_n>0$, $v_n=o(b_n)$ we have
\beq  \label{(4)}\sup_{x } \Big|\P\big\{S_n\in[x+v_n,x+v_n+\l]\big\}-\P\big\{S_n\in[x+v_n,x+v_n]\big\}\Big|=  \,o\Big(\frac{1}{b_n}\Big),\quad\,    n\to \infty.
\eeq
It is not clear how these assertions should imply our condition (A'), and conversely.
\end{remark}
\vskip 3 pt

\begin{remark}
   The key quantity
\beq \sup_{m \in\Z }\big|\P\{S_n=m+k\big\}-\P\{S_n=k \big\}\big| ,
\eeq for partial sums of integer valued r.v.'s was thoroughly investigated by Mukhin in several works,   in \cite{Mu3}, \cite{Mu1}, \cite{Mu} and \cite{Mu2a} notably, using
structural characteristics of the summands.    Mukhin wrote in \cite{Mu2}: \lq\lq ... getting from here more general sufficient conditions turns out to be difficult in view
of the lack of good criteria for relations \eqref{(3)} and \eqref{(4)}.  Working with asymptotic equidistribution properties are more convenient in this respect\,\rq\rq.
Sufficient conditions of different type for the validity of the l.l.t. are given in \cite{Mu2b}.
\end{remark}

\vskip 15 pt
\noi {\tf Characteristics of a random variable.}
 The local limit theorem  is often studied
 by using various structural characteristics, which are interrelated.  There exists a subsequent literature.
This unfortunately does not include a survey, and we  only report some of the background here.

\medskip

{\color{blue}
 {{(1)}} Let $X$ be an integer-valued random variable, that is, $X$ is  taking with values in $\N$, and we do not     specify whether  $X$ might take   values in some
 sub-lattice,  almost surely.
   \vskip 5 pt

The \lq\lq smoothness\rq\rq characteristic
\begin{eqnarray}\label{delta}  \d_X =\sum_{m\in \Z}\big|{\mathbb P}\{X=m\}-{\mathbb P}\{X=m-1\}\big|,
\end{eqnarray}
  was  introduced and much investigated by Gamkrelidze in \cite{Gam1},    \cite{G2} \cite{G1}, and  in  \cite{GamSher} with Shervashidze.

  \vskip 5 pt    Let us list some properties of this characteristic.

\begin{proposition}[\cite{Gam1}]\label{deltaX.prop}
  {\rm (1)}  $\d_X\le 2$.
\vskip 3pt {\rm (2)}  If $\d_X<2$, then there is an integer $m_0$ such that $\P\{X=m_0\}>0$ and  $\P\{X=m_0-1\}>0$. Thus the span of $X$ is 1.
\vskip 3pt {\rm (3)} $\d_X\ge 2\max_{m\in\Z} \P\{X=m\}$. If the distribution of $X$ is unimodal, namely if there is an $m_0$ such that $\P\{X=m+1\}\le \P\{X=m \}$ for $m\ge
m_0$, and $\P\{X=m-1\}\le \P\{X=m \}$ for $m\le 0$, then $\d_X= 2\max_{m\in\Z} \P\{X=m\}$. (So is the case  for stable laws.)
\vskip 3pt {\rm (4)} Let $X_1,X_2$ be two independent random variables.  Then $\d_{X_1+X_2}\le \min\big\{\d_{X_1}, \d_{ X_2}\big\}$.
\vskip 3pt {\rm (5)}
  {\rm (\cite{Gam1}, p.\,825)} Let $X_i,i\ge 1$ be a sequence of i.i.d. integer-valued random variables with maximal span 1, and let for each
   $n$, $S_n=X_1+\ldots +X_n$. Then  $\d_{S_n} \to 0$ as $n\to \infty$. If the random variables are moreover uniformly bounded, $|X_i|\le L$, then
\beq \label{delta.sn.bdd}
\d_{S_n}= \frac{2}{\sqrt{{\rm Var( X_1)}2\pi n}}+ \mathcal O\Big(\frac{1}{n}\Big).
\eeq
\end{proposition}

 \vskip 15 pt Let now $Y$ be an integer-valued random variable. The following inequalities are implicit in \cite{G2}, p.\,428. We have
 \beq \big|\P\{Y\equiv 0\ {\rm mod}\ h\}-h^{-1}\big|\le \d_Y, \qq \big|\P\{Y\equiv j\ {\rm mod}\ h\}-h^{-1}\big|\le (2-h^{-1})\, \d_Y,
 \eeq
for any $j=1,\ldots, h-1$, and $h\ge 2$.
\vskip 2 pt
Consequently, a sufficient condition for a sequence of   integer-valued random variables $Y_n$ to be asymptotically uniformly distributed   is that $\d_{Y_n} \to 0$ as $n\to
\infty$.
In particular    if for each $n$, $S_n=\sum_{k=1}^n X_k$, where $\{X_n, n\ge 1\}$ is  a sequence of independent, integer-valued
random variables, a sufficient condition for $S_n$ to be a.u.d (Definition \ref{aud}) is that $\d_{S_n} \to 0$ as $n\to \infty$.

 \vskip 5 pt
 \vskip 5 pt
The  smoothness  characteristic $\d_X$ is connected to the characteristic function $\p_X(t)=
 {\mathbb E\,} e^{i t X}$ through the relation
 \begin{eqnarray}\label{delta1a} (1-e^{it})\p_X(t)&=&\sum_{m\in \Z} \frac{(itm)}{m!}\big({\mathbb P}\{X=m\}-{\mathbb P}\{X=m-1\}\big).
 \end{eqnarray}
Hence
 \begin{eqnarray}\label{delta1b} |\p_X(t)|&\le & \frac{\d_X}{2|\sin (t/2)|  }\qq\quad (t\notin 2\pi\Z).
 \end{eqnarray}
The Lemma below (\cite{G1}, Lemma 2) is obtained as    a consequence of inequality \eqref{delta1b} and Cram\'er's inequality
{\color{red}(\cite{Cr70}, Lemma\,1, p.27,)} which we recall.
\vskip 3 pt
{\color{blue}
{\color{green}
 {\bf Cram\'er's inequality:}} If $f(t)$ is the characteristic function of a distribution function and if, for some $b>0$,
$$ |f(t)|\le\a, \qq \hbox{for ${\color{red}2b}>|t|\ge b$,}$$
where $0<\a<1$, then
\beq\label{cramer.ineq}
|f(t)|\le 1-(1-\a^2)(|t|^2/8b^2)
\eeq
for $|t|<b$.

}
\begin{lemma}\label{L2G1} Let  $\eta_i$, $i\ge 1$,  be independent copies of a random variable $\eta$ and let $T_n= \eta_1 + \ldots+ \eta_n$ for each $n$.  Let $n_0$ be
chosen so that $\d_{T_{n_0}}<\sqrt 2$. Let     $f(t )$ denote the characteristic function of $\eta$.
Then for all $|t|\le \pi$,
$$ |f(t )|\le e^{-ct^2}, $$
where
$$ c=\frac{1}{2\pi^2 n_0}\,\Big(1-\frac{\d^2_{T_{n_0}}}{2}\Big).$$
\end{lemma}
{\color{green}
\begin{remark}    The existence of an integer $n_0$ such that $\d_{T_{n_0}}<\sqrt 2$ is possible only if the span of $\eta$ is one (Proposition \ref{deltaX.prop}).
\end{remark}
}
 \vskip 3 pt     Lemma \ref{L2G1} is used in Gamkrelidze {\color{green}\cite{G1}}
        to prove the following  result.
   }

 \begin{theorem}
  \label{gamkrelidze} Let $  X$ be  a sequence
  of independent integer valued random variables satisfying the following conditions:

\begin{itemize}
\leftskip6mm
\itemsep=1pt
\item[{\rm (i)}\ ]     There exists a positive real $\l<  \sqrt 2$ and an integer $n_0$ such that
for all $k$, $  \d_{X_k^1+\ldots +X_k^{n_0}} \le \l$,  where $X_k^{ j}, 1\le j\le n_0$ are independent copies of $X_k$.

\item[{\rm (ii)}\ ]   The central limit theorem is applicable.

\item[{\rm (iii)}\ ]  ${\rm Var}(S_n)=\mathcal O(n)$.
\end{itemize}
Then the local limit theorem is applicable in the strong form.
\end{theorem}
{\color{green}
\begin{remark} The sufficient conditions obtained in this statement, as well as in some others, do not require any explicit speed of convergence in the central limit theorem,
which is a remarkable characteristic of the method used.
\end{remark}
}
{\color{blue}
  Let us explain how to bound from above the main  integral term ($I_4$ in \eqref{basic.bound}), by using Lemma \ref{L2G1}. Let $f(t,S_n)$, $f(t,X_j)$ respectively denote the
  characteristic function of $S_n$ and of $X_j$. {\color{green} Let $B_n^2= {\rm Var}(S_n)$. }Assumption (i) implies that
$$|f(\frac{t}{B_n},S_n )|\le \prod_{j=1}^n |f(\frac{t}{B_n},X_j )|\le e^{-cn(t/B_n)^2} $$
and one can take  $ c=\frac{1}{2\pi^2 n_0}\,(1- {\l^2 }/{2})$.
Let $\e >0$. Then  {\color{green}by assumption (iii)},
\ben\label{I3.G}
\int_{\e \le |t |\le \pi B_n}|f(\frac{t}{B_n},S_n )| \dd t &\le& \int_{\e \le |t |\le \pi B_n}e^{-cn(t/B_n)^2} \dd t
\cr &\le &
B_n\int_{\e/B_n \le |u|\le \pi }e^{-cnu^2} \dd u
\cr &\le &
 \frac{B_n}{\sqrt {2cn}}\int_{|v|\ge\sqrt {2cn} \e/B_n  }e^{-v^2/2} \dd v
\cr &\le &  \frac{B_n}{\sqrt {2cn}}e^{-cn\e^2/B_n^2}\,\le \,  C\,e^{-C'\e^2 }
\een
{\color{green}which is small for $\e$ large.}

also obtained under which a local limit theorem follows from an integral limit theorem in some form. However, as a rule, these conditions are hard to verify in the case of
summands that are not identically distributed.
  \vskip 10 pt

  {(2)}    Consider the following characteristic
\begin{eqnarray}\label{vartheta}  \t_X =\sum_{m\in \Z}
{\mathbb P}\{X=m\}\wedge{\mathbb P}\{X=m+1\} ,
\end{eqnarray}
 where $a\wedge b=\min(a,b)$. Note    that
   \begin{eqnarray}\label{vartheta1}
0\le \t_X<1 .\end{eqnarray}

  Indeed, let   $k_0$ be some integer such that
$  {\mathbb P}\{X=m_0\} >0$. Then
$$\sum_{m=m_0}^\infty  {\mathbb P}\{X=m\}\wedge{\mathbb P}\{X=m+1\}\le \sum_{m=m_0}^\infty  {\mathbb P}\{X=m+1\}=\sum_{m=m_0+1}^\infty   {\mathbb P}\{X=m \},
$$
and so  $0\le  \t_X\le       \sum_{m< m_0 }  {\mathbb P}\{X=m+1\} +\sum_{m=m_0+1}^\infty  {\mathbb P}\{X=m+1\}<1$.

\vskip 5 pt
{\color{blue} Note also   that
 \begin{equation}\label{corres} \sum_{m\in\Z} \Big(\P\{X=m \}+ \P\{X=m +1\}-\big|\P\{X=m \}+ \P\{X=m +1\}\big|\Big)=2\t_X,
 \end{equation}
since for positive reals $a$ and $b$, one has $a+b-|a-b|= 2(a\wedge b)$. But  the left-hand side is just $2 -\d_X$. Whence the relation
\beq\label{vartheta2}
\d_X =2(1- \t_X) ,
\eeq
which was quoted in  Mukhin
 \cite{Mu}, p.\,700.    We thus have the equivalence
 \begin{equation}\label{equiv}  \t_X>0\quad \Longleftrightarrow \quad \d_X<2.
\end{equation} \smallskip\par
An important  consequence is
\begin{corollary}
The condition $\d_X<2$ implies that $X$ has a Bernoulli component. \end{corollary}
 This follows from  Lemma \ref{bpr}. We add the interesting remark below.
\begin{remark} \label{deux.char}
  The  two characteristics seem through formula \eqref{vartheta2} equivalent. However by formula \eqref{delta1a}, Lemma \ref{L2G1} and Theorem \ref{gamkrelidze},
 the first is related to the method of characteristic  functions, whereas the second is related to the Bernoulli part extraction, which is known to be characteristic function
 free, and is used in the proof of the next result,  see also prerequisites to Theorem \ref{ger1}.
\end{remark}
}

\vskip 5 pt

 {\color{green} Davis and McDonald
     proved in \cite{MD} a variant of  Gamkrelidze's Theorem \ref{gamkrelidze},  based on   the Bernoulli part
extraction of a random variable. The   characteristic $\t_X$ is used in \cite{MD},  and it is required that $\t_X>0$. This method will be  explained in details later.
Davis and McDonald's theorem \cite[Th.\,1.1]{MD}   states as follows.
}

 \begin{theorem}\label{dmd.th}  Let $\{ X_j , j\ge 1\}$ be independent, integer valued random variables with partial sums
$S_n= X_1+\ldots +X_n$ and let
$f_j(k)=
{\mathbb P}\{X_j=k\}$.  Also for each $j$ and $n$, let
$$q(f_j)= \sum_{k} [f_j(k)\wedge f_j(k+1)], \qq Q_n=\sum_{j=1}^n q(f_j)  $$
and assume that $q(f_j)>0$ for each $j\ge 1$. Further assume that there exist numbers $b_n>0$, $a_n$  such that
 $$\lim_{n\to \infty}b_n= \infty, \qq \limsup_{n\to \infty} \ {b_n^2}/{Q_n}<\infty,$$
and
$$ \frac{S_n-a_n}{b_n} \ \ \buildrel{\mathcal L}\over{\Longrightarrow}\ \ \mathcal N(0,1).$$
  Then
$$ \lim_{n\to \infty} \sup_{k} \Big|b_n {\mathbb P}\{S_n=k\} -\frac{1}{\sqrt{2\pi}}e^{-  \frac{(k-a_n)^2}{2b_n^2}}\Big|=0. $$
 \end{theorem}

{\color{green} It is implicitely assumed that ${\mathbb P}\{X_j=v_k\}>0$ and ${\mathbb P}\{X_j=v_{k+1}\}>0$ for some $k\in \Z$.
This is a serious restriction,   a large class of examples or models being   built with random variables which are  precisely of  this type. This is however not a restriction
for the study of the almost sure local limit theorem.
 }

  \begin{remark} \label{ssprime1}Assume that $X$ has finite mean $\mu$ and finite variance
$\sigma^2$. Then we have the following inequality
 \begin{eqnarray}
\s^2\ge  \frac{ D^2   }{4} \t_X.
\end{eqnarray}
Indeed,
\begin{align*}
&{\mathbb P}\Big\{|X - \mu|\ge  \frac{D} {2} \Big\}  =\sum_{v_k\ge\mu+ {{D}\over{2}}
}{\mathbb P}\{X=v_k\}+ \sum_{v_k\le \mu- {{D}\over{2}} }{\mathbb P}\{X=v_k\}
\cr
&\ge
\!\!\sum_{v_{k+1}\ge\mu+ {{D}\over{2}} }{\mathbb P}\{X=v_{k})\wedge
{\mathbb P}\{X=v_{k+1}\}  + \!\!\!  \sum_{v_k\le \mu- {{D}\over{2}} }{\mathbb P}\{X=v_k\}\wedge
{\mathbb P}\{X=v_{k+1}\} \cr
&=  \sum_{v_{k}\ge\mu- {{D}\over{2}}
}{\mathbb P}\{X=v_{k}\}\wedge {\mathbb P}\{X=v_{k+1}\}  +\!\!\!  \sum_{v_k\le \mu-
{{D}\over{2}} }{\mathbb P}\{X=v_k\}\wedge {\mathbb P}\{X=v_{k+1}\}
\cr&\ge
 \t_X.
\end{align*}
Next by
 Tchebycheff's inequality,
${{D^2}\over{4}}{\mathbb P}\big\{|X - \mu|\ge  {{D}\over{2}}\big\}\le  \sigma^2$.
\end{remark}

 \vskip 15 pt

 {(3)} {\color{blue}The following characteristic for which we refer to    Mukhin \cite{Mu}, was already studied in earlier investigations on the local limit theorem.}
   \begin{eqnarray}\label{ssprime2} \mathcal D(X,d):=\inf_{a\in \R}{\mathbb E\,} \langle (X-a)d\rangle^2 .
 \end{eqnarray}
Here $d$ is a real number, $|d|\le 1/2$ and $\langle \a \rangle$ is the distance from $\a$ to the nearest integer.
Notice that $ \mathcal D(X,d)=0$ if and only if $X$ is lattice valued with span $1/d$.
The link with $\t_X$ is given in  Lemma 2 of  \cite{Mu},
\begin{eqnarray}\label{ssprime3} \mathcal D(X,d) \ge  \frac{ |d|^2   }{4} \t_X.
 \end{eqnarray}

{\color{blue}This characteristic is also close to those arising from Prohorov and Rozanov's Theorems, (see Theorem \ref{rozanov.I}, see also \cite{PrRo}, \cite{Ro})
\begin{eqnarray}\label{char.Pro.Roz}
\nu(X,h)=\min_{j,0,\ldots, h-1} \P\{X\not\equiv j\, ({\rm mod} \,h) ,  \end{eqnarray}
where $h\ge 2$. We have (\cite{Mu}, p.\,700),
\begin{eqnarray}\label{comp.char.Mukh.Pro.Roz}
\frac{1}{2h^{3}}\nu(X,h)\le \mathcal D(X,\frac1h)\le \frac14\nu(X,h).\end{eqnarray}}
\vskip 3 pt

Consider also the   characteristic
$$ H(X ,d) = {\mathbb E\,} \langle X^*d\rangle^2,$$  where
$X^*$
is  a symmetrization of $X$.
   In Mukhin \cite{Mu} and \cite{Mu1},
the two-sided inequality
 \begin{eqnarray}\label{fih} 1-2\pi^2 H(X ,\frac{t }{2\pi})  \le |\p_X(t)|\le 1-4  H(X ,\frac{t }{2\pi})   ,
\end{eqnarray}
is established.

{\color{blue}Mitalauskas and Statulyavichus used  in \cite{MiSt} and in other papers  the following refined characteristic
\begin{eqnarray}\label{char.Mita.Statu}
\a(X,a,q,M)= q^{-2}\sum_{r\in(-q/2,q/2]}r^2\P\{aX^*\equiv r\, ({\rm mod} \,q),\, |X^*|\le M\},\end{eqnarray}
where $a,q$ are positive coprime numbers, $(a,q)=1$, $a\le q/2$. \smallskip\par
For the shown values of $a$ and $q$, the following relation holds
 \begin{eqnarray}\label{comp.char.Mukh.Mita.Statu} H(X, a/q)=\a(X,a,q,\infty).
\end{eqnarray}
 }The following   is  the one-dimensional version of Theorem 5 in \cite{Mu}  and is stated without proof, however.

\begin{theorem}[Mukhin]\label{Mukhin.th.Hn} Let $X_1,\ldots, X_n$ have zero mean and finite third moments. Let
$$ B_n^2= \sum_{j=1}^n{\mathbb E\,} |X_j|^2 ,\qq H_n= \inf_{1/4\le d\le 1/2}\sum_{j=1}^n H(X_j
,d), \qq L_n= \frac{\sum_{j=1}^n{\mathbb E\,} |X_j|^3}{(B_n)^{3/2}} .$$ Then
 $\D_n\le CL_n\, \big( {B_n }/{ H_n}\big)$, recalling that $\D_n$ is  defined in Definition {\color{red} \ref{defllt}}. \end{theorem}

{\color{green}
\begin{remark} Under the assumption made, the central limit theorem is fulfilled with speed of convergence $L_n$, by Berry-Esseen inequality. However no condition is made on
$L_n$.
\end{remark}
}
     A manuscript
devoted to  estimates of the rate of convergence was announced in \cite{Mu} with  no resulting
publication, however.

\vskip 15 pt

\noi{\tf The Bernoulli part of a random variable.}
     This   is one instance of coupling method. A remarquable feature of this one is  that it is  \lq\lq characteristic  function
free\rq\rq.     This approach, due to McDonald \cite{MD}, is indeed purely probabilistic and  consists with extracting the Bernoulli
part of a random variable.       It seems, according to
\cite{AGKW}  that this idea first appeared in   Kolmogorov's work  \cite{KO2} on concentration functions.
  \vskip 1 pt  It is worth citing here  Kolmogorov (from his 1958's paper \cite{KO2} p.\ 29): \lq\lq ...
\!\!{\it Il semble cependant que nous restons toujours  dans une p\'eriode o\`u la comp\'etition de ces deux directions {\rm [characteristic functions or
direct methods from the calculus of probability]} conduit aux r\'esultats les plus f\'econds ...\rq\rq }.
 \smallskip\par

Dabrowski and McDonald   mentioned in \cite{DMD} that this idea is implicit in Galstyan \cite{Gal}, Mineka \cite{Min} and R\"osler \cite{Ros}.

 \smallskip\par
  Let $X$ be
  a random variable   such that  ${\mathbb P}\{X
\in\mathcal L(v_0,D)\}=1$. We do not   assume  that the  span $D$ is maximal. Further, {\it no} integrability condition on $X$ is necessary.    Put
$$ f(k)= {\mathbb P}\{X= v_k\}, \qq k\in \Z .$$
We assume that
 \begin{equation}\label{basber1}\t_X>0,
  \end{equation}
where $\t_X$ is defined in \eqref{vartheta}. Recall that $
\t_X<1$.
\vskip 5 pt
Let $0<\t\le\t_X$. One can associate to $\t$ and $X$  a
sequence $  \{ \tau_k, k\in \Z\}$     of   non-negative reals such that
\begin{equation}\label{basber0}  \tau_{k-1}+\tau_k\le 2f(k), \qq  \qq\sum_{k\in \Z}  \tau_k =\t.
\end{equation}
Just take
 $\tau_k=  \frac{\t}{\nu_X} \, (f(k)\wedge f(k+1))  $.
   Now   define   a pair of random variables $(V,\e)$   as follows:
  \begin{eqnarray}\label{ve} \qq\qq\begin{cases} {\mathbb P}\{ (V,\e)=( v_k,1)\}=\tau_k,      \cr
 {\mathbb P}\{ (V,\e)=( v_k,0)\}=f(k) -{\tau_{k-1}+\tau_k\over
2}    .  \end{cases}\qq (\forall k\in \Z)
\end{eqnarray}
   By assumption     this is   well-defined,   and the margin  laws verify
\begin{eqnarray}\begin{cases}{\mathbb P}\{ V=v_k\} &= \  f(k)+ {\tau_{k }-\tau_{k-1}\over 2} ,
\cr
 {\mathbb P}\{ \e=1\} &= \ {\vartheta} \ =\ 1-{\mathbb P}\{ \e=0\}   .
\end{cases}\end{eqnarray}
Indeed, ${\mathbb P}\{ V=v_k\}= {\mathbb P}\{ (V,\e)=( v_k,1)\}+ {\mathbb P}\{ (V,\e)=( v_k,0)\}=f(k)+ {\tau_{k }-\tau_{k-1}\over 2} .$
 Further  ${\mathbb P}\{ \e=1\}  =\sum_{k\in\Z} {\mathbb P}\{ (V,\e)=( v_k,1)\}=\sum_{k\in\Z} \tau_{k }={\vartheta}  $.

 \begin{lemma} \label{bpr} Let $L$
be a Bernoulli random variable    which is independent of  $(V,\e)$, and put
 $Z= V+ \e DL$.
We have $Z\buildrel{\mathcal D}\over{ =}X$.
\end{lemma}

\begin{proof}[Proof]
 Plainly,
\begin{eqnarray*}{\mathbb P}\{Z=v_k\}&=&{\mathbb P}\big\{ V+\e DL=v_k, \e=1\}+ {\mathbb P}\big\{ V+\e DL=v_k, \e=0\} \cr
&=&{{\mathbb P}\{ V=v_{k-1}, \e=1\}+{\mathbb P}\{
V=v_k, \e=1\}\over 2} +{\mathbb P}\{ V=v_k, \e=0\}
\cr&=& {\tau_{k-1}+ \tau_{k }\over 2} +f(k)-{\tau_{k-1}+ \tau_{k
}\over 2}
= f(k).
 \end{eqnarray*}
\end{proof}

Consider now independent random variables  $ X_j,j=1,\ldots,n$,    and assume that
\begin{equation}\label{basber.j}  \t_{X_j}>0, \qq \quad  j=1,\ldots, n.
\end{equation}
Let $0<\t_j\le \t_{X_i}$, $j=1,\ldots, n$. Iterated  applications of Lemma \ref{bpr} allow us to
 associate to them a
sequence of independent vectors $ (V_j,\e_j, L_j) $,   $j=1,\ldots,n$  such that
 \begin{eqnarray}\label{dec0} \big\{V_j+\e_jD   L_j,j=1,\ldots,n\big\}&\buildrel{\mathcal D}\over{ =}&\big\{X_j, j=1,\ldots,n\big\}  .
\end{eqnarray}

Further the sequences $\{(V_j,\e_j),j=1,\ldots,n\}
 $   and $\{L_j, j=1,\ldots,n\}$ are independent.
For each $j=1,\ldots,n$, the law of $(V_j,\e_j)$ is defined according to (\ref{ve}) with $\t=\t_j$.  And $\{L_j, j=1,\ldots,n\}$ is  a sequence  of
independent Bernoulli random variables. Set
\begin{equation}\label{dec} S_n =\sum_{j=1}^n X_j, \qq  W_n =\sum_{j=1}^n V_j,\qq M_n=\sum_{j=1}^n  \e_jL_j,  \quad B_n=\sum_{j=1}^n
 \e_j .
\end{equation}

\begin{proposition}[\it Decomposition of sums] \label{lmd}We have\begin{eqnarray*} \{S_k, 1\le k\le n\}&\buildrel{\mathcal D}\over{ =}&  \{ W_k  +  DM_k, 1\le k\le n\} .
\end{eqnarray*}
And  $M_n\buildrel{\mathcal D}\over{ =}\sum_{j=1}^{B_n } L_j$.
 \end{proposition}

\vskip 15 pt
  \noi {\tf A local limit theorem  with effective rate.}
 A careful application of the Bernoulli extraction part method allows one to obtain not only a local limit theorem, but in addition an effective rate of convergence. This is
 done in the recent work of Giuliano-Weber  \cite{GW3}. The proof is  completely elementary and the estimates obtained are quite simple and therefore easy to apply.
\smallskip\par    Let   $ X_j,j=1,\ldots,n$ be independent square integrable random variables  taking almost surely values in a common lattice   $\mathcal L(v_{
0},D )=\{v_k,k\in \Z\}$, where $v_{ k}=v_{ 0}+D k$, $k\in \Z$,   $v_{0} $ and $D >0$ are   real numbers. Let $S_n=X_1+\ldots +X_n$. Assume that condition \eqref{basber.j}  is
fulfilled.

\vskip 3 pt
Let also $0<\t_j\le \t_{X_i}$, $j=1,\ldots, n$. By Proposition \ref{lmd}, we have
$${\mathbb P} \{S_n =\kappa \}   =     \E_{(V,\e)}     {\mathbb P}_{\!L}
\big\{D M_n  =\kappa-W_n
\big\} $$
where $S_n'= \E_L S_n= W_n + \frac{D}{2} B_n $, and  $\E_{\!L}$, ${\mathbb P}_{\!L}$ (resp.\ $\E_{(V,\e)}$, ${\mathbb P}_{(V,\e)}$) stand  for the integration
symbols   and probability symbols relatively to the
$\s$-algebra generated by the sequence  $\{L_j , j=1, \ldots, n\}$   (resp.\ $\{(V_j,\e_j), j=1, \ldots, n\}$).
 Further ${\mathbb E\,} S'_n= {\mathbb E\,} S_n$, ${\mathbb E\,}  (S'_n) ^2=
 {\mathbb E\,}  S_n ^2- \frac{D^2\Theta_n }{4}$.
  Set
\begin{eqnarray*} \begin{cases}H_n&=  \ \ \sup_{x\in \R} \big|{\mathbb P}_{(V,\e)}
\big\{{S'_n-\E_{(V,\e)}S'_n  \over \sqrt{{\rm Var}(S'_n)}  }<x\big\} - {\mathbb P}\{g<x\} \big|,\cr
  \rho_n(h)&= \ \ {\mathbb P}\big\{\big|\sum_{j=1}^n \e_j-\Theta_n\big|>h\Theta_n\big\},   \cr
  \Theta_n&=\ \ \sum_{j=1}^n \t_j.
  \end{cases}
\end{eqnarray*}
\vskip 5 pt

\begin{theorem}[\cite{GW3},\,Th.\,1.1]
\label{ger1}
  For any  $0<h<1$, $0<\t_j\le \t_{X_j}$,
 and all   $\kappa \in \mathcal L( v_{
0}n,D )$
 \begin{eqnarray*}      {\mathbb P}\{S_n =\kappa \} &\le  &
      \Big( \frac{1+ h  }{ 1-h}\Big) \, {  D \over
\sqrt{2  \pi  {\rm Var}(S_n)  } }  e^{-\frac{(\kappa - {\mathbb E\,} S_n)^2}{  2(1  +  h){\rm Var}(S_n)   } }
 \cr & &\quad   + {C_1
 \over
\sqrt{   (1-h)\Theta_n} } \big(H_n +      \frac{1}{(1-h)\Theta_n}
 \big) + \rho_n(h)     ,\end{eqnarray*}
where $ C_1= \max ({ 8 /\sqrt{2  \pi  } },C_0  )  $ and $C_0$ is the same constant as in Theorem \ref{lltber}.
 \vskip   3 pt Further,
\begin{eqnarray*}
{\mathbb P} \{S_n =\kappa \}   &\ge &   \Big(\frac{   {      1- h   }}{       {  1 +h      }}\Big)  { D \over
\sqrt{2\pi {\rm Var}(S_n)     }}
 {  e^{- \frac{(\kappa -
 {\mathbb E\,} S_n)^2}{
 2(1-h){\rm Var}(S_n)    } }}
   \cr
   & &\ - \   {C_1
 \over
\sqrt{   (1-h)\Theta_n} } \Big(H_n  +   \frac{1}{(1-h)\Theta_n}  +  2\rho_n(h)
 \Big) - \rho_n(h).
  \end{eqnarray*}
 \end{theorem}

We shall hereafter deduce from Theorem \ref{ger1} the following corollary.

\begin{corollary}  \label{ger2}  Assume that $\frac{  \log \Theta_n }{\Theta_n}\le  {1}/{14} $. Then, for all $\kappa \in \mathcal L( v_{
0}n,D )$ such that
$$\frac{(\kappa - {\mathbb E\,} S_n)^2}{    {\rm Var}(S_n)  } \le   \big({\frac { \Theta_n} {14 \log \Theta_n} }\big)^{1/2}   ,$$  we have
\begin{eqnarray*}  \Big| {\mathbb P} \{S_n =\kappa \} -{ D    e^{- \frac{(\kappa - {\mathbb E\,} S_n)^2}{    2 {\rm Var}(S_n)    } }  \over \sqrt{2\pi {\rm Var}(S_n)     }}
\Big|  & \le &   C_2\Big\{ D\big({    {   \log
\Theta_n }     \over
 {    {\rm Var}(S_n)   \Theta_n} } \big)^{1/2}  +    {    H_n +  \Theta_n^{-1}
\over \sqrt{   \Theta_n} } \Big\} .
   \end{eqnarray*}
Here   $C_2=   2^{7/2} \, C_1 $.
   \end{corollary}

 \begin{proof}[Proof of Theorem \ref{ger1}.]

 We   denote again $X_j= V_j+D\e_jL_j$, $S_n= W_n  +  M_n$,
$j, n\ge 1$.
      Fix $0<h<1$ and let
  \begin{eqnarray}\label{dep00}A_n=\Big\{\big| B_n - \Theta_n\big|\le h\Theta_n
\Big\}, \qq\qq  \rho_n(h)= {\mathbb P}_{(V,\e)}(A_n^c)   .
\end{eqnarray}
For $\kappa \in \mathcal L(v_0,D)$,
  \begin{equation}
    {\mathbb P} \{S_n =\kappa \}
 \ =\ \E_{(V,\e)}   \Big( \chi(A_n)+\chi(A_n^c)\Big)    {\mathbb P}_{\!L}
\Bigl\{D \sum_{j= 1}^n \e_jL_j  =\kappa-W_n
\Bigr\} ,\end{equation}
and so
 \begin{equation}\label{dep0} \biggl|{\mathbb P} \{S_n =\kappa \}  -  \E_{(V,\e)}   \chi(A_n)    {\mathbb P}_{\!L}
\Bigl\{D \sum_{j= 1}^n \e_jL_j  =\kappa-W_n
\Bigr\}\biggr|\, \le \,   {\mathbb P}_{(V,\e)}(A_n^c)\, = \,  \rho_n(h)
. \end{equation}

  By  Lemma \ref{lltber}, since $\sum_{j= 1}^n \e_jL_j\buildrel{\mathcal D}\over{ =}\sum_{j=1}^{B_n } L_j$,
  \begin{eqnarray*}  \sup_{z}\, \Big|  {\mathbb P}_{\!L}\big\{\sum_{j=1}^{N } L_j=z  \big\} -{2\over \sqrt{2\pi
N}}e^{-{( z-(N/2))^2\over
 N/2}}\Big|
 \le {C_0\over N^{3/2}} .
 \end{eqnarray*}
    On $A_n$,      $   (1-h )\Theta_n \le  B_n   \le     (1+h )\Theta_n$,   so that
    \begin{eqnarray} \label{dep2} \Big|\E_{(V,\e)}     \chi(A_n) \Big\{     {\mathbb P}_{\!L}
\big\{D \sum_{j= 1}^n \e_jL_j  =\kappa-W_n
\big\} - {2e^{-{(\kappa-W_n-D(B_n/2))^2\over
D^2(B_n/2)}}\over
\sqrt{2\pi B_n}} \Big\}\Big|\cr  \le   C_0\ \E_{(V,\e)}  \chi(A_n)\cdot  B_n^{-3/2}\le     \frac{ C_0}{    (1-h )^{ 3/2}} \,\frac{
1}{(\sum_{i=1}^n \t_i )^{ 3/2}
 }
  .
\end{eqnarray}
Inserting this in  (\ref{dep0}), we get
 \begin{equation}\label{dep01}  \Big|{\mathbb P} \{S_n =\kappa \}  -  \E_{(V,\e)}   \chi(A_n)   {2e^{-{(\kappa-W_n-D(B_n/2))^2\over
D^2(B_n/2)}}\over
\sqrt{2\pi B_n}}  \Big|  \ \le \  \frac{ C_0}{     (1-h )^{ 3/2}(\sum_{i=1}^n \t_i )^{ 3/2}
 }
 + \rho_n(h)
.
  \end{equation}

       Then
\begin{equation}\label{prime} {\rm Var}(S'_n)={\rm Var}(S_n)-
\frac{D^2}{4}\sum_{i=1}^n \t_{ i}= \sum_{i=1}^n\Big(\s_i^2-\frac{D ^2\t_{ i}}{4}\Big).
\end{equation}

 Put $T_n= {S'_n-\E_{(V,\e)}S'_n  \over \sqrt{{\rm Var}(S'_n)}  }$.
 As $\E_{(V,\e)}S'_n={\mathbb E\,} S_n$, we note that
\begin{eqnarray*}  {(\kappa-W_n- D(B_n/2))^2\over  D^2(B_n/2)}
 & = &{ {\rm
Var}(S'_n) \over  D^2 (B_n/2) }\Big({\kappa-{\mathbb E\,} S_n \over \sqrt{{\rm Var}(S'_n)}  } - T_n \Big)^2,
\end{eqnarray*}
and   rewrite (\ref{dep01}) as follows,
\begin{eqnarray}\label{dep21} \Big|{\mathbb P} \{S_n =\kappa \}  -  \Upsilon_n\Big|&\le &  \frac{ C_0}{     (1-h )^{ 3/2}} \,\frac{
1}{\Theta_n^{ 3/2}
 }
 + \rho_n(h),
  \end{eqnarray}
where
\begin{eqnarray}\label{dep210}  \Upsilon_n=  \E_{(V,\e)}   \chi(A_n) {2e^{-{
 {\rm Var}(S'_n) \over  D^2 (B_n/2) }\big({\kappa-{\mathbb E\,} S_n \over \sqrt{{\rm Var}(S'_n)}  } - T_n\big)^2}\over
\sqrt{2\pi B_n}}
.
 \end{eqnarray}
  Set for $-1< u\le 1$,
 $$ Z_n(u)=   \E_{(V,\e)}       e^{-{2{\rm Var}(S'_n)\over  D^2  (1 + u ) \Theta_n   } \big({\kappa-{\mathbb E\,} S_n
\over
\sqrt{{\rm Var}(S'_n)}  }-  T_n\big)^2}.  $$
Then
\begin{eqnarray} \label{dep3}      { 2Z_n(-h) - 2\rho_n(h) \over
\sqrt{2\pi    (1 +h )\Theta_n}}\ \le \   \Upsilon_n  &\le &{ 2Z_n(h) \over
\sqrt{2\pi    (1 - h )  \Theta_n}} .
 \end{eqnarray}

  Let $Y  $ be a centered
  random variable.  By the transfer formula, for any positive reals
$a$ and
$b$, we have
\begin{eqnarray}\label{tech} \Big|{\mathbb E\,}  e^{-a(b-Y)^2} - \frac{1}{ \sqrt{1+2a}}\, e^{-  \frac{b^2}{2+ 1/a} }\Big|&\le &  4\sup_{x\in \R} \big|{\mathbb P}
 \{Y<x \} - {\mathbb P}\{g<x\} \big| .
\end{eqnarray}
Applying this to $Z_n(u)$ with
 $ a= {2{\rm
Var}(S'_n)\over D^2 (1+u ) \Theta_n}$, $b=
\frac{\kappa-  {\mathbb E\,} S_n}{\sqrt{{\rm Var}(S'_n)}}  $ gives
  with \eqref{prime},
\begin{eqnarray*} \frac{b^2}{2+ 1/a} &= &\frac{(\kappa-  {\mathbb E\,} S_n)^2}{ {\rm Var}(S'_n)\big( 2+  \frac{D^2 (1+u )\Theta_n }{2{\rm
Var}(S'_n)}\big)}  \,=\,\frac{(\kappa-  {\mathbb E\,} S_n)^2}{   2{\rm Var}(S_n)(1  +  \d(u)) },  \end{eqnarray*}
where
$\d(u)= \frac{D^2      \Theta_n  u}{4 {\rm Var}(S_n) } $.
 Further \begin{eqnarray*}\frac{1}{ \sqrt{1+2a}}& =& \Big(\frac{1}{  {1+ {4{\rm Var}(S'_n)\over D^2
  (1+u ) \Theta_n}}}\Big)^{1/2}  \,=\,
 \frac{D}{ 2    }\Big(\frac{   { \Theta_n   (1+ u  )}}{       {{\rm Var}(S_n)  (1+ \d(u)  ) }}\Big)^{1/2} .
 \end{eqnarray*}

This along with   \eqref{tech}   provides the bound,
   \begin{eqnarray} \label{fb}    \Big| Z_n(u) -   \frac{D}{ 2    }\Big(\frac{   { \Theta_n   (1+ u  )}}{       {{\rm Var}(S_n)  (1+
\d(u)  ) }}\Big)^{1/2} e^{- \frac{(\kappa-  {\mathbb E\,} S_n)^2}{   2{\rm Var}(S_n)(1  +  \d(u)) } }  \Big| &\le  &    4H_n , \end{eqnarray}
 with
$$ H_n=\sup_{x\in \R} \Big|{\mathbb P}_{(V,\e)}
\big\{T_n<x\big\} - {\mathbb P}\{g<x\} \Big| .$$
Besides, it follows from (\ref{ssprime1}) that
 $ 0\le \d(h) \le h$, for $h\ge 0$.
   By reporting (\ref{fb}) into (\ref{dep3})  we get,
   \begin{eqnarray*}   \Upsilon_n  &\le &
  { 8H_n   \over \sqrt{2  \pi  (1-h)\Theta_n} }      +   \Big( \frac{1+ h  }{ 1-h}\Big) \, {  D \over
\sqrt{2  \pi  {\rm Var}(S_n)  } }  e^{-\frac{(\kappa-  {\mathbb E\,} S_n)^2}{  2(1  +  h){\rm Var}(S_n)   } }
  .
\end{eqnarray*}
  And by combining   with (\ref{dep21}),
  \begin{eqnarray} \label{dep5}    {\mathbb P}\{S_n =\k\} &\le  &
      \Big( \frac{1+ h  }{ 1-h}\Big) \, {  D \over
\sqrt{2  \pi  {\rm Var}(S_n)  } }  e^{-\frac{(\kappa-  {\mathbb E\,} S_n)^2}{  2(1  +  h){\rm Var}(S_n)   } }
 \cr & &\quad   + { 8H_n
 \over
\sqrt{2  \pi  (1-h)\Theta_n} } +   \frac{ C_0}{     (1-h )^{ 3/2}} \,\frac{ 1}{ \Theta_n  ^{ 3/2}
 }
 + \rho_n(h)     .
\end{eqnarray}
     Similarly  by using
(\ref{dep3}),
\begin{eqnarray} \label{dep6}
 \Upsilon_n
  &\ge  &
\Big(\frac{   {      1- h   }}{       {  1 +h      }}\Big)  { D \over
\sqrt{2\pi {\rm Var}(S_n)     }}
 {  e^{- \frac{(\kappa-
 {\mathbb E\,} S_n)^2}{
 2(1-h){\rm Var}(S_n)    } }}
   - {8  H_n+ 2\rho_n(h) \over
\sqrt{2\pi    (1 +h )\Theta_n}}  .
 \end{eqnarray}

By combining  with  (\ref{dep21}), we obtain
\begin{eqnarray}\label{dep5a}  {\mathbb P} \{S_n =\kappa \}   &\ge &   \Big(\frac{   {      1- h   }}{       {  1 +h      }}\Big)  { D \over
\sqrt{2\pi {\rm Var}(S_n)     }}
 {  e^{- \frac{(\kappa-
 {\mathbb E\,} S_n)^2}{
 2(1-h){\rm Var}(S_n)    } }}
   - { 8  H_n+ 2\rho_n(h) \over
\sqrt{2\pi    (1 +h )\Theta_n}}\cr & &\ -\frac{ C_0}{     (1-h )^{ 3/2}} \,\frac{
1}{\Theta_n^{ 3/2}
 }
 - \rho_n(h),
  \end{eqnarray}
 As  $ C_1= \max ({ 8
 /\sqrt{2  \pi  } },C_0  )  $, we deduce
\begin{eqnarray} \label{f1}    {\mathbb P}\{S_n =\k\} &\le  &
      \Big( \frac{1+ h  }{ 1-h}\Big) \, {  D \over
\sqrt{2  \pi  {\rm Var}(S_n)  } }  e^{-\frac{(\kappa-  {\mathbb E\,} S_n)^2}{  2(1  +  h){\rm Var}(S_n)   } }
 \cr & &\quad   + {C_1
 \over
\sqrt{   (1-h)\Theta_n} } \big(H_n +      \frac{1}{(1-h)\Theta_n}
 \big) + \rho_n(h)     .
\end{eqnarray}
And
\begin{eqnarray}\label{f2}  {\mathbb P} \{S_n =\kappa \}   &\ge &   \Big(\frac{   {      1- h   }}{       {  1 +h      }}\Big)  { D \over
\sqrt{2\pi {\rm Var}(S_n)     }}
 {  e^{- \frac{(\kappa-
 {\mathbb E\,} S_n)^2}{
 2(1-h){\rm Var}(S_n)    } }}
   - \cr & &\
  {C_1
 \over
\sqrt{   (1-h)\Theta_n} } \big(H_n  +   \frac{1}{(1-h)\Theta_n}  +  2\rho_n(h)
 \big) - \rho_n(h).
  \end{eqnarray}
This achieves the proof.
\end{proof}

\begin{proof}[Proof of Corollary \ref{ger2}]

 In order to estimate $\rho_n(h)$ we use the following Lemma  (\cite{di}, Theorem 2.3)
\begin{lemma} \label{primo}
Let $X_1, \dots, X_n$     be independent random variables, with $0 \le X_k \le 1$ for each $k$.
Let $S_n = \sum_{k=1}^n X_k$ and $\mu = {\mathbb E\,} S_n$. Then for any $\e   >0$,
 \begin{eqnarray*}
(a) &&
 {\mathbb P}\big\{S_n \ge  (1+\e  )\mu\big\}
    \le  e^{- \frac{\e  ^2\mu}{2(1+ \e  /3) } } .
\cr (b) & &{\mathbb P}\big\{S_n \le  (1-\e  )\mu\big\}\le    e^{- \frac{\e  ^2\mu}{2}}.
 \end{eqnarray*}
\end{lemma}
We deduce,   noticing that $ e^{- \frac{\e  ^2\mu}{2}}\le  e^{- \frac{\e  ^2\mu}{2(1+ \e  /3)}}$,
 \begin{equation*}  \rho_n(h)
 \ =\ {\mathbb P}\big\{\sum_{k=1}^n \e_k>(1+ h) \Theta_n\big\}+ {\mathbb P}\big\{\sum_{k=1}^n \e_k<(1- h) \Theta_n\big\}
 \ \le \  2 e^{- \frac{h^2\Theta_n}{2(1+ h/3)}}.
 \end{equation*}
   By assumption  $\frac{  \log \Theta_n }{\Theta_n}\le  {1}/{14} $.   Thus
$h_n:=\sqrt{\frac{7 \log \Theta_n} {2\Theta_n}}\le 1/2$ and so  $\frac{h_n^2\Theta_n}{2(1+ h_n/3)}\ge  {(3/2)\log
\Theta_n}
 $. It follows that
\begin{eqnarray}  \label{hn}\rho_n(h_n)   &\le  &  {2}\, { \Theta_n^{ - 3/2}}. \end{eqnarray}
Let $C_2=   2^{7/2}\max( C_1 , 1)$.
 Further
  \begin{align*}   {C_1
 \over
\sqrt{   (1-h_n)\Theta_n} }& \big(H_n +      \frac{1}{(1-h_n)\Theta_n}
 \big) + \rho_n(h_n) \, \le \,   2^{  1/2}C_1{  H_n
 \over
\sqrt{   \Theta_n} }  +    {2^{  3/2}C_1 +2
 \over
 \Theta_n^{  3/2} }
\cr &  \, \le \,   {  2^{  1/2}\max(C_1, 1)
 \over
\sqrt{   \Theta_n} }\big( H_n  +    {2  +\sqrt 2
 \over
 \Theta_n  } \big)
 \, \le \,    {  C_2
 \over
\sqrt{   \Theta_n} }\big( H_n  +    {1
 \over
 \Theta_n  } \big)   .
\end{align*}
 Therefore
 \begin{eqnarray*}      {\mathbb P}\{S_n =\k\} &\le  &         {  D  ( 1+4h_n )\over
\sqrt{2  \pi  {\rm Var}(S_n)  } }  e^{-\frac{(\kappa-  {\mathbb E\,} S_n)^2}{  2(1  +  h_n){\rm Var}(S_n)   }}     +   {C_2  \over \sqrt{   \Theta_n}
}\big( H_n +    {1 \over  \Theta_n  }   \big)    .
\end{eqnarray*}
 Besides
\begin{align*}
 & {C_1
 \over
\sqrt{   (1-h_n)\Theta_n} } \big(H_n  +   \frac{1}{(1-h_n)\Theta_n}  +  2\rho_n(h_n)
 \big) + \rho_n(h_n)
\cr & \, \le\,      { 2^{1/2}C_1
 \over
\sqrt{   \Theta_n} } \big( H_n   +    {  6
 \over
 \Theta_n  }\big) +  { 2
 \over
 \Theta_n^{3/2}  }
 \cr & \le     \,     { 2^{1/2}\max( C_1 , 1)
 \over
\sqrt{   \Theta_n} } \big( H_n   +    {  6+\sqrt 2
 \over
 \Theta_n  }\big)
 \le     { C_2 \over \sqrt{   \Theta_n} } \big( H_n   +    {  1 \over \Theta_n  }\big)
  .
  \end{align*}
Consequently,
\begin{eqnarray*}   {\mathbb P} \{S_n =\kappa \}   &\ge &      { D(1-2h_n) \over \sqrt{2\pi {\rm Var}(S_n)     }}      {  e^{- \frac{(\kappa-   {\mathbb E\,} S_n)^2}{
2(1-h_n){\rm Var}(S_n)    } }}   -    { C_2   \over \sqrt{   \Theta_n} } \big(H_n +    {1 \over \Theta_n  }\big) .
  \end{eqnarray*}
 \smallskip\par
(1) If  $  \frac{(\kappa-  {\mathbb E\,} S_n)^2}{  2 {\rm
Var}(S_n)  }
 \le \frac{   1+h_n   }{ h_n}$,
 then by using the   inequalities $e^u\le 1+3u$ and $Xe^{-X}\le e^{-1}$ valid for $0\le u\le 1$,   $X\ge 0$, we   get
 \begin{eqnarray*}  e^{- \frac{(\kappa-  {\mathbb E\,} S_n)^2}{   2(1-h_n){\rm Var}(S_n)    } }
&= &    e^{- \frac{(\kappa-  {\mathbb E\,} S_n)^2}{    2 {\rm Var}(S_n)    } } e^{\frac{(\kappa-  {\mathbb E\,} S_n)^2}{  2 {\rm
Var}(S_n)  }\frac{ h_n}{   1+h_n   } }
\cr&\le &    e^{- \frac{(\kappa-  {\mathbb E\,} S_n)^2}{    2 {\rm Var}(S_n)    } }   \Big\{ 1+ 3\frac{(\kappa-  {\mathbb E\,} S_n)^2}{  2 {\rm
Var}(S_n)  }\frac{ h_n}{   1+h_n   }
\Big\}
 \cr&\le &      e^{- \frac{(\kappa-  {\mathbb E\,} S_n)^2}{    2 {\rm Var}(S_n)    } }+ \frac{ 3h_n}{   e(1+h_n  ) }  \le  e^{- \frac{(\kappa-  {\mathbb E\,} S_n)^2}{    2
 {\rm Var}(S_n)    } }+   2h_n .
\end{eqnarray*}
Hence,
\begin{eqnarray*}    {  D  ( 1+4h_n )\over
\sqrt{2  \pi  {\rm Var}(S_n)  } } e^{- \frac{(\kappa-  {\mathbb E\,} S_n)^2}{   2(1-h_n){\rm Var}(S_n)    } }
 &\le &     {  D  ( 1+4h_n )\over
\sqrt{2  \pi  {\rm Var}(S_n)  } }   \big\{  e^{- \frac{(\kappa-  {\mathbb E\,} S_n)^2}{    2 {\rm Var}(S_n)    } }+ 2h_n
\big\}
  \cr&\le &   {  D  e^{- \frac{(\kappa-  {\mathbb E\,} S_n)^2}{    2 {\rm Var}(S_n)    } }   \over
\sqrt{2  \pi  {\rm Var}(S_n)  } }   +{  10h_n D   \over
\sqrt{2  \pi  {\rm Var}(S_n)  } }
.\end{eqnarray*}
Therefore,
 recalling that $ h_n=\sqrt{\frac{7 \log \Theta_n} {2\Theta_n}}$,
\begin{eqnarray*}      {\mathbb P}\{S_n =\k\} -     {  D  e^{- \frac{(\kappa-  {\mathbb E\,} S_n)^2}{    2 {\rm Var}(S_n)    } }   \over
\sqrt{2  \pi  {\rm Var}(S_n)  } }  &\le  & {  10h_n D   \over
\sqrt{2  \pi  {\rm Var}(S_n)  } }   +   C_2\, {     H_n +  \Theta_n^{-1}
\over \sqrt{   \Theta_n} }
 \cr   &\le   & C_2\Big\{ D\big({    {   \log \Theta_n }     \over
 {    {\rm Var}(S_n)   \Theta_n} } \big)^{1/2}  +    {     H_n +  \Theta_n^{-1}
\over \sqrt{   \Theta_n} } \Big\} .
\end{eqnarray*}
since $    5\sqrt{ 7     / \pi  }   \le C_2$.

 \smallskip\par
(2) If  $  \frac{(\kappa-  {\mathbb E\,} S_n)^2}{  2 {\rm
Var}(S_n)  } \le \frac{   1    }{2 h_n}$,
 then as $e^{- u}\ge 1-3u$ if $0\le u\le 1$, we get
\begin{eqnarray*}  e^{- \frac{(\kappa-  {\mathbb E\,} S_n)^2}{   2(1-h_n){\rm
Var}(S_n)    } } &\ge  &    e^{- \frac{(\kappa-  {\mathbb E\,} S_n)^2}{    2 {\rm Var}(S_n)    } } e^{-h_n\frac{(\kappa-  {\mathbb E\,} S_n)^2}{    {\rm
 Var}(S_n)  }  }
 \ge      e^{- \frac{(\kappa-  {\mathbb E\,} S_n)^2}{    2 {\rm Var}(S_n)    } }   \Big\{ 1- 3h_n\frac{(\kappa-  {\mathbb E\,} S_n)^2}{    {\rm
 Var}(S_n)  }
\Big\}
 \cr&\ge &      e^{- \frac{(\kappa-  {\mathbb E\,} S_n)^2}{    2 {\rm Var}(S_n)    } }- \frac{ 3h_n}{   e  }  \ge  e^{- \frac{(\kappa-  {\mathbb E\,} S_n)^2}{    2 {\rm
 Var}(S_n)    } }-   2h_n .
\end{eqnarray*}

Hence,
  \begin{eqnarray*}   { D(1-2h_n) \over \sqrt{2\pi {\rm Var}(S_n)     }}      {  e^{- \frac{(\kappa-   {\mathbb E\,}  S_n)^2}{  2(1-h_n){\rm Var}(S_n)    } }}
 &\ge &    { D(1-2h_n) \over \sqrt{2\pi {\rm Var}(S_n)     }}   \big\{ e^{- \frac{(\kappa-  {\mathbb E\,} S_n)^2}{    2 {\rm Var}(S_n)    } }-   2h_n
\big\}
 \cr &\ge &   { D \over \sqrt{2\pi {\rm Var}(S_n)     }}    e^{- \frac{(\kappa-  {\mathbb E\,} S_n)^2}{    2 {\rm Var}(S_n)    } } -{3h_n D  \over \sqrt{2\pi {\rm Var}(S_n)
 }}
    .\end{eqnarray*}
Consequently,
  \begin{eqnarray*}   {\mathbb P} \{S_n =\kappa \} -{ D \over \sqrt{2\pi {\rm Var}(S_n)     }}    e^{- \frac{(\kappa-  {\mathbb E\,} S_n)^2}{    2 {\rm Var}(S_n)    } }
  &\ge &     -{3h_n D  \over \sqrt{2\pi {\rm Var}(S_n)     }}   -    { C_2   \over \sqrt{
\Theta_n} } \big(H_n +    {1 \over \Theta_n  }\big)
\cr &\ge &-C_2\Big\{ D\big({    {   \log \Theta_n }     \over
 {    {\rm Var}(S_n)   \Theta_n} } \big)^{1/2}  +    {     H_n +  \Theta_n^{-1}
\over \sqrt{   \Theta_n} } \Big\} .
   \end{eqnarray*}
\end{proof}

The following corollary  of Theorem \ref{ger1} implies  when restricted to the i.i.d. case,  Gnedenko's Theorem \ref{gnedenko}. One can also deduce  a strong form of the
local limit theorem as  in Gamkrelidze \cite{G2} and provide an effective bound. See  \cite{GW3}.

\begin{corollary}\label{ger3}
Assume that
\begin{equation}\label{c2a}\lim_{n\to\infty}  \big(\frac{ {\rm Var}(S_n)     }{   \Theta_n  }\big)^{1/2}  \big(    H_n +  \frac{ 1     }{
\Theta_n  }\big)  =0.
\end{equation}
Then
  \begin{eqnarray}\label{r2a} \lim_{n\to\infty} \sup_{ \k\in \mathcal L( v_{
0}n,D ) }\Big| \sqrt{ {\rm Var}(S_n)     }{\mathbb P} \{S_n
=\kappa \} -{ D    e^{- \frac{(\kappa-  {\mathbb E\,} S_n)^2}{    2 {\rm Var}(S_n)    } }  \over \sqrt{2\pi     }} \Big|    & = & 0.
  \end{eqnarray}
\end{corollary}
   Condition \eqref{c2a} is for instance satisfied if
$${\rm (i)}\ \lim_{n\to\infty}   {\rm Var}(S_n)  =\infty, \qq {\rm (ii)}\ \lim_{n\to\infty}       H_n =0, \qq  {\rm (iii)}\
\limsup_{n\to\infty}
\frac{ {\rm Var}(S_n)     }{
\Theta_n  }  <\infty,
    $$
  If $X_i$ are i.i.d., then   $\lim_{n\to \infty}H_n= 0$  as a consequence of    \cite[Th.\,8 p.\,118]{P}.  Condition (ii) is thus satisfied. As conditions (i) and (iii)
  trivially hold, Corollary \ref{ger3} applies and   Gnedenko's theorem follows from   \eqref{r2a}.

\begin{proof}[Proof of Corollary \ref{ger3}]
It follows from Corollary \ref{ger2}  and assumption \eqref{c2a} that \begin{eqnarray*}\label{llta} \lim_{n\to\infty} \sup_{\frac{(\kappa-  {\mathbb E\,} S_n)^2}{    {\rm
Var}(S_n)  } \le  ({\frac { \Theta_n} {14 \log \Theta_n} } )^{1/2}}\Big| \sqrt{ {\rm Var}(S_n)     }{\mathbb P} \{S_n
=\kappa \} -{ D    e^{- \frac{(\kappa-  {\mathbb E\,} S_n)^2}{    2 {\rm Var}(S_n)    } }  \over \sqrt{2\pi     }} \Big|    & = & 0.
   \end{eqnarray*}
Now if $\frac{(\kappa-  {\mathbb E\,} S_n)^2}{    {\rm Var}(S_n)  } >  ({\frac { \Theta_n} {14 \log \Theta_n} } )^{1/2}$, then
$\exp\{- \frac{(\kappa-  {\mathbb E\,} S_n)^2}{    2 {\rm Var}(S_n)    } \}\le \exp\{- \frac{1}{    2   }({\frac { \Theta_n} {14 \log \Theta_n} } )^{1/2}\}$.
By using the first part of Theorem \ref{ger1},  with $h=h_n$   and \eqref{hn},
 \begin{eqnarray*}   & &  \sqrt{   {\rm Var}(S_n)  }  {\mathbb P}\{S_n =\k\} \,\le  \,
      \Big( \frac{1+ h_n  }{ 1-h_n}\Big) \, {  D \over
\sqrt{2  \pi   } } e^{-\frac{(\kappa-  {\mathbb E\,} S_n)^2}{  2(1  +  h_n){\rm Var}(S_n)   } }
 \cr & &\quad   + {C_1   \over
\sqrt{   (1-h_n) } } \, \Big( {      {\rm Var}(S_n)
 \over
 \Theta_n  } \Big)^{1/2}\big(H_n +      \frac{1}{(1-h_n)\Theta_n}
 \big) +   \sqrt{   {\rm Var}(S_n)  }\rho_n(h_n)
 \cr &\le &
      {  3D \over
\sqrt{2  \pi   } } e^{-\frac{1}{ 3  }({\frac { \Theta_n} {14 \log \Theta_n} } )^{1/2} } +  C_1 \sqrt 2\, \Big( {      {\rm Var}(S_n)
 \over
 \Theta_n  } \Big)^{1/2}\big(H_n +      \frac{1}{(1-h_n)\Theta_n}
 \big)
 \cr & &\quad     + 2  \sqrt{   {\rm Var}(S_n)  } \,  { \Theta_n^{ - 3/2}}    .
\end{eqnarray*}
Thus by assumption \eqref{c2a},
\begin{eqnarray*}\label{llta1} \lim_{n\to\infty} \sup_{\frac{(\kappa-  {\mathbb E\,} S_n)^2}{    {\rm Var}(S_n)  } >  ({\frac { \Theta_n} {14 \log \Theta_n} } )^{1/2}}\Big|
\sqrt{ {\rm Var}(S_n)     }{\mathbb P} \{S_n
=\kappa \} -{ D    e^{- {(\kappa-  {\mathbb E\,} S_n)^2}/{    2 {\rm Var}(S_n)    } }  \over \sqrt{2\pi     }} \Big|    & = & 0.
   \end{eqnarray*}
   Consequently,
   \begin{eqnarray*} \lim_{n\to\infty} \sup_{ \k\in \mathcal L( v_{
0}n,D ) }\Big| \sqrt{ {\rm Var}(S_n)     }{\mathbb P} \{S_n
=\kappa \} -{ D    e^{-  {(\kappa-  {\mathbb E\,} S_n)^2}/{    2 {\rm Var}(S_n)    } }  \over \sqrt{2\pi     }} \Big|    & = & 0.
   \end{eqnarray*}
   \end{proof}

\begin{remark} Let $\psi:\R\to \R^+$ be even, convex and such that   $\frac
{\psi(x)}{x^2}$  and $\frac{x^3}{\psi(x)}$  are non-decreasing on $\R^+$. Further assume that
 $ {\mathbb E\,} \psi( X_j )<\infty $.
 Then under the conditions of Corollary \ref{ger2}, we have the following   strengthening
\begin{eqnarray*}   \Big| {\mathbb P} \{S_n =\kappa \} -{ D e^{- \frac{(\kappa-  {\mathbb E\,} S_n)^2}{    2 {\rm Var}(S_n)    } } \over \sqrt{2\pi {\rm Var}(S_n)     }}
\Big|         & \le &    C_3\Big\{D\big({    {   \log \Theta_n }     \over{    {\rm Var}(S_n)   \Theta_n} } \big)^{1/2}  +    {    L_n+  \Theta_n^{-1}\over \sqrt{   \Theta_n}
} \Big\} ,
   \end{eqnarray*}
where $C_3$ an explicit  constant  and $L_n= \big({  \sum_{j=1}^n{\mathbb E\,} \psi (X_j)  }\big)/
\,{   \psi (\sqrt
{ {\rm Var}(S_n )})}$.
\end{remark}

\vskip 3  pt
\begin{problem} A related difficult question can be presented as follows. Let $\{k_j, j\ge 1\}$ be  an increasing sequence of positive integers and $\{p_j, j\ge 1\}$ be a
sequence of reals in $]0,1[$.   Let
$\b_j$ are independent binomial random variables with
\begin{equation}\label{kj.llt}
\begin{cases}
{\mathbb P}\{ \b_j=1\}= p_j,\cr
{\mathbb P}\{ \b_j=0\}= 1- p_j
\end{cases}
\qq (j\ge 1). \end{equation}
 Describe the   LLT for the sequence
$S_n= k_1\b_1+
\ldots + k_n\b_n$, $n\ge 1$.
\end{problem}

{\color{red}
\vskip 15 pt

 \noi{\tf Local limit theorems and the Landau--Kolmogorov inequalities.} In this subsection we will follow \cite{R}. Let $f:{\mathbb Z}\to{\mathbb R}$ and
for $p\in[1,\infty)$ set
\[ \Vert f\Vert_p = \left(\sum_{k\in{\mathbb Z}} |f(k)|^p\right)^{{1}\over{p}}\quad{\rm and }
\quad \Vert f\Vert_{\infty} = \sup_{k\in{\mathbb Z}} |f(k)|.\]
Define also
\[  \Delta^0f(k)=f(k)\qquad{\rm and}\qquad \Delta^{n+1}f(k)=\Delta^nf(k+1) - \Delta^nf(k).\]
Suppose that $X,Y:\Omega\to{\mathbb Z}$ are random variables and put $F(k)=\P(X\leq k),\>G(k)=\P(Y\leq k)$. It turns out that
the following variant of the Landau-Kolmogorov inequality holds
\begin{equation}
\label{elk}
\Vert \Delta F - \Delta G\Vert_{\infty} \leq C\sqrt{\Vert F - G\Vert_{\infty}}\cdot\sqrt{\Vert\Delta^3 F\Vert_1 +\Vert\Delta^3 G\Vert_1,}
\end{equation}
for some universal constant $C$. Analogous inequalities hold for Wasserstein, total variation and local metrics (\cite{R}, Theorem 2.2).

In particular, if we set in (\ref{elk}) $X=\sum_{k=1}^n \varepsilon_i,$ where $\{\varepsilon_i\}$ are Bernoulli i.i.d., with the
probability of success $p$, and
\[\P\{Y=k\}={{1}\over{\sqrt{2\pi}}}\int_{{{k-1/2-\E(X)}\over{\sqrt{{\rm Var}(X)}}}}^{{{k+1/2-\E(X)}\over{\sqrt{{\rm Var}(X)}}}} e^{-{{x^2}\over{2}}}dx, \]
then by the Berry--Esseen theorem $\sqrt{n}\Vert F - G\Vert_{\infty} =O(1)$. Further by Proposition 3.8 in \cite{R}
\[ n\Vert\Delta^3 F\Vert_1 \leq \left( {{2p+1}\over{1-p}} + {{2(1-p)+1}\over{p}}\right)\]
and by a simple calculus $n\Vert\Delta^3 G\Vert_1=O(1)$. Therefore with (\ref{elk}) we obtain
\[ n^{3/4}\sup_{k\in{\mathbb Z}}|\P\{\sum_{i=1}^n\varepsilon_i = k \} - \P\{Y=k\}| = O(1).\]

Thus this approach does not yield the proper rate. However, there are many dependent sequences for which the Berry-Esseen theorem hold
(possibly with worse rate, see e.g. \cite{Su00}). Also we can substitute the normal distribution with the discrete translated Poisson
distribution (see Lemma 4.1 in \cite{R}). Further, the measure of smoothness ($\sqrt{\Vert\Delta^3 F\Vert_1 +\Vert\Delta^3 G\Vert_1}$)
can be calculated for exchangeable or markovian sequences. So this method allows to obtain improved LLTs for e.g. the magnetization
in the Curie--Weiss model (\cite{R}, Theorem 4.5), for the number of isolated  vertices  and triangles in the  Erd\H os--R\'enyi model
(\cite{R}, Theorems 4.8, 4.11) and for embedded sum of independent random variables (for example as in Proposition \ref{lmd})
(\cite{R}, Theorem 4.14).

}

\vskip 15 pt

 \noi{\tf Mod-$ \boldsymbol \phi$ convergence and local limit theorem.}  This approach initiated by Delbaen, Jacod, Kowalski and Nikeghbali \cite{DKN}, \cite{JKN},  \cite{KN}
 is an attempt to refine convergence in law of normalized sequences of random variables
$$ X_n = \frac{Y_n-m_n}{\sqrt{\s_n}}$$
by looking more carefully at the limiting behavior of the characteristic functions $\phi_n(t) = {\mathbb E\,} e^{it Y_n}$ without normalizing. This behavior may bring in some
specific interesting cases, additional informations in presence of a convergence  in law of $X_n$, to the normal law for instance. This is best  illustrated on examples.

 \smallskip\par   Consider first the following example. Let $\varpi_n$ be a random variable counting the number of distinct cycles in a uniformly chosen permutation $\s$ of
 $\{1,\ldots, n\}$. Goncharov proved that $X_n=(\varpi_n-\log n)/\sqrt{ \log n}$ converges in law to $\mathcal N(0,1)$, and Kolchin proved the local limit theorem, see
 \eqref{lltrand.perm}. The characteristic function of $\varpi_n$ is given by
\begin{equation*}
{\mathbb E\,} \big( e^{it\varpi_n}\big) \,=\, \prod_{j=1}^n \big( 1- j^{-1}+ j^{-1}e^{it}\big)
\end{equation*}
Obviously the product diverges as $n\to \infty$, $t\notin 2\pi \Z$. However
\begin{equation*}
{\mathbb E\,} \big( e^{it\varpi_n}\big) \,=\, \prod_{j=1}^n \big( 1+(e^{it}-1)j^{-1}\big) (1+ j^{-1})^{1-e^{it}}\, e^{(e^{it}-1)H_n},
\end{equation*}
where $H_n= 1+1/2 + \ldots+ 1/n$.
One observes that the second term is the characteristic function of a Poisson random variable $P_{H_n}$ of parameter $H_n$. The first term converges since
$$ \prod_{j\ge 1} \Big(1 + \frac{z}{j}\Big) \Big(1 + \frac{1}{j}\Big)^{-z} =\frac{1}{\Gamma (1+z)}, $$
for any $z\in \C$. So that for any $t\in \R$.
\begin{equation*}
{\mathbb E\,} \big( e^{it\varpi_n}\big) \,\sim\, \frac{1}{\Gamma (1+e^{it})}\,e^{it P_{H_n}}, \qq \qq n\to \infty.
\end{equation*}
This suggests a decomposition of type $\varpi = X_n+ Y_n$ where $X_n$ is Poisson and $Y_n$ is independent of $X_n$. However $1/{\Gamma (1+e^{it})}$ is not a characteristic
function.
The authors called this behavior Mod-Poisson convergence with limiting function $1/{\Gamma (1+e^{it})}$.

  \smallskip\par  Consider a second example. Let $\o(k)$ be the prime divisor function (counting the number of primes dividing $k$). Let for each integrer $n$, $N_n$ be
  uniformly distributed in  $\{1, \ldots, n\}$. Erd\H{o}s  and Kac proved that $( \o(N_n)-\log\log n)/\sqrt{\log\log n}$ converges in law to $\mathcal N(0,1)$.
  R\'enyi and Tur\'an proved
\begin{equation*} \frac{1}{n} \sum_{k\le n} e^{i\o(k)} ={\mathbb E\,} \big( e^{itP_{\log\log n}}\big) \Phi(t) \big( 1+ o(1)\big)\qq \quad n\to \infty,
\end{equation*}
where
\begin{equation*}\Phi(t) \,=\,\frac{1}{\Gamma (e^{it})}\,\prod_p \Big(1-\frac{1}{p}\Big)^{e^{it}}\Big(1+ \frac{e^{it}}{p-1}\Big) .\end{equation*}
Moreover  the infinite product is also the limiting function for
$$X_n \buildrel{\mathcal D}\over{=}\sum_{p\le n} B_p$$
where $B_p $ are independent binomial random variables with ${\mathbb P}\{B_p=1\}= 1/p$.

 \smallskip\par  Finally consider a third  example extracted from random matrix theory.  Keating and Snaith proved that if $X_n$ is a random matrix taking values in the
 unitary group $U(n)$, distributed according to the natural Haar measure, and $P_n(T) = \det(1 - TX_n)$, then
\begin{equation}
{\mathbb E\,} e^{it|P_n(1)|}= e^{-(\log n)t^2/2}\, \frac{G(1 + it/2)2 )^2}{G(1+it)} (1+o(1)),
\end{equation}
locally uniformly for $t \in\R$. Here $G(z)$ is the Barnes function, which is
holomorphic of order 2 and such that $G (1) = 1$ and $G(z + 1) = \Gamma(z)G(z)$, see \cite{Ada}.

\bigskip
Introduce a definition. Let $\mu$ be a probability law on $\R^d$, $d\ge 1$, with characteristic function $\phi$. Let $(X_n)_{n\ge 1}$ be $\R^d$-valued random vectors with
characteristic functions $(\phi_n)_{n\ge 1}$.

\begin{definition}\label{mod.phi} We say that there is a  mod-$\phi$ convergence if there exist $A_n\in{\rm GL}_d(\R)$, $n\ge 1$, such that
\begin{eqnarray*}
{\rm (H1)}& & \hbox{$\phi$ is integrable on $\R^d$},
\cr
{\rm (H2)}& & \hbox{Denoting $\Sigma_n=A_n^{-1}$ we have $\Sigma_n\to 0$, and the vectors $Y_n=A_n^{-1}(X_n)$}
\cr & &  \hbox{converge in law to $\m$},
\cr {\rm (H3)}& & \hbox{For all $k\ge 0$, we have}
\cr  & &
\qq \qq \qq \sup_{n\ge 1} \int_{|t|\ge a\atop |\Sigma_n^*t)|\le k} \big|\phi_n\big(\Sigma_n^*t\big) {\rm d} t \quad \to \ 0\quad \hbox{as $a\to \infty$}.
\end{eqnarray*}
\end{definition}

Assumption (H1) implies that $\m$ admits a density, call it $\a$.
We note that $ {\mathbb E\,} e^{it.Y_n}= \phi_n\big(\Sigma_n^*t\big)$.
Assumption (H3)  is for instance satisfied if for each $k$, there exists an integrable function $h_k(t)$ such that
$$ \phi_n\big(\Sigma_n^*t\big)\le h_k(t),$$
for all $n$ and all $t$ such that $|\Sigma_n^*t)|\le k$.

\vskip   3 pt
Delbaen, Kowalski and Nikeghbali proved the following local limit theorem.

\begin{theorem} Assume {\rm mod}-$\phi$ convergence for the sequence $(X_n)_{n\ge 1}$. Then for $f$ continuous and compactly supported we have
$${\mathbb E\,} f(X_n) = \a(0) |\det(A_n)|^{-1} \Big(\int_{\R^d} f(x) {\rm d} x\Big) \big( 1+ o(1) \big), \qq \quad n\to \infty.  $$
\end{theorem}

No rate of convergence is however indicated. In \cite[Th.\,4]{KN}, a sharper version of this result with rate of convergence is given, however under  very restrictive
conditions.
\vskip   3 pt It is worth quoting, in relation with the above, that Hwang introduced in \cite[p.\,451]{Hw3} an analytic assumption  on the probability generating functions of
integer-valued random variables $(X_n)_n$, that is, on the power series
$\sum_{n\ge 1}{\mathbb P}(X_N =n)z^n=\E(z^{X_N})$,
which is very closely related to mod-Poisson convergence. This assumption is used as a basis to deduce results on Poisson approximation of the sequence, with  applications.

\vskip 15 pt
\noi{\tf Local limit theorem under tightness conditions.} Recently  Dolgopyat \cite{Dolgo} obtained a local limit theorem for sums of independent random vectors satisfying
appropriate tightness assumptions. In particular, his result holds in dimension one for independent uniformly bounded summands, improving  Prohorov's Theorem.
Let $\{X_j, j\ge 1\}$ be independent $\R^d$-valued random variables such that
\begin{eqnarray}\label{dolgo1}& & \E X_j =0,
\end{eqnarray}
 \begin{eqnarray}\label{dolgo2} \E |X_j|^3\le m,
\end{eqnarray}
and assume there exists a constant $\e_0>0$ such that for each $s\in \R^d$,
\begin{eqnarray}\label{dolgo3}\qq \qq\qq \E \langle X_j,s\rangle^2\ge \e_0|s|^2.
\end{eqnarray}
Let $S_N= \sum_{j=1}^N X_j$ and let $V_N$ denotes its covariance matrix
\begin{equation}  V_{N, l_1,l_2} =  \sum_{j=1}^N \E \big(X_{j,l_1}X_{j,l_2} \big) .
\end{equation}
Under condition \eqref{dolgo2}, condition \eqref{dolgo3} is equivalent to the existence of $\e_1, \e_2>0$ such that for each proper affine subspace $\Pi \subset \R^d$, we
have
\begin{eqnarray}\label{dolgo4} \P\{ d(X_j,\Pi)\le \e_1\}\le 1-\e_2.\end{eqnarray}

 \vskip 3 pt  If $\mathcal H$ is a proper subgroup of $\R^d$ we call the sequence $\{X_N\}$ arithmetic, otherwise it is called non-arithmetic.

Call a closed subgroup $H$ of $\R^d$     sufficient if there is a deterministic sequence $a_N$  such that $S_N -a_N$    mod ($H$)  converges almost surely.  The minimal
subgroup, denoted by
$\mathcal H$, is defined as the intersection of all sufficient subgroups. {\color{green} The Haar measure on $\mathcal H$ is defined as follows. $\mathcal H$ is isomorphic to
the product of $\Z^{d_1} \times \R^{d-d_1}$. Then $\lambda_\mathcal H$ is the product of the counting measure on the first factor and the Lebesgue measure on the second
factor normalized as follows.}
\vskip 3 pt
 Choose a set $D$ so that each $x\in \R^d$ can be uniquely written as $x=h+\theta$ where $h\in\mathcal H$, $\theta\in D$. $\lambda_\mathcal H$ is normalized so that
\begin{eqnarray}\label{dolgo5}
\int_ {\R^d} g(x)\dd x =\int_\mathcal H\int_D   g(h + \theta)\dd \lambda_\mathcal H(h)\dd \lambda_D(\theta).    \end{eqnarray}
  where $\lambda_D$ is the Lebesgue measure on $D$ normalized to have total volume $1$.

\vskip 3 pt  Given a random variable $Y$, let $C_Y$ be the convolution operator
\begin{equation}
  C_Y (g) (x) =\E g(x+Y).\end{equation}
\begin{proposition}\label{dolgo.prop}{\rm (a)} If $H$ is sufficient then $\mathbb R^d/H$ is compact.
\vskip 2 pt {\rm (b)} The minimal subgroup is sufficient.
\end{proposition}
By the above Proposition,  there exists a bounded sequence $a_N$ such that $S_N -a_N$   mod ($H$)  converges almost surely. Fix such a sequence and denote the limiting random
variable by S.

\vskip 3 pt   We denote by $\mathcal C(\R^d)$ (respectively $\mathcal C_0(\R^d)$) the space of continuous (respectively $r$ times differentiable) functions on $\R^d$. The
subscript $0$ indicates that we consider only
functions of compact support in the corresponding space.

\vskip 3 pt {\color{green}The following result is proved.}
 \begin{theorem}\label{dolgo.LLT}
For each $g\in \mathcal C_0(\R^d)$, for each sequence $z_N=\mathcal O(\sqrt N)$ such that $z_N -a_N \in \mathcal H$ we have
\begin{equation}
\lim_{N\to \infty} \Big[ \frac{\E  g(S_N- z_N)}{u_N(z_N)}\Big] \,=\, \int_\mathcal H \mathcal C_\mathbb S (g) (h) \dd \lambda_\mathcal H (h)
\end{equation}
where $\lambda_\mathcal H$ is the Haar measure on $\mathcal H$ and $u_N(z)$  is the density of the normal random
variable with zero mean and covariance $V_N$.

In particular, in the non-arithmetic case, for each sequence $z_N=\mathcal O(\sqrt N)$ we have,
\begin{equation}
\lim_{N\to \infty} \Big[ \frac{\E  g(S_N- z_N)}{u_N(z_N)}\Big] \,=\, \int_{\R^d} g(x) \dd x.
\end{equation}
\end{theorem}

 Several   examples of applications are provided in \cite{Dolgo}.

[LLT and H\"older-Continuity]$\surd$}  //  {\color{blue}Postnikov$\surd$} // Maller //  {\color{blue}Breuillard [Diophantine measures and LLT]$\surd$}//
{\color{red}Richter$\surd$}// {\color{blue}Rogozin$\surd$}.]]

{\color{red}
\vskip 15 pt
\noi{\tf Local limit theorems and domains of attraction.}
Let $\{X_k\}$ be i.i.d. random sequence such that there exist normalizing sequences $B_n>0,$ $A_n,$ $n\geq 1,$ and for
some  non-degenerate probability measure $\mu$ we have
\[
{{X_1+X_2+\ldots +X_n-A_n}\over{B_n}}\to \mu,\qquad {\rm in}\>\, {\rm distribution}.
\]
Thus $\mu$ is necessarily $\alpha$--stable, $\alpha\in(0,2]$ and we say that ${\mathcal L}(X_1)$ is in the
domain of attraction of $\alpha$--stable law, ${\mathcal L}(X_1)\in{\mathcal DA}(\alpha)$ in short.
The following characterization can be found  on p. 313 in \cite{F}.
\begin{theorem}
\label{dat1}
Suppose ${\mathcal L}(X_1)$ is non-degenerated. Then,\\
(i)$\>$ ${\mathcal L}(X_1)\in {\mathcal DA}(2)$ iff $\E(X_1^2I(|X_1|\leq x))$ is slowly varying;\\
(ii)$\>$ ${\mathcal L}(X_1)\in {\mathcal DA}(\alpha),$ $\alpha\in(0,2)$ iff $x^{\alpha}\P\{|X_1|>x\}$ is slowly varying and
\[ {{\P\{X_1\leq -x \}}\over{\P\{X_1>x\}+ \P\{X_1\leq -x\}}}\to q,\qquad {{\P\{X_1>x\}}\over{\P\{X_|>x\} + \P\{X_1\leq -x\}}}\to p,\qquad {\rm as\quad} x\to\infty\quad
(p+q=1);\]
(iii)$\>$ up to a constant factor normalizing sequence $B_n$ satisfies $n\E(X_1^2I(|X_1|\leq B_n))\sim B_n^2$.
\end{theorem}

The following theorem contains (in full)  both the Gnedenko (\cite{G}) and the Shepp-Stone (\cite{Sh64},\cite{Sto})
local limit theorems (see also \cite{IBLIN} Th. 4.2.1, p.121; \cite{BGT} Th. 8.4.2, p.\,351 and \cite{Na73}, Th. 2).
\begin{theorem}
\label{dat2}
Suppose $\{X_k\}$ is i.i.d.
and $g_{\alpha}$ is a stable density, $\alpha\in(0,2]$.
In order that there exist sequences $\{B_n\},$ $B_n\to\infty,$  $\{A_n\}$ and such that:\\
I)$\>$ for ${\mathcal L}(X_1)$ lattice on $\{v_0+mD\},$ $v_0\in{\mathbb R},$ $m\in{\mathbb Z},$ where $D$
is maximal step we have uniformly in $m\in{\mathbb Z}$
\[ B_n{\P}\{S_n=nv_0+mD\} = dg_{\alpha}\left({{nv_0+mD-A_n}\over{B_n}}\right) + o(1);\]
II)$\>$ or for ${\mathcal L}(X_1)$ non-lattice and any continuous function $G$ with a compact support we have
uniformly in $r\in{\mathbb R}$
\vskip 6pt\noindent
\[ {B_n} { \E}(G(S_n - r)) = g_{\alpha}\left({{r-A_n}\over{B_n}}\right) \int G(y)dy + o(1);\]
it is necessary and sufficient that ${\mathcal L}(X_1)\in {\mathcal DA}(\alpha)$ and the step $D$ is maximal.
\end{theorem}

\vskip 15 pt
\noi{\tf Refinements of the Gnedenko local limit theorem.}
Let $H_k(x)$ denote Hermite--Chebyshev polynomials, i.e.:
\[ H_k(x) = (-1)^k\cdot e^{{x^2}\over{2}}
{{\partial^k}\over{\partial x^k}}\bigg(
 e^{-{x^2}\over{2}}\bigg) =
k!\sum_{i=0}^{\lfloor{{k}\over{2}}\rfloor}
{{(-1)^i x^{k-2i}}\over{i!(k-2i)!2^i}}.\]
Suppose $\{X_k\}$ is a centered i.i.d. random sequence and $\varphi(\theta)=\E(e^{i\theta X_1})$. Define cumulants by
\[ \gamma_k =  {{1}\over{i^k}}
{{\partial^k}\over{\partial\theta^k}}\log{\varphi(\theta)}|_{\theta=0}. \]
We have $\gamma_0=\gamma_1=0,$ $\gamma_2={\rm Var}(X_1)=\sigma^2,$ $\gamma_3=E(X_1^3)=\mu_3$.
Put
\[ {\frak R}_{\nu}(x) =
\sum_{(k_1,k_2,\ldots,k_{\nu})\atop k_i\geq 0, \sum_{i=1}^{\nu} ik_i=\nu}
H_{\nu+2\sum_{i=1}^{\nu} k_i}(x)
\prod_{m=1}^{\nu} {{1}\over{k_m!}}
\Bigl({{\gamma_{m+2}}\over{(m+2)!\sigma^{m+2}}}\Bigl)^{k_m}.
\]
The following result is due to Esseen (see \cite{Es45}, Theorem 5, p.63). The mulitidimensional case can be found in
\cite{Bi70} and \cite{BR86}(Theorem 22.1, p.231).
\begin{theorem}\quad
\label{aet1}
If ${\mathcal L}(X_1)$ is lattice on
$\{v_0+mD\},$ $\E(X_1)=0,$
$\E\vert X_1 \vert^l<\infty,$ $l\geq 3$
then uniformly in $m$
\[
{{\sigma\sqrt{n}}\over{d}}
{\P}\{S_n=nv_0+mD\}
= \phi(y_{nm})\Bigl(1 + \sum_{\nu=3}^l
{{{\frak R}_{\nu-2}(y_{nm})}\over{n^{{\nu-2}\over{2}}}}\Bigl)
+o\left({{1}\over{n^{{l-2}\over{2}}}}\right),
\]
where  $\phi(x)= {{1}\over{\sqrt{2\pi}}}e^{{-x^2}\over{2}}$ and $y_{nm} = {{nv_0+mD}\over{\sigma\sqrt{n}}}$.
In particular for $l=3$ we have
\[
{\P}\{S_n=nv_0+mD\}
= \phi(y_{nm})\Bigl(1  + (y_{nm}^3-3y_{nm}){{\mu_3}\over{6\sigma^3\sqrt{n}}}\Bigl) + o\left({{1}\over{\sqrt{n}}}\right).
\]
\end{theorem}

Now, define
\[ M(s) = {{\partial}\over{\partial s}}\log\rho(s) =
 \sum_{k\geq 2}  {{\gamma_k}\over{(k-1)!}}s^{k-1}
\!\!\!\!\!\!\!, \qquad
\sigma_s^2 = {{\partial}\over{\partial s}} M(s) =
 \sum_{k\geq 2}  {{\gamma_k}\over{(k-2)!}}s^{k-2}
\!\!\!\!\!\!\!, \]
where $\rho(s)=\E(\exp\{s X_1\})<\infty$ for $s\in(s_{-}, s_{+})=:{\mathbb S}$
and $\gamma_k = {{\partial^k}\over{\partial s^k}}\log\rho(s)\vert_{s=0}$.
Consider the equation
\begin{equation}
\label{lde1}
\sigma t = M(s).
\end{equation}
Since $\sigma^2>0$ the unique solution $s$ of (\ref{lde1})
can be expressed for small $t$ in
terms of
$t$ (Cf. \cite{PE95}, p. 180)
\[ s = {{t}\over{\sigma}} - {{\gamma_3}\over{2\sigma^4}}t^2
-\left({{\gamma_4}\over{6\sigma^5}} - {{\gamma_3^2}\over{2\sigma^7}}\right)t^3
-\left({{\gamma_5}\over{24\sigma^6}} - {{5\gamma_3\gamma_4}\over{12\sigma^8}}
+{{5\gamma_3^3}\over{8\sigma^{10}}}\right)t^4
- \cdots .\]
Therefore
\begin{eqnarray*}
\lefteqn {
\ln\rho(s)  - s\sigma t
}\\ &  & =
-{{t^2}\over{2}} + {{\gamma_3}\over{6\sigma^3}}t^3
+ {{\gamma_4\sigma^2 - 3\gamma_3^2}\over{24\sigma^6}}t^4
+ {{\gamma_5\sigma^4 - 10\gamma_3\gamma_4\sigma^2 + 15\gamma_3^3}
\over{120\sigma^9}}t^5 + \cdots
\\&  &
= -{{t^2}\over{2}} + \Lambda(t)t^3
\!\!.
\end{eqnarray*}
Here $\Lambda(t)$ stands for the Cram\'er series (cf. \cite{Cr38}; \cite{PE95}, (5.77), p.180).
The following result is due to Richter (\cite{Ri}, Theorem 3).
\begin{theorem}
\label{ldt2}\quad
Suppose ${\mathcal L}(X_1)$ is lattice on
$\{v_0+mD\},$ $\E(X_1)=0,$
then for $y_{nm}>1,$ $|y_{nm}|=o(\sqrt{n})$
\[
{{\sigma\sqrt{n}}\over{d}}
{\P}\{S_n = nv_0+mD\}
= \phi(y_{nm})
\exp\bigg\{ {{y_{nm}^3}\over{\sqrt{n}}}
\Lambda\bigg({{y_{nm}}\over{\sqrt{n}}}\bigg) \bigg\}\bigg
(1 + O\bigg({{y_{nm}}\over{\sqrt{n}}}\bigg)\bigg).
\]
\end{theorem}
In particular, for $X_k$ with ${\P}\{X_k=-p_2\}=p_1, {\P}\{X_k=p_1\}=p_2,$ $p_1+p_2=1,$
using Stirling's formula  Khintchin in \cite{Kh29} obtained

\begin{eqnarray*}
\lefteqn{\!\!\!\!\!\!\!\!\!\!\!\!
\P\{\sum_{\nu=1}^n X_{\nu}=k-np_2\}={{n}\choose{k}}p_{1}^{n-k}p_{2}^k}\\
& &  = {{1}\over{\sigma\sqrt{2\pi n}}}
e^{-{{(y_{nk})^2}\over{2\sigma^2}}}e^{-{{(y_{nk})^3}\over{\sqrt{n}}}\sum_{\nu=3}^{\infty}
{{p_1^{k-1}-(-p_2)^{k-1}}\over{\nu(\nu-1)\sigma^{k-1}}}({{y_{nk}}\over{\sqrt{n}}})^{\nu-3}}\left(1+O\left({{y_{nk}}\over{\sqrt{n}}}\right)\right),
\end{eqnarray*}
where $\sigma^2=p_1p_2,$ $y_{nk}={{k-np_2}\over{\sigma\sqrt{n}}}$.

For multidimensional case of Theorem \ref{ldt2} see \cite{Ri58} and for dependent case we recall
\cite{CS85} and \cite{SS00}.

Now, assume $\E(X_1)=0$ and consider the equation $M(s)=t$. Let $s=s(t)$ be its unique solution for $t\in (\tau_{-}, \tau_{+})$.
We have the following local limit theorem for large deviations (see \cite{BR60}; \cite{Pe65}, Theorem 6).
\begin{theorem}
\label{ldt3}
Suppose ${\mathcal L}(X_1)$ is lattice on $\{mD\}$, $m\in\Z$ and $\E(X_1)=0$. Then, $\sigma_s^2>0$ for $s\in{\mathbb S}$ and
for $nt$ taking values of the form $mD,$ $m\geq 0$
\[ \sup_{\epsilon\leq t\leq \tau_{+}-\epsilon}
\left\vert \sigma_s\sqrt{2\pi n}{{\P\{S_n=m\})}\over{(\rho(s)e^{-st})^n}}-D\right\vert = O\left({{1}\over{n}}\right),\]
where $\epsilon>0$.
\end{theorem}

\vskip 15 pt
\noi{\tf Local limit theorem for the number of renewals.}
Let $X_k$ be nonnegative i.i.d. random variables and $\E(X_k)=a, {\rm Var}(X_k)=\sigma^2$. Put
$\nu(x)=n$ if $S_n<x$ and $S_{n+1}\geq x,$ $S_0=0$. A. Nagaev in \cite{Na70}
obtained the following local limit theorem for the number of renewals.
\begin{theorem}
\label{rent1}
If $x\to\infty$ then, uniformly in $n$
\[ {\P}\{\nu(x) = n \} ={{1}\over{\sqrt{2\pi x\sigma^2/a^3}}}e ^{-{{(n-x/a)^2}\over{2x\sigma^2/a^3}}} + o\left({{1}\over{\sqrt{x}}}\right).
\]
\end{theorem}
}

{\color{red}
\vskip 15 pt
\noi{\tf Local limit theorems for densities and refinements.}
Let $B_n^{-1}(S_n-A_n)$ has density $p_n(x)$ where $S_n=X_1+\ldots+X_n$ and $X_k$ are i.i.d.
\begin{theorem}[\cite{IBLIN}, Theorem 4.3.1]
\label{dent1}
In order that for some choice of constants $A_n,$ $B_n>0$
\[ \lim_n \sup_x |p_n(x)-g_{\alpha}(x)| =0,  \]
where $g_{\alpha}$ is the density of $\alpha$--stable distribution $(0<\alpha\leq 2),$ it is
necessary and sufficient that ${\mathcal L}(X_1)\in{\mathcal DA}(\alpha)$ and $\sup_x p_{n_0}(x)<\infty$ for some
$n_0\geq 1$.
\end{theorem}
Let $F_n(x)={\P}\{S_n-A_n\leq xB_n\}$ and $F_n = a_n F_n^{\rm ac}+ b_nF_n^s$ where $F_n^{\rm ac}(F_n^s)$ is
absolutely continuous (singular) part.
\begin{theorem}[\cite{IBLIN}, Theorem 4.4.1]
\label{dent2}
In order that for some choice of constants $A_n,$ $B_n>0$
\[ \lim_n \int_{-\infty}^{\infty} |p_n(x)-g_{\alpha}(x)|\,dx =0,  \]
where $g_{\alpha}$ is the density of $\alpha$--stable distribution $(0<\alpha\leq 2),$ it is
necessary and sufficient that ${\mathcal L}(X_1)\in{\mathcal DA}(\alpha)$ and $a_{n_0}>0$ for some
$n_0\geq 1$.
\end{theorem}
For $\alpha=2$  Theorem \ref{dent2} is due Prokhorov (\cite{Pr52}).
For finite dimensional generalization see \cite{MH64} and \cite{Bl02} for infinite dimensional counterexample.
The Orlicz spaces case is treated in \cite{Bo19}. For examples that CLT does not imply LLT see \cite{Sto13} (14.9 and 17.7(i)).

Let $B_n=\sigma\sqrt{n}$ and $A_n=0$ where $\sigma^2={\E}(X_1^2)$ and ${\E}(X_1)=0$. The analog of Theorem \ref{aet1}
for densities is as follows.
\begin{theorem}[\cite{IBLIN}, Theorem 4.5.2]
\label{dent3}
Suppose that ${\E}|X_1|^l<\infty$ and $\sup_x p_{n_0}(x)<\infty$ for some
$n_0\geq 1$. Then,
\[ p_n(x) =
 \phi(x)\Bigl(1 + \sum_{\nu=3}^l
{{{\frak R}_{\nu-2}(x)}\over{n^{{\nu-2}\over{2}}}}\Bigl)
+ o\left({{1}\over{n^{{l-2}\over{2}}}}\right),
\]
where  $\phi(x)= {{1}\over{\sqrt{2\pi}}}e^{{-x^2}\over{2}}$ and ${\mathfrak R}_k$ are as in Theorem \ref{aet1}.
In particular for $l=3$ we have
\[
p_n(x)
= \phi(x)\Bigl(1  + (x^3-3x){{\mu_3}\over{6\sigma^3\sqrt{n}}}\Bigl) + o\left({{1}\over{\sqrt{n}}}\right).
\]
where $\mu_3=\E(X_1^3)$.
\end{theorem}
Let $\E(\exp\{s_0 |X_1|\})<\infty$ for some $s_0>0$
and $\Lambda(t)$ stands for the Cram\'er series (see Theorem \ref{ldt2}).
The following result is due to Richter (\cite{Ri}, Theorem 2).
\begin{theorem}
\label{dent4}\quad
Suppose ${\mathcal L}(X_1)$ has bounded and continuous density. Then, for $x\geq 1,$ $x=o(\sqrt{n})$
\[
p_n(x) = \phi(x)
\exp\bigg\{ {{x^3}\over{\sqrt{n}}}
\Lambda\bigg({{x}\over{\sqrt{n}}}\bigg) \bigg\}\bigg
(1 + O\bigg({{x}\over{\sqrt{n}}}\bigg)\bigg).
\]
\end{theorem}
}

{\color{red}
\vskip 15 pt
\noi{\tf Local limit theorem for residues of linear transforms of sums.}
Let $\{X_k\},$ $k\in\{0\}\cup{\mathbb N},$
be a sequence of non-degenerate random variables taking values on the real line ${\mathbb R}$ and
defined on a probability space $(\Omega, {\mathcal F}, \P)$.
For $\lambda\neq 0$ consider the residues of the linear transform of $S_n=X_1+\ldots+X_n$
\[
\Lambda_n(a)=\left\{ \lambda S_n  +  a\right\}=
(\lambda S_n+a) {\>\rm mod\>} 1= \lambda S_n +a -
\lfloor \lambda S_n +a \rfloor,\quad
\lambda,a\in{\mathbb R},
\]
where $\lfloor x\rfloor = \sup\{k\in{\mathbb Z}\,;\, k \leq x < k + 1 \},$
${\mathbb Z}=\{\ldots,-1,0,1,\ldots\}$.
Observe that $\Lambda_n(a)\in[0,1),$ $\Lambda_n(k)=\Lambda_n(0)$ for
$k\in{\mathbb Z}$ and $\Lambda_n(a)\equiv \{a\}$ if $S_n$
takes only multiples of ${{1}\over{\lambda}}\cdotp$
Let ${\mathcal U}_0$ denote Lebesque measure on interval
$[0,1)\cap{\mathbb R}$ while ${\mathcal U}_{{1}\over{m}}$ is the uniform measure
on $\{0,{{1}\over{m}},\ldots,{{m-1}\over{m}}\},$ $m\in{\mathbb N}$.
Put $\widehat\mu(\theta)=\E[e^{i\theta X_0}]$ and define $\lambda(\widehat\mu)$ via
$\lambda(\widehat\mu)\cdot\inf\{{n\in{\mathbb N}}\,;\,\vert\widehat\mu(2\pi n\lambda)\vert=1\}  = 1,$
where $\lambda(\widehat\mu)=0$ if such $n$ does not exist.
If $\lambda(\widehat\mu)>0$ set
$b=\lambda(\widehat\mu){\rm Arg}(\widehat\mu(2\pi(\lambda(\widehat\mu))^{-1}\lambda))$.
It turns out that the law ${\mathcal L}(\Lambda_n)$ is in the domain of attraction
of ${\mathcal U}_{\lambda(\widehat\mu)}$ under the following condition
\begin{equation}
\label{ecn}
\vert \E[e^{i\theta S_{n}}]\vert \leq (1-c(1-\vert \E[e^{i\theta X_0}]\vert^2))^{\lfloor {{n}\over{q}}\rfloor}
,\qquad\theta\in{\mathbb R},\quad n\in{\mathbb N},
\end{equation}
for some $c\in(0,1]\cap{\mathbb R}$ and an integer $q\geq 2$.
Condition (\ref{ecn}) is satisfied for a functional $f$ defined on digits of continued fractions expansion (see Lemma \ref{incfr})
and is used for local limit theorems for $\psi$--mixing Markov chains (see Lemma 1.5 in \cite{Na61}
and Proposition 7 in \cite{MPP}).
\begin{theorem}[\cite{Sz10a}, Theorem 1]
\label{tmod1}
Suppose that identically distributed, non-degenerate random variables $\{X_k\},$ satisfy
(\ref{ecn}).
If $\lambda(\widehat\mu)=0$ then $\Lambda_n(a)$ converges to  ${\mathcal U}_0$.
If $\lambda(\widehat\mu)>0$ and $S_n$
then $\Lambda_n(-nb\lambda)$ converges to ${\mathcal U}_{\lambda(\widehat\mu)}$.
\end{theorem}
Theorem \ref{tmod1}
includes results in \cite{DW} and \cite{NM66} (for identically distributed random variables) and reveals similar dichotomy
as in the classical local limit theorem (cf. Theorem 10.17 in \cite{3MaSe}).
In view of this it is classified as a local limit theorem (with the normalizing sequence $\equiv 1$).
\begin{corollary}
\label{cmod}
Suppose integer valued, identically distributed, non-degenerate random variables $\{X_k\}$ satisfy (\ref{ecn})
and $I\subset [0, 1)$.
Then either $\P\{\Lambda_n(a) \in I \}\to_n  {\mathcal U}_0(I)$
or $\P\{ \Lambda_n(-nb\lambda\} \in I )\to_n  {\mathcal U}_{{1}\over{m}}(I)$, for some integer $m>1$ and $b\in \R$.
\end{corollary}
It seems that the limiting behaviour of random sums mod 1 was for the first time investigated by L\'evy (\cite{Le39}).
}

{\color{green}\subsection{The ergodic case}}  \label{1.14}
{\color{blue}
Aside the study of intrinsic properties of the various kind of ergodic sums, such as $L^p$-convergence, mixing properties, spectral representation, almost everywhere
convergence, recurrence, variational inequalities,  and spectral regularization inequalities, the central limit theorem and  the local limit theorem are tools of
investigation in this area, so different from Probability theory. Results obtained   in this far-off setting   show
the validity of these limit theorems (for bounded variation functions) for a broad variety of transformations.
This is a very active research's area, covered by an abundant and still   flourishing  literature.

\vskip 2 pt   We essentially consider  local limit theorems for piecewise monotonic transformations, and   the used methods.
Local limit theorems related to sub-shifts of finite type, toral automorphims are not considered, we refer to Roger's Thesis \cite{Roge} for instance, for a good introduction
and results.  We do not  consider either related results obtained in other important specific areas, such as Sina\"i's billard. We refer for an introduction to  Bunimovich,
Sinai and Chernov \cite{BCS}, for instance.
  \vskip 3 pt
Concerning the (numerous) results  for interval expanding maps, we refer, among many other papers,   to Rousseau-Egele \cite{Rou}, Broise \cite{Bro} nice and informative
synthesis' article  completing Rousseau-Egele's,   Calderoni, Campanino  and Capocaccia \cite{CCC},  Gou\"ezel \cite{Gou} for extensions to non uniformly expanding maps. A
typical historical example is the transformation $(2x)\, {\rm mod}(1)$, for which Kac \cite{K} established in 1940 the validity of the central limit theorem for a class of
Lipschitz functions.

 \vskip 3 pt\noi{\tf Local limit theorem for interval expanding maps.} Let $I=[0,1]$ endowed with the normalized Lebesgue measure $m$. Here we discuss an important family of
 expanding   transformations for which the central and local limit theorems hold for   large classes of functions.
\vskip 2 pt  Lasota and Yorke have \cite{LY} shown that if $T\colon I\to I$ is an expanding  transformation, which is piecewise of class $C^2$, then there exists a measure
$\m$ on $I$, which is $T$-invariant, absolutely continuous with respect to $m$ and whose density has bounded variation.  Put  for  $f\in L^0(m)$ and $n\ge 1$,
$S_nf=\sum_{k=0}^{n-1} f\circ T^k$. The existence of the central limit theorem (for functions with bounded variation) was proved by Wong \cite{Wo}, provided that the limit
$\s^2  =\lim_{n\to \infty}\frac{{\rm Var}(S_n)}{n}$ exists and is positive.

\vskip 2 pt The approach used for proving the central limit theorem
 and
 the local limit theorem (non lattice and lattice case) is spectral, and based on a method of spectral decomposition of the operator, already used for   Markov chains by
 Doeblin \cite{Do} and Fortet \cite{For}, whom in particular studied the continued fraction transformation, and the transformation $(2x)\, {\rm mod}(1)$. This transformation
 was also studied by Kac \cite{K}, by means of  another approach based on Probability theory and Fourier analysis.

\vskip 3 pt In  a nicely written paper \cite{Rou}, Rousseau-Egele investigated    these questions more precisely (central limit theorem with remainder term,  and local limit
theorem) for  expanding  piecewise monotonic transformations.  Let $T\colon I\to I$ be an application enjoying the following property. There exists   a countable subdivision
{\color{green}$(I_j)$} of $I$ where $I_j=(a_{j-1},a_j)$ is an open interval  such that:
\vskip 2 pt {\rm (1)} (local inversion  condition) The restriction of $T$ to $I_j$ is strictly monotone and can be extended into an  application of  class $C^2$ on $\bar
I_j$.
\vskip 2 pt {\rm (2)} $T(I_j)$ is composed of finitely many distinct intervals.
\vskip 2 pt {\rm (3)}  (dilatation condition)  There exists an integer $n$ such that $\g= \inf_{x\in I} \big|(T^n)'(x) \big| >1$.
\vskip 3 pt
The Perron-Frobenius operator associated to $T$ is the operator defined by
\begin{equation}\label{Per.Fro}
\int_I \phi f. g \,\dd m\,=\,\int_I   f. g\circ T \,\dd m \qq \qq (f\in L^1(m), g\in L^\infty (m))
\end{equation}
This is a positive contraction of $L^1(m)$ and one has $\phi f=f$ if and only if the measure $\m=fm$ is $T$-invariant. Condition (1) allows one to write this operator under
the explicit form
\begin{equation}\label{Per.Fro.1}
 \phi f (x) = \sum_j f(\s_j x)\p_j(x)\chi_j(x),
\end{equation}
where $\s_j$ is the inverse of $T$ over $J_j=T(\bar I_j)$, $\p_j(x)=|\s'_j(x)|$ and $\chi_j$ is the indicator function of $J_j$.

\vskip 2 pt    The variation of  a function $f\colon I\to \C$ is defined as follows:
\begin{equation}\label{bdd.var}
 v(f):= \sup_{0=t_0<\dots<t_n=1 \atop   n\in \N}\sum_{i=1}^{n}
  |f(t_{i})-f(t_{i-1})| .
 \end{equation}
(By Jordan's theorem,   $f$ is of bounded variation, $v(f)<\infty$, if and only if  it can be written as a difference of two nondecreasing functions.) Definition
\eqref{bdd.var}  can be  extended to $f\in L^1(m)$  by taking $v(f) $ as being the infimum of $v(g)$, for $g$ in the class modulo $m$ of $h$. Let   $\mathcal V$ denote  the
sub-space of $L^1(m)$ of functions $f$ such that  $v(f)<\infty $. Then $\mathcal V$ equipped with the norm $\|f\|_{\mathcal V}=  v(f) + \|f\|_{1}$ is a Banach space, which is
dense in $L^1(m)$.

  \begin{theorem}[\cite{Rou},\,Th.\,5]\label{rouss.egel} Let $f$ be of the form
$$ f(x) \,=\, \eta + k(x)$$
where $k(x)$ is integer-valued, with  non-integer  integral, $\eta$ being some real. Assume that conditions (1), (2), (3) are fulfilled. Further assume that $T$ is weakly
mixing, and that for $f\in \mathcal V$, the functional equation $f =k+ \p\circ T- \p$ has no solution  in $\p\in \mathcal V$, $k\in \R$. Then
\begin{align*}\lim_{n\to \infty} \Big| \s\sqrt n\m\big\{ z+S_nf -n\m(f)\in \D\big\}- \frac{1}{\sqrt{2\pi}}e^{-\frac{z^2}{2\s^2 n}}\nu(\D -z -n (\eta -\m(f))\Big| \,=\, 0,
\end{align*}
uniformaly for every real $z$ and any finite interval $\D$,  where $\nu $ is the counting measure on $\Z$. \end{theorem}
\vskip 2 pt    The whole strategy used in the proof consists with studying the spectrum of $\phi$, acting as an operator on     $\mathcal V$, in view of applying the
Ionescu-Tulcea and Marinescu Theorem to $\phi$ and the pair $(L^1(m),\mathcal V)$, in order to get  a decomposition of the operator $\phi ^n$ of the type $\phi
^n=\sum_{i=1}^p \l_i^n\phi_i+\psi^n$, $n\ge 1$.  This decomposition next allows one to estimate the characteristic function of $S_n(f)$, $f\in \mathcal V$.
\vskip 2 pt Theorem \ref{rouss.egel} applies to the following examples (conditions (1), (2), (3) are fulfilled, moreover these transformations are strongly mixing):
\begin{example}[$\b$-transformations]
   $ Tx =\{\b x\}$, $\b>1$.
The associated Perron-Frobenius operator is defined by
$$\phi f(x)= \frac{1}{\b} \sum_{j=0}^{[\b]-1} f\Big( \frac{x+j}{\b}\Big)+\frac{1}{\b}  f\Big( \frac{x+[\b]}{\b}\Big)\,\chi_{[0, [\b]]}(x).$$
  \end{example}

   \begin{example}[Continued fraction transformation]
   $ Tx =\big\{\frac{1}{x}\big\}$, $T(0)= 0$.
The associated Perron-Frobenius operator is
$$\phi f(x)=\sum_{j=1}^\infty f\Big( \frac{1}{j+x}\Big)
  \Big(\frac{1}{j+x}\Big)^2.$$ \end{example}
 \vskip 5 pt  This important example is more developed below, where further  another (non-spectral) approach is introduced, and a local limit theorem   stated with hints of
 proof.}

\vskip 5 pt\noi{\color{red}{{\tf The continued fractions case.}
Any irrational number $x\in (0,1]$ can be uniquely expressed as a
simple non-terminating continued fraction
$x=[0\,;\,a_1(x),a_2(x),a_3(x)\ldots],$ where
\[ [0\,;\, x_1,\ldots,x_n] =
{{1}\over{x_1 + {{1}\over{x_2 + {{1}\over{x_3+ {{ }\atop{\ddots+
{{1}\over{x_n}}}}}}}}}}, \quad n\in{\mathbb N}=\{1,2,\ldots\}
\]
and $a_n(x)\in{\mathbb N}$ (\cite{Kh64}, Theorem 14).
It follows from Euclid's algorithm
that the continued fraction transformation $\,T\,$ defined by
\[ T(x)=
{{1}\over{x}} - \left\lfloor{{1}\over{x}}\right\rfloor
= {{1}\over{x}} - \max\{n\in{\mathbb N}\,;\, nx \leq 1 \}
,\quad x\in(0,1],
\]
generates natural numbers $a_n(x),$ called partial quotients, or digits via
\[
a_1(x)= \left\lfloor{{1}\over{x}}\right\rfloor,\quad a_{n+1}(x)= a_1(T^n(x))\]
(cf. \cite{IK02}, p.14).
Let ${\rm P}$ denote Gauss' measure, i.e.
\[ {\rm P}(A)={{1}\over{\ln{2}}}\int_A {{\lambda(dx)}\over{1+x}},\]
where $A\subseteq(0,1)$ and $\lambda$ is Lebesgue: measurable set and measure, respectively.
It is well-known that ${\rm P}$ is an invariant measure for $T$ (cf. \cite{IK02}, Theorem 1.2.1).
Besides, if ${{\mathcal F}_{\!\!{k}}^{\,m}}$ denote the $\sigma$--algebra generated
by $a_k,a_{k+1},\ldots,a_{m}$ then
by Corollary 1.3.15 in \cite{IK02}
\begin{equation}
\label{psi}
\psi_n =
\sup_{k\in{\mathbb N}}
\sup\{
\vert {{{\rm P}(A\cap B)}\over {{\rm P}(A)\cdot {\rm P}(B)}}-1\vert;
\> {\rm P}(A)\cdot {\rm P}(B)>0,
\> A\in{{\mathcal F}_{\!\!{1}}^{\,k}}
,B\in {\mathcal F}_{n+k}^{\infty}\} \leq C\varrho^n,
\end{equation}
where $C\leq {{2\ln{2}}\over{7-4\sqrt{2}}},$
$\varrho = {{7-4\sqrt{2}}\over{2}}$
and $\psi_1\leq 2\ln{2}-1 < 0.39$
(i.e. $\{a_n\}$ is an exponentially fast $\psi$--mixing sequence (see \cite{Br07})).

In order to state normal local limit theorem
fix a Borel function $f$ taking values on the real line ${\mathbb R}$ and
such that $\E[f(a_1)]=0$.
Define also the sequence $b_n$ as follows:
if $\>\E[f^2(a_1)]<\infty$ then
$b_n = \sigma\sqrt{n}$ where
$\sigma^2 = \E[f^2(a_1)]
+ 2\sum_{k=2}^{\infty} \E[f(a_1)f(a_k)]$
(see \cite{IBLIN}, Theorem 18.5.2)
and if $\>\E[f^2(a_1)]=\infty$ then
$b_n = \sup\{x>0\, ;\, x^{-2}\E[f^2(a_1)I_{[\vert f(a_1)\vert\leq x]}]
\geq {{1}\over{n}} \},$
{\color{green} meaning {\color{orange}  that $b_n^2\sim n\E[f^2(a_1)I_{[\vert f(a_1)\vert\leq b_n]}]$}}.
\begin{theorem}
\label{lltcfr}
Assume $\,\E[f^2(a_1)I_{[\vert f(a_1)\vert\leq x]}]$ is a slowly varying function
and $\E[f(a_1)]=0$.

If
$f(a_1)$ is lattice distributed on $\{v_0+mD: m\in{\mathbb Z}\},$ ${\mathbb Z}=\{\ldots,-1,0,1,\ldots\},$
for some $v_0\in{\mathbb R}$ and $D>0$ the maximal span
then
\begin{equation}
\label{eg1}
\sup_{m\in{\mathbb Z}}\vert\sqrt{2\pi} b_n {\rm P}\{\sum_{\nu=1}^n f(a_{\nu})=nw+kd\}
- D\exp\{-{{(nv_0+mD)^2}\over{2b_n^2}}\} \vert = o(1);
\end{equation}
if $f(a_1)$ is non-lattice distributed then
\begin{equation}
\label{es1}
\sup_{r\in{\mathbb R}}\vert \sqrt{2\pi}b_n
\E[G(\sum_{\nu=1}^n f(a_{\nu})  - r)] -
\exp\{-{{r^2}\over{2b_n^2}}\}
\int G(y)dy \vert
= o(1),
\end{equation}
for any continuous function $G$ with a compact support.

Conversely, assume that (\ref{es1}) or (\ref{eg1})
holds for some sequence $b_n,$
then $\E[f^2(a_1)I_{[f(a_1)\leq x]}]$ is a slowly varying function
in the sense of Karamata,  the span $D$ is maximal and $\E[f(a_1)]=0$.
\end{theorem}

It should be noted that, using spectral perturbation method,
the normal LLT in more general case of Gibbs--Markov maps but under more restrictive assumption
on regularity of the probability tail of marginal distribution was obtained in \cite{AD01}.
For stable, non-normal LLT for continued fraction see e.g. \cite{Le52}, \cite{He87}, \cite{He00}
or \cite{AD01s}.

The  proof of  Theorem \ref{lltcfr} (see \cite{Sz10} for details) does not make  use of the spectral perturbation method (cf. \cite{Na61}, \cite{Sev}, \cite{HH01}).
The general strategy (see \cite{GK}, \S49)
is to use Fourier transform inversion formula
over distributional limit theorem (Theorem 3 in \cite{Sz09}). This is possible by the following lemma inspired by
Lemma 1.5 in \cite{Na61}.
\begin{lemma}
\label{l2cfr}
Let $f$ be a Borel function. Then
\begin{equation}
\label{incfr}
\vert\E_{\lambda}[e^{{\rm i}\theta S_n}]\vert
\leq
(1-{{1}\over{324}}
(1-\vert \E_{\lambda}[e^{{\rm i}\theta f(a_1)}]\vert^2))^{{n-1}\over{2}},\qquad \theta\in{\mathbb R},\>
n\in{\mathbb N},
\end{equation}
where $\E_\lambda$ is the expectation with respect to the Lebesgue measure.
\end{lemma}
Moreover, this method is applicable for sums of additive functionals of non-stationary (non-homogeneous), non-lattice Markov chains
as one can see from the proof of Theorem 1 in \cite{MPP}). We apply this result for the sequence $\{f_k(a_k)\}$ where $f_k$
are Borel function. Define incomplete quotients via $s_n=(s_{n-1}+a_n)^{-1}$ and $s_0=0$. By Corollary 1.2.8 in \cite{IK02}
(see also p.36) $s_n$ is a Markov chain. Moreover, by the proof of Lemma 1 on p.432 in \cite{Sz10} $\{s_n\}$ satisfies
condition (2) in \cite{MPP} with $a=3^{-1}$ and $b=2$.
Now, observe that $f_k(a_k)=f_k(\lfloor 1/s_k\rfloor)$. Define $\sigma_n^2={\mathbb E}(\sum_{k=1}^n f_k(a_k))^2,$
$\tau_n^2=\sum_{k=1}^n {\mathbb E}(f_k^2(a_k)),$ $\phi_k(t)={\mathbb E}(\exp\{itf_k(a_k)\})$ and set $\gamma=b^{-1}a^4$.
By Theorem 1 in \cite{MPP} we get the following.
\begin{corollary}
\label{lltcfrn}
Suppose that ${\mathbb E}(f_k(a_k))=0,$ $k\geq 1$ and
\begin{equation}
\tau_n^{-2}\sum_{k=1}^n {\mathbb E}(f_k(a_k)I(|f_k(a_k)|\geq \epsilon\tau_n))\to 0 \qquad\epsilon>0.
\end{equation}
Assume also that there is $\delta>0$ and $n_0\in{\mathbb N}$ and some real function $g$ such that $\exp\{-g\}$ is integrable on $\mathbb R$
such that for $1\leq |u|\leq \delta\tau_n$ and $n>n_0$ we have
\begin{equation}
{{\gamma}\over{8}}\sum_{k=1}^n (1-|\phi_k(u/\tau_n)|^2) > g(u).
\end{equation}
Finally, assume that for $u\neq 0$ there is $c(u)$, an open interval $O_u$ containing $u$ and an $n_0=n_0(u)$ such that for all
$t\in O_u$, and $n>n_0$ we have
\begin{equation}
{{\gamma}\over{8\ln{\tau_n}}}\sum_{k=1}^n (1-|\phi_k(t)|^2)\geq c(u)>1.
\end{equation}
Then,
\begin{equation}
\label{es1n}
\sup_{r\in{\mathbb R}}\vert \sqrt{2\pi}\sigma_n
{\mathbb E}[G(\sum_{\nu=1}^n f_{\nu}(a_{\nu})  - r)] -
\exp\{-{{r^2}\over{2\sigma_n^2}}\}
\int G(y)dy \vert
= o(1),
\end{equation}
for any continuous function $G$ with a compact support.
\end{corollary}
}}

{\color{blue} \noi {\tf Local limit theorems for Riesz-Raikov sums.}
These sums are of the type
$$ \sum_{0\le n <N} f(\theta^nt),$$
where $\theta >1$ and $f$ belongs to a class of relatively smooth functions. Suppose that  $f$ is an H\"older continuous function and that
$\theta$ is a Pisot-Vijayaragavan number.
In Petit \cite{Pe}, it is proved that
$$ \lim_{N\to \infty} N^{-1} \int_0^1 \big| S_nf(t)\big|^2 \dd t = \s^2<\infty.$$
If $\s^2>0$, then the central limit theorem holds,
$$ \frac{1}{\s \sqrt N}\, S_N f \ \buildrel{\mathcal D}\over {\Rightarrow} \  \mathcal N (0, 1), \qq N\to \infty.$$
Further a speed of convergence is obtained and a local limit theorem established.

}

{\color{green}\subsection{Miscellaneous results}}  \label{1.15}
{\color{red}
\vskip 5 pt (1)
  {\it Linear random fields.} Suppose $X_j=\sum_{i\in\Z^d}a_i\epsilon_{j-i},$ where $\epsilon_i$ are i.i.d with $\E(\epsilon_i)=0$ and
$\E\epsilon_i^2<\infty$. Assume also that $\sum_{i\in\Z^d}a_i^2<\infty.$ Let $\Gamma_n^d$ be a sequence of finite subset of
$\Z^d$ such that ${\rm card}(\Gamma_n^d)\to_n\infty$.
Let $1\leq k \leq J_n$ be integers and $n(k)=(n_1(k),\ldots,n_d(k))\in\Z^d,$ $m(k)=(m_1(k),\ldots,m_d(k))\in\Z^d,$
where $n_i(k)\leq m_i(k),$ $1\leq i\leq d$. Any set of the form $\Gamma_n^d(k)=\prod_{i=1}^d [n_i(k), m_i(k)]\cap\Z^d$ is called a discrete
rectangle.
Set $S_n=\sum_{j\in\Gamma_n^d} X_j$ and $B_n={\rm Var}(S_n)$.
\begin{theorem}[\cite{FPS20}, Theorem 2.2]
\label{tfps}
Suppose that $B_n\to\infty$. If a) $\sum_{i\in\Z^a}|a_i|<\infty$ assume $\epsilon_0$ is non-lattice; if b) $\sum_{i\in\Z^a}|a_i|=\infty$
assume {\color{green}(Cram\'er's condition \eqref{cramer.cond})} $\limsup_{|t|\to\infty}|\E(\exp\{it\epsilon_0\})|<1$ and that the sets $\Gamma_n^d$ are constructed as a
pairwise disjoint union of
$J_n$ discrete rectangles such that
\[\lim_{n\to\infty}J_n^{2/d}\log{(B_n)}\inf_{i\in\Z^d}\left|\sum_{j\in\Gamma_n^d} a_{j-i}\right|^{-2/d}=0.\]
Then, for all complex-valued  continuous functions $h$ such that $|h|\in L^1(\R^d)$
\[\lim_{n\to\infty} \sup_{u\in\R^d}\left| \sqrt{2\pi}B_n\E(h(S_n-u)) - \exp\{-u/2B_n^2\}\int h(x)\lambda(dx)\right|=0,\]
where $\lambda$ is the Lebesgue measure.
\end{theorem}
\vskip 5 pt (2)
  {\it Operator form.}  Let $\{ \xi_k\}_{k\in {\Bbb Z}_{+}}$ be a homogeneous Markov chain with
the transition probability $P(x,A),\ x\in\R,$  $A$ - Borel set, with $\mu$ - the initial distribution
and $\pi$ - the stationary distribution.
Let $f$ be a Borel function $X_k=f(\xi_k)$
and $S_n=\sum_{k=1}^n X_k,\ n=1,2,\dots\quad .$
Let us define the following {\em uniform recurrence} condition $({\Psi})$,
(cf. \cite{INN85}):
there exist $a>0$ and $b<\infty$ such that for every Borel
set $A,$ $\mu(A)>0$,
\begin{equation}
\label{urc}
a\mu(A)\leq P(x,A) \leq
b\mu(A),\qquad \mu-{\rm a.e.}
\end{equation}
Under this condition $\sigma^2=\lim n^{-1}{\rm Var}(S_n)$ exists in $(0,\infty)$.
Let ${\mathfrak n}$ denote the density function
of the standard normal law. We have the following local limit theorems in operator form.
\begin{theorem}[\cite{Sz08b}, Theorem 4]
\label{stl}
  Assume (\ref{urc}). Suppose $X_0$ is lattice on
$\{v_0+mD\}$
for some Borel function $f$ such that $f\in L^2(\mu),$ $E_{\pi}[f]=0$
and $0\leq h\in L^{\infty}(\mu)$
is not vanishing on $\{v_0+mD\}.$ Then
\[ \sup_{km\in{\mathbb Z}}\,{\rm ess}\sup_x\left\vert
{{\sigma\sqrt{n}}\over{d}}
\E[I_{[S_n=nv_0+mD]}h(\xi_n)\,\vert\, \xi_0 = x]
- \phi({{nv_0+mD}\over{\sigma\sqrt{n}}})
\E_{\pi}[h] \right\vert = o(1).\]
\end{theorem}
\begin{theorem}[\cite{Sz08b}, Theorem 3]
\label{st}
  Assume (\ref{urc}). Let $X_0$ be non-lattice
for some $f\in L^2(\mu)$, $E_{\pi}[f]=0$. Then for any
continuous function $G$ with a compact support and
$0\leq h\in L^{\infty}(\mu)$
\[\sup_{r\in{\mathbb R}}\,{\rm ess}\sup_x\vert
\sigma\sqrt{n}\E[G(S_n - r)h(\xi_n)\,\vert\, \xi_0=x]
- \phi({{r}\over{\sigma\sqrt{n}}})\int G(y)dy\,
\E_{\pi}[h]\vert = o(1).\]
\end{theorem}

\vskip 5 pt (3)
  {\it Nonconventional sums.} Suppose $\{\xi_k\}$ is a homogeneous Markov chain with transitions $P(x,A)$ and stationary distribution $\pi$. Consider non-conventional sums
\[ S_n = \sum_{k=1}^n F(\xi_k,\xi_{2k}\ldots, \xi_{lk}),\]
where $l$ is fixed integer and $F$ is a Borel function on $\R^d$.
Assume $\E(F(\xi_k,\xi_{2k}\ldots, \xi_{lk}))=0$ and  $\E(F^2(\xi_k,\xi_{2k}\ldots, \xi_{lk}))<\infty$.
For any real numbers $x_1,\ldots, x_{l-1}$ set
\[ A_{x_1,\ldots,x_{l-1}} = \{h\geq 0\,|\,F(x_1,\ldots, x_{l-1},x)\in \{kh\, ;\, k\in\Z\}\},\]
$\pi$ almost surely,
and
\[ B_{x_1,\ldots,x_{l-1}} = \{h\geq 0\,|\,F(x_1,\ldots, x_{l-1},x)-F(x_1,\ldots, x_{l-1},y)\in \{kh\, ;\, k\in\Z\}\},\]
$\pi^2=\pi\times\pi$ almost surely.
If $B_{a_1,\ldots,x_{l-1}}\neq\emptyset$ define
\[ h(x_1,\ldots,x_{l-1}) = \sup\{h\,|\, h\in  B_{x_1,\ldots,x_{l-1}}\}. \]
We call the case lattice one if there exists $h>0$ such that
\[ h(x_1,\ldots,x_{l-1}) = h \in  A_{x_1,\ldots,x_{l-1}} \]
$\pi^{l-1}=\times_{i=1}^{l-1}\pi$ almost surely.
If
\[ \pi^{l-1}\{(x_1,\ldots,x_{l-1})\,|\,  B_{x_1,\ldots,x_{l-1}} = \emptyset \}>0\]
then we call the case non-lattice.
\begin{theorem}[\cite{HK16}, Theorem 2.6]
\label{tnc}
 Assume that there exists integer $m$ such that (\ref{urc}) holds for $m$-step transition probability $P^m(x,A)$.
 In  the lattice case
\[ \sup_{k\in{\mathbb Z}}\left\vert
{{\sigma\sqrt{n}}\over{h}} \P\{S_n=k\}
- \phi({{kh}\over{\sigma\sqrt{n}}}) \right\vert = o(1).\]
In the non-lattice case for
any continuous function $G$ with a compact support
\[\sup_{r\in{\mathbb R}}\vert
\sigma\sqrt{n} \E (G(S_n - r))
- \phi({{r}\over{\sigma\sqrt{n}}})\int G(y)dy\,\vert = o(1).\]
\end{theorem}
For extensions to dynamical systems see Chapter 2 in \cite{HK18}.
\vskip 5 pt (4)
  {\it A generalization of Richter's theorem.}
Let $\{T_n\}_{n\geq 1}$ be a sequence of random variables taking only integral values with maximal span 1.
Assume that the moment generating function
$M_n(s)=\sum_{k\in{\Z}}P(T_n=k)e^{ks}$ satisfy
$M_n(s)= \exp\{\phi_nu(s)+v(s)\}(1+O(\kappa_n^{-1}))$ as $n\to\infty,$
uniformly for $|s|\leq \rho,$ $s\in\C, \rho>0,$ and that the following hold:
\vskip 0pt\noindent
1. $\lim_n \phi_n=\infty;$
\vskip 0pt\noindent
2. $u(s)$ and  $v(s)$ do not depend on $n$, are analytic for $|s|\leq \rho$ and $u''(0)\neq 0$;
\vskip 0pt\noindent
3. $\lim_n \kappa_n=\infty;$
\vskip 0pt\noindent
4. there exist constants $0<\epsilon\leq \rho$ and $c=c(\epsilon,\rho),$ where $\epsilon>0$ may be taken
arbitrarily small but fixed, such that
\[ \left|{{M_n(r+it)}\over{M_n(r)}}\right| = O(e^{-c\phi_n}), \]
uniformly for $-\rho\leq r\leq \rho$ and $\epsilon\leq |t|\leq \pi$, as $n\to\infty$.

 Set
$u_m=u^{(m)}(0),$ $v_m=v^{(m)}(0),$ $\mu_n=u_1\phi_n,$  $\sigma_n^2=u_2\phi_n$. We have the following
generalizations of Theorem \ref{ldt2} (see also \cite{CS85}, \cite{Kolc} and \cite{Jo14}).
\begin{theorem}[\cite{Hw2}, Theorem 1]
\label{thw2}
If $k=\mu_n+x\sigma_n$, $x=o(\sqrt{\phi_n})$, then
\begin{equation}
\label{hw2e1}
\sigma_n\P\{T_n= k \} ={\phi(x)}\exp\{{{x^3}\over{\sqrt{\phi_n}}}\Lambda({{x}\over{\sqrt{\phi_n}}})\}\left(1
+ \sum_{1\leq m\leq\nu}{{\Pi_m(x)}\over{\sigma_n^m}}
+ O\left( {{|x|^{\nu+1}+1}\over{(\sqrt{\phi_n})^{\nu+1}}} + {{|x|+1}\over{\kappa_n\sqrt{\phi_n}}}\right) \right),
\end{equation}
where $\nu$ is non-negative integer (depending upon the error term $\kappa_n^{-1}$) and $\Pi_m(x)$ are polynomials of degree
$m$ such that $\Pi_{2i}$ has only even powers of $x$ and $\Pi_{2i-1}$ has only odd powers of $x$, $i=1,2,\ldots$.
\end{theorem}
For example,
\[ \Pi_0(x)= 1,\qquad \Pi_1(x) = \left(v_1 - {{u_3}\over{2u_2}} \right)x, \]
\[ \Pi_2(x) = \left({{v_2}\over{2}} + {{v_1^2}\over{2}} - {{v_1u_3}\over{u_2}} + {{5u_3^2}\over{8u_2^2}}
- {{u_4}\over{4u_2}}\right)x^2 + {{1}\over{8}} - {{v_2}\over{2}} - {{v_1^2}\over{2}} + {{v_1u_3}\over{2u_4}}
- {{5u_3^2}\over{24u_4^2}}. \]
}

{\color{blue}
\subsection{Poisson approximation--Another coupling method.}   \label{appendix-2}
Here we   prospect the role of the Poisson distribution in the study of the local limit theorem and briefly describe a few of the most remarkable results.
Under certain simple conditions, binomial distributions can be approximated   by Poisson distributions. This is for instance expressed by Le Cam's inequality below. This
yields if necessary the importance of this law.
The Poisson distribution, however, not always attracted attention. Citing Katti and Rao's review of Haight's 1967 book, one find \lq\lq {\it If someone had asked us last year
to write a book on the Poisson distribution, our immediate reaction would have been, \lq\lq What for?\rq\rq. We are accustomed to the Poisson starting from the first
introductory course in statistics and prone to fed that this is just one of  the simple distributions to be used to illustrate  new theories that one might develop. It was
thrilling to discover how wrong such a judgment can be.}\rq\rq \,Things  much changed since  and  the study of Poisson distribution,  Poisson--Dirichlet distribution,
scale-invariant Poisson processes \cite {ABT}, the  important  approximation properties to sums of independent random variables, together with  the study of logarithmic
combinatorial structures  \cite{ABT1} form   a whole active research's domain with many concrete applications.  
\vskip 4 pt
\noi {\it Poisson approximation to   binomial distribution.}
  Let $X_{n1}, \ldots X_{nn}$ be $n\ge 1$ independent random variables taking values in $\N$ and let $S_n=X_{n1}+ \ldots +X_{nn}$.  Put
\begin{eqnarray}\label{poisson.notation}
{\rm (a)} & & \hbox{$p_{ni}(k)={\mathbb P}\{X_{ni}=k\}$ \qq\quad $p_{ni} =1-p_{ni}(0)$,}
\cr &\cr {\rm (b)} & & \hbox{$\lambda_n=\sum_{i=1}^n p_{ni} $ \qq    ($i=1,\ldots, n$, $n\ge 1$, $k\ge 0$).}
\end{eqnarray}

 \vskip  8pt If for each $n$, $X_{ni}$ are i.i.d. Bernoulli random variables with \hbox{$p_n\! =\! {\mathbb P}\{X_{ni}\! =\! 1\}$}, then
 \begin{equation} \label{poisson} \lim_{n\to \infty} {\mathbb P}\{S_n=k\}\,=\, e^{-\lambda}\frac{\lambda^k}{k!}, \qq \quad k=0,1,\ldots,
\end{equation}
if $np_n\to \lambda$ as $n\to \infty$, $\lambda$ being a positive real.  As $S_n$ has the binomial distribution $\mathcal B(n,p_n)$, the distribution  $\mathcal B(n,p_n)$
tends to the Poisson distribution $\mathcal P(\lambda)$.
 The case $np_n=\lambda$ was established by Poisson.

   The convergence can be strengthened and in fact Wang \cite{Wan} showed the following generalized local limit theorem,   improving upon   Johnson and Simons 1971's analog
   result \cite{JS} for the subcase $np_n=\lambda$.

\begin{theorem}\label{wang}  Assume that
\beq \label{wang.assumption}
\aligned
{\rm (i)} & \quad \hbox{$\lambda_n=\lambda + \mathcal O (1/n)$ for some real $\lambda>0$,}
\cr
{\rm (ii)} & \quad \hbox{$\frac{p_{ni}(k)}{p_{ni}}=\a(k)$ are independent from $i$ and $n$ for all $k\ge 1$,}
\cr
{\rm (iii)} & \quad  \hbox{$\max_{1\le i\le n} p_{ni}\to 0$ as $n\to \infty.$}
\cr
\endaligned
\eeq
Then,
\begin{equation} \label{wang1} \lim_{n\to \infty} \sum_{k=0}^\infty h(k) \big|{\mathbb P}\{S_n=k\}-{\mathbb P}\{Y=k\}\big| \,=\, 0, \qq \quad k=0,1,\ldots,
\end{equation}
for all non-negative functions $h$ with $\sum_{k=0}^\infty   h(k){\mathbb P}\{Y=k\}<\infty$, where $Y$ has compound Poisson distribution with parameter $\lambda$ and
compounding distribution $\{ \a(k), k\ge 1\}$.
\end{theorem}

In Wang  \cite{Wan}, other modes of convergence are also considered and compared.

\vskip 4 pt
\noi {\it Le Cam's inequality.}
Let $S_n$ be the sum of $n$ independent binomial random variables  $X_i$ with  ${\mathbb P}\{X_i=1\}= p_i$ and set
$$ \lambda = p_1+\ldots+ p_n.$$
Le Cam \cite{LC} established in  1960 the following remarkable inequality

 \begin{theorem}  \label{lecam}
\begin{equation*}  \sum_{k=0}^\infty \Big|\,{\mathbb P}\{S_n=k\} -e^{-\lambda}\frac{\lambda^k}{k!} \,\Big|\, \le \, 2 \sum_{i=1}^n p_i^2.
\end{equation*}
\end{theorem}

 Taking $p_i=\lambda/n$ and noting that the right side simplifies to $2\lambda^2/n$, we get the classical Poisson limit \eqref{poisson} with a control of the error term.  Le
 Cam's inequality identifies the sum $\sum_{i=1}^n p_i^2$ as a quantity governing the quality of the Poisson approximation.

\bigskip
Using a coupling method, Hodges and Le Cam \cite{HLC} gave an elementary proof of this result. Coupling methods are  elegant and powerful, the Bernoulli extraction part of a
random variable can be ranged in this category of methods. Other powerful methods are the semi-group method and the Chen-Stein method. These are nicely described in several
papers, notably the one of Steele \cite{Ste}.

 \smallskip\par

Hodges and Le Cam proof goes as follows.
Let $Y_1, Y_2,  \ldots$ be independent Poisson random variables with ${\mathbb E\,} Y_i=p_i$ and chosen so that ${\mathbb P}\{X_i=Y_i\}$ is as large as possible. More
precisely
\begin{equation}\label{jointdistriBinomialPoisson}
\begin{cases}
 {\mathbb P}\{ X_i=Y_i=1\}&=\quad p_ie^{-p_i}\cr
 {\mathbb P}\{ X_i=1\,,\,Y_i=0\}&=\quad  p_i\big( 1-e^{-p_i}\big)\cr
 {\mathbb P}\{ X_i=Y_i=0\}&= \quad e^{-p_i} -p_i\big( 1-e^{-p_i}\big)\cr
 {\mathbb P}\{ X_i=0\,,\,Y_i=y\}&=\quad \frac{p_i^ye^{-p_i}}{y!}\qquad y=2,3,\ldots\cr
\end{cases}\end{equation}
It is necessary to have ${\mathbb P}\{ X_i=Y_i=0\}\ge 0$ which is fulfilled if $p_i\le 0,8$. Let $T_n= \sum_{i=1}^n Y_i $, $n\ge 1$. From the well-known additivity property
of Poisson distributions, $T_n$ has Poisson distribution  $\mathcal P( \sum_{i=1}^n \lambda_i)$.

Put
\begin{equation} \label{DBP}
D=\sup_{u}\, \big|{\mathbb P}\{S_n\le u\}-{\mathbb P}\{T_n\le u\}\big|.
 \end{equation}
 Then
\begin{equation*}
D \le 2 \sum_{i=1}^n p_i^2 .\end{equation*}

\vskip 4 pt
\noi {\it Measures of disparity.}  Several measures of disparity between the  distribution of two non-negative integer valued random variables $X$ and $Y$ are considered in
the literature. A first natural measure of disparity is
\begin{equation} \label{d-disparity}
d(X,Y) = \sup_{A} \big|{\mathbb P}\{X\in A\}-{\mathbb P}\{Y\in A\}\big| ,\end{equation}
 where the supremum is taken over all steps $A$ of non-negative integers.
 One verifies that $d(X,Y) $ can be alternatively written
 \begin{equation} \label{d-disparity.alternative}
d(X,Y) = \frac12 \, \sum_{k=0}^\infty \big|{\mathbb P}\{X=k\}-{\mathbb P}\{Y=k\}\big|.\end{equation}
Another standard disparity measure is
\begin{equation} \label{d.0-disparity}
d_0(X,Y) = \sup_{k\ge 0} \,\big|{\mathbb P}\{X=k\}-{\mathbb P}\{Y=k\}\big|.\end{equation}
Obviously $d_0(X,Y)\le d(X,Y)$.  It is easily seen that $d$, $d_0$ are in fact metrics in the space of all probability distributions on $\N$, in particular the triangle
inequality is satisfied. Simple but quite useful bounds for $d$ and $d_0$ are
\begin{eqnarray} \label{disparity-bound}
d(X,Y) &\le & {\mathbb P}\{X\neq Y\},
\cr &\cr d_0(X,Y) &\le & \max \big\{ {\mathbb P}\{X<Y\}, {\mathbb P}\{X>Y\}\big\}.
\end{eqnarray}

 Consider the case when $X$ is a sum of independent of non-negative integer valued random variables $X_i$, and $Y$ a suitably approximating Poisson variable. Le Cam's
 inequality can be reformulated as follows:
\vskip   3 pt {\it Let $X_1, \ldots, X_n$ be independent Bernoulli random variables with respective success probabilities $p_1, \ldots, p_n$ and let $Y$ be Poisson with mean
$\sum_{i=1}^n p_i$. Then}
\begin{eqnarray} \label{d-lecam}
d\Big(\sum_{i=1}^n X_i,Y\Big) &\le & \sum_{i=1}^n p_i^2.
\end{eqnarray}
We also state Franken's result \cite{Fr} concerning $d_0$.

 \vskip 2 pt {\it Let $X_1, \ldots, X_n$ be independent non-negative, square integrable integer-valued random variables,  and let $Y$ be Poisson with mean $\sum_{i=1}^n
 {\mathbb E\,} X_i$. Then}
\begin{eqnarray} \label{d.o-franken}
d_0\Big(\sum_{i=1}^n X_i,Y\Big) &\le & \frac{2}{\pi} \sum_{i=1}^n\big( {\mathbb E\,} X_i^2 + {\mathbb E\,} X_i(X_i-1)\big).
\end{eqnarray}

An extension of these results to the case of dependent non-negative integer-valued random variables was obtained in  Serfling \cite[Th.\,1]{Se2}.

}

{\color{green}\subsection{Some applications}}  \label{1.12}
The field of applications of the local limit theorem is not so well known from non-specialists, but  reveals to be impressive. We  only introduce, through  concrete examples
or problems, to some classes of applications, which are just  parts of  major domains of research.
\vskip 5 pt (1)
  {\it Asymptotic enumeration.}
Let $a_n(k)$ be a double sequence of non-negative real numbers. We associate the normalized double sequence
\begin{equation} \label{enum}
p_n(k)=\frac{a_n(k)}{\sum_ja_n(j)}.
\end{equation}
We say that $a_n(k)$ is asymptotically normal with mean $\m_n$ and variance $\s_n^2$ if
\begin{equation} \label{enumclt}
 \lim_{n\to \infty}\sup_{x\in \R}\Big|\sum_{k\le \s_nx+\m_n}p_n(k)-\frac{1}{\sqrt{2\pi}}\int_{-\infty}^xe^{-t^2/2} {\rm d} t \Big|=0.
 \end{equation}
Equivalently, we  say that the sequence $a_n(k)$ satisfies the central limit theorem.

\vskip   3 pt
If for some set of reals numbers $S$, {\color{green}there are numbers $\m_n$ and  $\s_n $ such that}
\begin{equation} \label{enumllt}
 \lim_{n\to \infty}\sup_{x{\color{green} \in S}}\Big|\s_n p_n([\s_nx+\m_n])-\frac{1}{\sqrt{2\pi}} e^{-x^2/2} {\rm d} t \Big|=0,
 \end{equation}
we also say that the sequence $a_n(k)$ satisfies the local limit theorem on $S$.
 \smallskip\par  The first results of this kind were   obtained  for Stirling numbers by Moser and Wyman \cite{MoserWyman}, Harper \cite{harper} and Lieb \cite{lieb}.
The existence of a  central limit theorem for a sequence $a_n(k)$  of numbers arising in enumeration, provides a qualitative feel for their behavior, but a local limit
theorem is better since it provides asymptotic information about $a_n(k)$ for $|k-\m_n|= \mathcal O(\s_n)$. However the step from \eqref{enumclt} to \eqref{enumllt} is often
not easy. It would become easier when in addition, some smoothness of $a_n(k)$ exists. In \cite[Th.\,1]{Be}, Bender proved the following central limit theorem based on the
nature of the generating function $\sum a_n(k) z^nw^k$, and gave several applications (to a dimer problem, to ordered sets of partitions, and Eulerian numbers).

\begin{theorem}\label{enumcltbender}
Let $f(z,w)$ have power series expansion
$$f(z,w)= \sum_{n,k\ge 0} a_n(k)z^nw^k $$
with non-negative coefficients. Suppose that:
\begin{eqnarray*}
\hbox{\rm (i)} & & \hbox{there exists an $A(s)$ continuous and non-zero near $0$,} \cr
\hbox{\rm (ii)} & & \hbox{an $r(s)$ with bounded third derivative near $0$,} \cr
\hbox{\rm (iii)} & & \hbox{a non-negative $m$, and $\e,\d>0$ such that} \cr
 & &\qq \qq \Big(1-\frac{z}{r(s)}\Big)^mf(z, e^s) -\frac{A(s)}{1- {z}{r(s)}}\  \cr
  & & \hbox{is analytic and bounded for}\ \ |s|<\e \quad |z|<|r(0)| + \d.
\end{eqnarray*}
Define
\begin{equation} \m=-\frac{r'(0)}{r(0)} ,\qq \s^2=\m^2- \frac{r''(0)}{r(0)}
\end{equation}
If $\s\neq 0$, then \eqref{enumclt} holds with $\m_n= n\m$ and $\s^2_n=n\s^2$.
\end{theorem}

Bender also indicated
 that certain smoothness conditions on the $a_n(k)$ are sufficient to conclude \eqref{enumclt} from \eqref{enumllt}.
 Recall that a sequence $a(k)$ is unimodal if it is first non-decreasing and then non-increasing. He showed in particular that if \eqref{enumclt} holds, $\s_n\to \infty$, and
 $p_n(k) $ are {\it unimodal}
for sufficiently large $n$, then for every $\e>0$, the local limit theorem holds on $S=\{x: |x|\ge \e\}$.
\vskip   3 pt Further if  $a_n(k) $ are {\it log concave}: $(a_n(k))^2\ge  a_n(k-1)a_n(k+1)$, then \eqref{enumllt} holds from \eqref{enumclt} with $S=\R$. Furthermore, the
rate of convergence depends only on $\s_n$ and the rate  of convergence in \eqref{enumclt}

\begin{remark} Canfield \cite[p.\,286]{Ca} has noted that concluding \eqref{enumclt} from \eqref{enumllt} requires in fact  that certain $a_n(k)$ are non-zero.
\end{remark}

 \smallskip\par  For the case when $P_n(x) = \sum_{k}a_n(k)x^k$ are polynomials whose roots are all reals and non-positive, Harper \cite{harper} showed that if $\s_n\to
 \infty$, then the local limit theorem holds with $S=\R$.  Bender  showed  how to recover simply Harper's result  from a general  central limit theorem and using   Newton's
 inequality to show that $p_n(k)$ are log concave.
 Applications  to Eulerian numbers, Stirling numbers of the first kind and of the second kind, separated samples, are also given.
 \smallskip\par
 Bender's  central and local limit theorems were since extended to  wider classes of generating functions, and to some general combinatorial structures, see for instance
 Flajolet and Soria \cite{FSo}, Gao and  Richmond \cite{GR} and the references therein.

\vskip 5pt (2)
  {\it Coefficients of polynomials of binomial type.}
A sequence of polynomials $(P_n)_{n\ge 0}$ satisfying
\begin{equation} \label{binom.pol}
P_n(x+y)= \sum_{k=0}^n {n\choose k} P_k(x)P_{n-k}(y),
\end{equation}
and having $P_0(1)=1$, is said to be of {\it binomial type}. For a sequence of polynomials to be of binomial type, it is necessary and sufficient that there exists a (formal)
power series lacking constant term such that
\begin{equation} \label{binom.pol.pow.ser}
\exp( xg(u))= \sum_{n=0}^\infty P_n(x) \big(\frac{u^n}{n!}\big).
\end{equation}
Expressing $g(u)$ as an exponential generating function
\begin{equation} \label{expo.gen.funct.}
g(u)=  \sum_{n=0}^\infty g_n \big(\frac{u^n}{n!}\big),
\end{equation}
expansion \eqref{binom.pol.pow.ser} reveals the following identity for the coefficients $s(n,k) $ of $P_n(x)$
\begin{equation} \label{snk}
s(n,k)= \big(\frac{1}{k!}\big)\sum {n\choose v_1,\ldots, v_k} g_1\ldots g_k.
\end{equation}

Canfield \cite[Th.\,II]{Ca} proved the following local limit theorem.

\begin{theorem}\label{lltbinpol} Let $s(n,k)$ be the coefficients of polynomials of binomial type (it is not assumed that $g(u)$ is a polynomial). Assume  that the $s(n,k)$
satisfy the central limit theorem \eqref{enumclt}  for some numbers $\m_n$ and $\s_n\to \infty$, and that the following conditions are satisfied:
\begin{eqnarray*}\label{lltbinpolcond}
\hbox{\rm (i)}& &\hbox{$s(n,k) $ are log concave: $(s(n,k))^2\ge  s(n,k-1)s(n,k+1)$. }\cr
\hbox{\rm (ii)}& &\hbox{For each $n$, $\{k: s(n,k)\neq 0\}$ is  a bounded  interval}.\end{eqnarray*}
Then the $s(n,k)$ satisfy the local limit theorem with $S=\R$. Consequently,
\begin{equation}s(n,\m_n+x\s_n) \sim \frac{P_n(1)}{\s_n \sqrt{2\pi}}\, e^{-x^2/2},\qq\hbox{as $n\to \infty$,}
\end{equation}
where $x= \mathcal O(1)$.\end{theorem}
\begin{remark} We  note that condition (ii) excludes from the scope of the above Theorem the important case when the set $\{ k: s(n,k)= 0\} $ is infinite.
\end{remark}
{\color{red}{
\vskip 5 pt (3)
  {\it Asymptotics for  Legendre polynomials.} Let $\{\b_k\}$ be i.i.d. Bernoulli random variables with ${\mathbb P}\{\b_k=1\}=p$ and $\{\tilde{\b_k}\}$
  be an indepe{\color{green}n}dent copy of $\{\b_k\}$. Set ${\mathcal B_n} = \b_1+\ldots+\b_n$ and  $\tilde{\mathcal B}_n = \tilde{\b_1}+\ldots+\tilde{\b_n}$. Thus
\begin{equation} \label{binsym}
{\mathbb P}\{{\mathcal B}_n -  \tilde{\mathcal B}_n = 0\}=\sum_{k=0}^n {\mathbb P}\{\mathcal B_n=k\} {\mathbb P}\{-\tilde{\mathcal B}_n=-k\}=
\sum_{k=0}^n b^2(k,n,p).
\end{equation}
On the other hand by Exercise 11 in Ch. 2, \S5 in \cite{Ri68} we have
\[{\mathcal P_n}(1+2x) = \sum_{\nu=0}^n{{n}\choose{\nu}}^2x^{\nu}(1+x)^{n-\nu},\]
where ${\mathcal P_n}$ are  Legendre polynomials. If we substitute $x\rightarrow{{x}\over{1-x}}$ then we get
\[
(1-x)^n{\mathcal P_n}({{1+x}\over{1-x}}) =  \sum_{\nu=0}^n {{n}\choose{\nu}}^2
x^{\nu}.
\]
Finally we put in the latter $x={{p^2}\over{(1-p)^2}},$
$p\in(0,{{1}\over{2}})$ so that
\begin{equation} \label{polyleg}
(1-2p)^n{\mathcal P_n}({{1-2p+2p^2}\over{1-2p}})
=\sum_{k=0}^n b^2(k,n,p) = {\mathbb P}\{{\mathcal B_n} - \tilde{\mathcal B}_n = 0\}.
\end{equation}
For $p={{1-x+\sqrt{x^2-1}}\over{2}}$ in (\ref{polyleg}) we have $\sigma^2= 1-x^2+x\sqrt{x^2-1}$.
Therefore, by Theorem \ref{gnedenko} we obtain the Laplace-Heine theorem (see Theorem 8.21.1 in \cite{Sz39})
\begin{theorem}\label{laphei} For $x>1$ we have
${\mathcal P_n}(x) =
{{1}\over{\sqrt{2\pi n}}}(x+\sqrt{x^2-1})^n\sqrt{{x+\sqrt{x^2-1}}\over{\sqrt{x^2-1}}}(1+o(1)).$
\end{theorem}
Similarly one can use asymptotic expansions (see  Theorem 1 in \S51, \cite{GK}) to obtain the following
\begin{theorem}\label{laphei.a}
For $x>1$ we have
\[{\mathcal P_n}(x) =
{{1}\over{\sqrt{2\pi n}}}(x+\sqrt{x^2-1})^n\sqrt{{x+\sqrt{x^2-1}}\over{\sqrt{x^2-1}}}\left(1
+ {{1}\over{8n}}{{x-2\sqrt{x^2-1}}\over{\sqrt{x^2-1}}} + o\left({{1}\over{n}}\right)\right).
\]
\end{theorem}
This approach was used by Maejima and Van Assche to obtain asymptotics for more general Jacobi polynomials and other
orthogonal polynomials (see \cite{mae85}).
}}

\vskip 5 pt (4)
  {\it Allocation   problems.} Allocation (occupancy) problems   originates in statistical physics (Maxwell--Boltzmann statistics). In the 1930's nonparametric statistical
  tests led to an interest in allocation problems.
Consider for instance the following model: $n$ balls are distributed in $m$ boxes in such a way that all the $m^n$  arrangements are equally probable. Denote by $\m_r(m,n)$
($r=0,1,\ldots, n)$ the number of boxes containing exactly $r$ balls.

Sevastyanov and  Chistyakov \cite{ChSe} proved that if $n = [ma]$, $0< a<\infty$, and
$m\to\infty$, then  for every $s = 0,1,\ldots$   the random variables $\m_r(m,n)$ ($0 < r < s$) have an asymptotic $(s +1)$-dimensional normal distribution.
Kolchin \cite{Kolc1} later showed   the following result: If $n = [ma]$ where $0< a<\infty$,
then for any $r = 0,1,\ldots$ we have
\begin{equation} \lim_{m\to\infty}\s_r\sqrt m \Big| {\mathbb P}\{\m_r(m,n)=j\} -\phi\Big(\frac{j-mp_r}{\s_r\sqrt m}\Big)\Big| =0
\end{equation}
uniformly in $j$ for $|j-mp_r|< K\sqrt m$, $0 < K < \infty$, where $\phi(x)$ is the normal
density function, $p_r=e^{-\a}\a^r/r!$ and
$$\s_r^2=p_r\big( 1-p_r-\frac{p_r}{\a}(r-\a)^2\big) .$$
The proof is based  on Gnedenko's theorem.
\vskip 5 pt (5)
  {\it Random permutations.} Among the $n!$ permutations of $(1,2,\ldots,n)$ a permutation is chosen at random in such a way that all the $n!$ permutations are equally
  probable. Denote by $a(n)$ the number of cycles in this permutation. Then we have
\begin{equation}    {\mathbb P}\{a(n)=j\}   = S(n,j)/n!\qq  (j = 1,2,\ldots ,n)
\end{equation}   where $S(n,j)$,     $j = 1,2,\ldots ,n$  is a Stirling number of the first kind.   Goncharov proved that
\begin{equation}\lim_{n\to \infty}{\mathbb P}\Big\{\frac{\a(n)-\log n}{\sqrt {\log n}}\le x\Big\} =\Phi(x) =\frac{1}{\sqrt{2\pi}}e^{-u^2/2} {\rm d} u.
\end{equation}
   Kolchin \cite{Kolc2}   proved the corresponding local limit theorem, that is
\begin{equation}\label{lltrand.perm}\lim_{n\to \infty} \Big|\sqrt{\log n}{\mathbb P}\{a(n)=j\}-\phi\Big(\frac{j-\log n}{\sqrt {\log n}}\Big)\Big| =0,
 \end{equation}
uniformly with respect to $j$ in the interval $|j-\log n|<K (\log n)^{7/ 12}$. This result provides an asymptotic formula for $S(n,j )$, namely,
\begin{equation} S(n,j)\sim \frac{n!}{\sqrt {\log n}}\phi\Big(\frac{j-\log n}{\sqrt {\log n}}\Big)
\end{equation}
as $n\to \infty$ if $|j-\log n|<K (\log n)^{7/ 12}$.   Golomb \cite{Go} observed that if $L_n$  is the expected length of
the  longest    cycle   in  a   random  permutation of  $(1,2,\ldots ,n)$,  then    the     limit $\lim_{n\to \infty}\frac{L_n}{n}= \lambda$ exists.   Shepp and  Lloyd
\cite{LS} proved that
\begin{equation}
 \lambda =  \int_0^\infty\exp\Big\{-x-\int_x^\infty \frac{e^{-y}}{y}{\rm d} y\Big\} {\rm d} x.
\end{equation}

  Let {\color{green}$S_n$ denote the symmetry group}. The order of a permutation $g\in S_n$ is  the smallest positive integer $m$ for which $g^m=e$, the identity permutation.
  Let us choose a permutation at random among the $n!$ permutations of $S_n$ so that all the
 $n!$ permutations are equally probable, and let  $\nu_n$   denotes the order of the random permutation. Erd\H os and P. Tur\'an \cite{ET} showed that
\begin{equation}
\lim_{n\to \infty}{\mathbb P}\Big\{\frac{\log \nu_n-(\log n)^2/2}{\sqrt{(\log n)^3/3}}\Big\} = \Phi(x).
\end{equation}
 See for instance the monograph by Arratia, Barbour and Tavar\'e \cite{ABT2} for recent developments.

\vskip 5 pt (6)
  {\it Random mappings.} This research's area covers allocation theory, random permutations, random mappings, branching processes, random trees, and random forests. We refer
  to the book by Kolchin   and to   the one  by   Kolchin,  Sevastyanov, and   Chistyakov \cite{KSC}.
The main  tool used in   Kolchin's book  is Gnedenko's local limit theorem, and some of its extensions.

\vskip 5 pt
{\color{blue} (7)
  {\it Number theory.}
 Applications of the local limit theorem to additive problems are given in Postnikov's monograph \cite{Po}, Section 2.3. In \cite{Pos}, see also \cite{Po}, Section 2.4,
 asymptotic estimates   for the number of solutions of the diophantine equations
\ben
{\rm (E1)}& & x_1+\ldots +x_n= N
\cr {\rm (E2)}& & x^2_1+\ldots +x^2_n= N
\een
are obtained applying the local limit theorem to related sums of independent random variables and estimating resulting trigonometric sums. Later Ismatullaev and Zuparov
\cite{IsZu}
  considered   the number of solutions of the diophantine equation
 \ben
  {\rm (E3)}& & x^s_1+\ldots +x^s_n= N
\een
in which $n,N,s\ge 1$, and   sharpened and generalized Postnikov's result by using the local limit theorem.

\vskip 2 pt
{\color{red} Kubilius (\cite{Kub64}), \cite{Kub69}) proved local limit theorems for additive arithmetical functions.}
Local limit theorems of large deviations for   additive functions are proved in Manstavichyus  and  Skrabutenas \cite{MansSkra}. Let  $h(m)$ be additive integer-valued
function. The letter $p$ denoting a prime number,
{\color{red} let $\m_n=B_n^{-1}\max_{p\le n}  |h(p)|$,}
where $B_n^2 =\sum_{p\le n}p^{-1} h^2(p)\to \infty$  and $\m_n\to 0$ as $n\to \infty$.   An asymptotic formula {\color{red} for relative
error
$$R_n=\frac{B_n}{n}\big| \{m\le n: h(m) =k\} \big| \phi^{-1}(y_{nk}),$$
where $y_{nk}={{k-A_n}\over{B_n}},$ $A_n=\sum_{p\le n}p^{-1} h(p)$}
is obtained in terms of the Cram\'er series (see Theorem \ref{ldt2}), and which is valid uniformly in $k$ in the region
$|k-A_n|=o(\m_n^{-1})$.
\vskip 2 pt

In Stepanauskas \cite{Step}, additive functions of the form $f(n)=\o(n)+\rho(n)$, where $\o(n)$ counts the number of   prime divisors of $n$, and $\rho(n)$ is an
integer-valued additive function  {\color{green}obeying} to some strong growth assumptions. Two local limit theorems are proved, providing
  formulae for $|\{n\le x f(n)=k\}|$ of Sathe-Selberg type, with typical ranges for $k$ of the form $c_1 \log \log x\le k\le  c_2 \log \log x$, where $c_1, c_2$ are positive
  constants.

We also cite Tulyaganov \cite{Tuly}, where a local limit theorem for multiplicative arithmetic functions is obtained.
}}

{\color{red}\section{Almost sure versions} \label{2}

 {\color{blue}   In the following studies, we obtain almost sure local limit theorems  in various different but fundamental cases. We essentially follow the classical
 approach and first proceed with a second order theory  of  the underlying  system of random variables. This  is an important and very informative part of the  study. Next
 classical convergence criteria can be applied to derive almost sure limit theorems.   It turns out the study of the correlation function of these systems is each time hard.
}
\vskip 20 pt

The Theorem below was proved by Chung and Erd\H{o}s \cite[Theorem 6]{CE}, and can be viewed as a primary form of  an almost sure local limit theorem, see \cite{DK}.

\begin{theorem}\label{CE1}
Let $X_1, X_2, \ldots$ be a sequence of i.i.d.\ integer valued random variables with $\E X_1=0$ and put $S_k=\sum_{i=1}^k X_i$.
Assume that every integer $a$ is a possible value of $S_k$  for all sufficiently large k. Then for every integer $a$
\begin{equation}\label{CE2}
\lim_{n \to \infty} \frac1{\log M_n} \sum\limits_{k = 1}^n \frac{I\{S_k=a\}}{M_k} = 1 \quad \text{a.s.,}
\end{equation}
where
$$M_k=\sum_{i=1}^k \P (S_i=a).$$
\end{theorem}

\bigskip
Note that under the assumptions of Theorem \ref{CE1} we have $\lim_{k\to\infty} M_k=\infty$ and also
$$ \sum_{k=1}^n \frac{\P(S_k=0)}{M_k} = \sum_{k=1}^n \frac{m_k}{M_k} \sim \log M_n \quad \text{as} \ n\to\infty.$$
Hence setting $\xi_k= I \{S_k=a\}/\P(S_k=a)$ and
$$ d_k= m_k/M_k, \quad D_n=\sum_{k=1}^n D_k,  $$
relation (\ref{CE2}) can be written equivalently as
\begin{equation}\label{CE3}
\lim_{n \to \infty} \frac1{\log D_n} \sum\limits_{k = 1}^n d_k \xi_k = 1 \quad \text{a.s.}
\end{equation}
stating that the logarithmic averages of $\xi_k$ converge a.s.\ to 1. Note that $\xi_k$ is a strongly dependent sequence of random variables and the ordinary strong law
$$ \lim_{n\to\infty} \frac{1}{M_n} \sum_{k=1}^n m_k \xi_k=1 \quad \text{a.s.} $$
is not valid. Apart from an unproved remark in L\'evy \cite[p.\ 270]{Le37}, Theorem \ref{CE1} is the starting point of a recent and highly interesting theory dealing with
almost sure versions of classical weak limit theorems involving logarithmic averaging.

other results of \cite{CE}, can easily be corrected for symmetric $X_1$.

{\color{green}\subsection{ASLLT for i.i.d.\ square integrable lattice  random variables}}\label{2.1}

Let $X_1, X_2, \ldots$ be i.i.d.\ lattice random variables with maximal span 1, mean 0 and variance $\sigma^2$ and put $S_n=\sum_{k=1}^n X_k$. Then $\P (S_n=0)\sim
\frac{1}{\sigma \sqrt{2\pi n}}$ and thus by Theorem \ref{CE1}
we have
\begin{equation}\label{CE3.a}
\lim_{n \to \infty} \frac1{\log n} \sum\limits_{k = 1}^n \frac{I\{S_k=0\}}{\sqrt{k}} =   \frac{1}{\sigma} \sqrt{\frac{2}{\pi}}  \quad \text{a.s.}
\end{equation}
Under assuming finite variances, the following theorem yields a much more general ASLLT for i.i.d.\ lattice valued random variables.
\begin{theorem} \label{t1[asllt]} Let $X$ be a square integrable  random variable taking values on the lattice $\mathcal L(v_0,D)= \{v_0+kD, k\in \mathbb Z\}$  with maximal
span $D$. Let $\m ={\mathbb E\,} X$,
$\s^2={\rm Var} (X)$ which we assume to be positive
 (otherwise $X$ is degenerate).  Let also $ \{X_k, k\ge 1\}$ be independent copies of
$X$, and put
$S_n=X_1+\ldots +X_n$, $n\ge 1$.   Then
$$ \lim_{ N\to \infty}{1\over    \log N } \sum_{ n\le
N}  {  1 \over \sqrt n} {\bf 1}_{\{S_n=\kappa_n\}} \buildrel{a.s.}\over {=}{D\over
\sqrt{ 2\pi}\s}e^{-  {\k^2/ ( 2\s^2 ) } },$$
  for any  sequence of integers $\k_n\in \mathcal L(nv_0,D)$, $n=1,2,\ldots$      such that
  \begin{equation}\label{eq2}  \lim_{n\to \infty} { \k_n-n\m   \over    \sqrt{  n}  }= \k .
\end{equation}
 \end{theorem}

Concerning the sequence $\k_n$,  note that by  Gnedenko's Theorem \ref{gnedenko},  if
$ \k_n \in \mathcal L(nv_0,D)$ is a sequence which   verifies condition (\ref{eq2}), then
 \begin{equation}\label{llt1a}  \lim_{n\to \infty}  \sqrt n {\mathbb P}\{S_n=\k_n\}={D\over  \sqrt{ 2\pi}\s}e^{-
{ \k ^2\over  2   \s^2} } .
\end{equation}

Theorem \ref{t1[asllt]} was announced  in Denker and Koch \cite{DK} (Corollary 2). The proof however contains a gap. A complete proof was given in Weber \cite{W}. The
argument used depends on sharp inequalities  for the correlations
$$|\P (S_m=k_m, S_n= k_n)-\P (S_m=k_m)\P (S_n=k_n)|$$
which are also established in \cite{W}, and   are much harder to get than the correlation inequalities for
$$|\P (S_m<k_m, S_n<k_n)-\P (S_m<k_m)\P (S_n<k_n)|$$
needed for the proof of global a.s.\ central limit theorems, see e.g.\ Lacey  and Philipp \cite{LP}.}

\vskip 5 pt\noi  {\tf The correlation inequality.} This one  states as follows.


{\color{blue}
\begin{theorem} \label{t2[asllt]} Assume that
  \begin{equation}\label{basber}\hbox{ $ {\mathbb P}\{X=k\}\wedge{\mathbb P}\{X=k+1\}\,>0$  \qquad  for some   $k\in\Z$.}
 \end{equation}  Then there  exists a constant
$C $,
\begin{eqnarray} \label{basber.est}& &\sqrt{nm} \, \Big|{\mathbb P}\{S_n=\k_n, S_m=\k_m\}-{\mathbb P}\{S_n=\k_n \}{\mathbb P}\{  S_m=\k_m\} \Big|
\cr  &\le &  C \Big\{     {   1\over
 \sqrt{n\over m}-1   } +      {n^{1/2}  \over
 (n-m) ^{3/2}} \Big\} .
\end{eqnarray}
 \end{theorem}

 \begin{corollary} \label{cort1}      Let $0<c<1$. Under assumption \eqref{basber},  there   exists a constant
$C_{  c} $  such that for all $1\le m\le cn$,
\begin{eqnarray*} \sqrt{nm} \,\Big|{\mathbb P}\{S_n=\k_n, S_m=\k_m\}-{\mathbb P}\{S_n=\k_n \}{\mathbb P}\{  S_m=\k_m\} \Big|
 &\le &  C_{c} \, \sqrt{{m\over n}}. \end{eqnarray*}
 \end{corollary}

\begin{remark}\label{bbl}\rm
Condition (\ref{basber}) seems to be somehow artificial. It is for instance clearly not   satisfied if  ${\mathbb P}\{ X\in
{\mathcal N}\} =1$   where
${\mathcal N}=\{\nu_j, j\ge 1\}$ is an increasing sequence of integers such that $\nu_{j+1}-\nu_j >1$ for all $j$.
  This already defines  a large class of examples. However, condition (\ref{basber}) is   natural in our setting.
 By the local limit theorem (\ref{llt}), under condition (\ref{eq2}),
   $$ \lim_{n\to \infty}\sqrt n {\mathbb P}\{S_n=\ell_n\}= {D\over  \sqrt{ 2\pi}\s}e^{-
{\k^2\over  2  \s^2} }, \qq (\ell_n\equiv \k_n\ {\rm or}\ \ell_n\equiv\k_n+1)  .
$$
Then for some $n_\k<\infty$, ${\mathbb P}\{S_n=\k_n\}\wedge {\mathbb P}\{S_n=\k_n+1\}>0 $ if $n\ge n_\k$. Changing $X$ for   $X'=S_{n_\k}$, we see that $X'$ satisfies
(\ref{basber}).
\end{remark}

When  ${\mathbb E\,} |X|^{2+\e}<\infty$
for some positive $\e$,
 a
similar result can be  obtained  more easily.
\begin{eqnarray}\label{GW} & &\sqrt{nm} \,\Big|{\mathbb P}\{S_n=\k_n, S_m=\k_m\}-{\mathbb P}\{S_n=\k_n \}{\mathbb P}\{  S_m=\k_m\}
\Big|
\cr
 &\le &   C \Big({{1} \over {\sqrt{{n\over m}}-1}}+ \sqrt {n\over n-m}\ {{1} \over {(n-m)^\alpha}}\Big), \end{eqnarray}
with  $\a=\e/2$. See \cite{GW2}. Condition (\ref{basber}) is  not needed.
 Corollary \ref{cort1}  follows in that case directly from
(\ref{GW}).   When $\e=0$,   this   can obviously no longer be applied, and another approach has to be implemented. Notice that even
when
$\e>0$, \eqref{basber.est} is  stronger than \eqref{GW}.


\vskip 20 pt
{\color{green}\subsection{ASLLT for   i.i.d. stable lattice random variables}}
\label{2.2}


{\color{green} We refer to   Giuliano, Szewczak
\cite{GS}}. Let $G$ be  a stable  distribution with exponent $\alpha$ ($0<\alpha \le
2$, $\alpha \not =1$) and denote by $g$ the $\alpha$--stable density function related to the distribution function $G$.  Let $(X_n)_{n \ge   1}$ be a sequence of i.i.d.
random
variables   such that their common distribution $F$ is in the domain
of attraction of $G$. This means that, for a suitable choice of constants $a_n$ and
$b_n$, the distribution of
$$T_n := \frac {X_1 + \cdots X_n -a_n }{b_n}$$
converges weakly to $G$.  It is well known (see \cite{IBLIN}, p.\ 46)
that in such a case we have $b_n = L(n) n^{1/\alpha}$, where $L$ is
slowly varying in Karamata's sense. For $\alpha > 1$ we shall assume that $X_1$ is centered; by Remark 2 p.\,402 of \cite{AD}, this implies that $a_n =0$, for every $\alpha$.
Let $\phi$  be the characteristic function of $F$.
By \cite{AD}, Theorem 1, for $\alpha \not =1$ it  has the form
\begin{equation}\label{funzionecaratteristica}
f(t) = \exp\left\{-c|t|^\alpha h\Big(\frac{1}{|t|}\Big)(1-i\beta {\rm sign}(t)  \tan\frac{\pi\alpha}{2})+ o\left(|t|^\alpha h\Big(\frac{1}{|t|}\Big)\right)\right\},
\end{equation}
where $c= \Gamma(1-\alpha)(c_1+c_2) \cos \frac{\pi\alpha}{2}>0$ and $\beta =\frac{c_1-c_2}{c_1+c_2}\in[-1,1]$
are two constants and $h(x)=l(x)$ if $\alpha\in(0,2)$ and $c=\frac{1}{2},$ $\beta=0,$ $h(x)=E[X^2 1_{\{|X|\le  x\}}]$
if $\alpha=2$.

 \smallskip\par

 We assume that $X_1$ takes values in the lattice $\mathcal{L}(v_0,D)= \{v_0+ mD, \, m \in \mathbb{Z}\}$ where $D$ is the maximal span of the distribution;
 hence  $S_n := X_1 + \cdots X_n $ takes values in the lattice $\mathcal{L}(nv_0,D)= \{nv_0+ mD, \, m \in \mathbb{Z}\}.$

 \smallskip\par
  We further assume that
\begin{equation}
\label{e2a}
x^\alpha {\mathbb P}\{X>x\} = \big(c_1 +o(1)\big)l(x); \qquad x^\alpha
{\mathbb P}\{X\le  -x\} = \big(c_2 +o(1)\big)l(x),\quad
\alpha\in(0,2],
\end{equation}
 where $l$ is slowly varying as $x
\to \infty$ and $c_1$ and $c_2$ are two suitable non--negative
constants, $c_1+c_2>0$, related to the stable distribution $G$.

\bigskip Put for $\e>0$,
\begin{equation}\label{M} M(x)=M_\e(x) = \sup_{\frac{1}{\e  }\le  y \le  x}h(y) , \qquad x \ge   \frac{1}{\e  }
\end{equation}
where  $h(x)=l(x)$ if $\alpha\in(0,2)$ and  $h(x)=E[X^2 1_{\{|X|\le  x\}}]$ if  $\alpha=2$.
We have $M_\e(x) < \infty$ for every $x$.

\begin{theorem}\label{princ} Let $(X_n)_{n \ge  1}$ be
a centered, independent and identically lattice distributed
(i.i.l.d.) random sequence with span $D=1$. Assume that
\eqref{e2a} holds with $\alpha\in (1,2]$ and that there exists
$\gamma \in (0,2)$ such that for some $\e>0$
\begin{equation}\label{princ.a1}\sum_{k=a}^b \frac { L(k)\big\{1+ M\big(k^{1+\frac{1}{\alpha}}\big)+L^\eta (k)\big\}}{k} \le  C(\log^\gamma b-
\log^\gamma a),
\end{equation}
for some $\eta \in (0,1]$. Further assume    that the function $h$
 appearing in   \eqref{funzionecaratteristica}
 verifies
\begin{equation}\label{princ.a2}\liminf_{x\to \infty}h(x)=:\ell>0.
\end{equation}

Then
$$\lim_{N \to \infty}\frac{1}{\log
N}\sum_{n=1}^N \frac{b_n}{n}1_{\{S_n = \kappa_n\}}= g(\kappa).
$$
\end{theorem}

\begin{example}\rm  Let $h(x) = \log^\sigma x$, with $0<\sigma< \frac{\alpha}{1+\alpha}$. Notice that
$L(n)= \frac{b_n}{n^{\frac{1}{\alpha}}}$
for sufficiently large $n$. It is not difficult to check that (\ref{princ.a1}) holds (see \cite[Example 4.2]{GS}) for details).
\end{example}

If $L \equiv$ a constant  (i.e. $F$ belongs to the domain of {\it normal}  attraction of $G$, according to the definition on p.\ 92 of \cite{IBLIN}), then the assumptions in
Theorem \ref{princ}   are automatically satisfied. Whence,

\begin{corollary}\label{carino}
If $(X_n)_{n \ge  1}$ is a centered i.i.l.d. random sequence with span $D=1$, \eqref{e2a} holds
with $\alpha\in (1,2)$ and $L \equiv c$, then $$\lim_{N \to \infty}\frac{1}{\log
N}\sum_{n=1}^N \frac{c}{n^{1- \frac{1}{\alpha}}}1_{\{S_n = \kappa_n\}}=
g(\kappa).$$
\end{corollary}

Again as in Section \ref{2.1} the key ingredient for the proof is an efficient correlation inequality which is also delicate to establish.
This time we use properties of characteristic functions of distributions in the domain of attraction of stable laws.

\subsubsection{Preliminaries} \bigskip\noindent
Formula (\ref{funzionecaratteristica}) implies that
\begin{align}\label{logaritmo}
\log \big| f(t)\big|& = \mathfrak{Re}\big(\log f(t)\big)= -c|t|^\alpha h\big(\frac{1}{|t|}\big)\big(1+o(1)\big),
\\ \arg\big( f(t)\big)& =\mathfrak{Im}\big(\log f(t)\big)= -c|t|^\alpha h\big(\frac{1}{|t|}\big)\big(-\beta{\rm sign}(t)\tan \frac {\pi \alpha}{2}  +o(1)\,\big)\nonumber
 \end{align}
hence
\begin{align} \label{argomentosulogaritmo}\lim_{t \to 0}\Big|\frac {\arg \big(f(t)\big)}{\log|f(t)|}\Big|=\Big|\beta \tan \frac {\pi \alpha}{2}\Big|.
 \end{align}
We notice that $L(n)= h^{\frac{1}{\alpha}}(b_n)$ for $\alpha\in(0,2)$ (by Remark 2 p.\,402 in \cite{AD}), while
  $L(n)=\sqrt{E\big[X^21_{\{|X|\le  b_n\}}\big]}$ for $\alpha=2$.

  \begin{remark}\label{equivalenza} \rm Let $\widetilde h\sim h$ as $x \to + \infty$. Then, by \eqref{logaritmo},
 \[ \log \big|f(t)\big|
 =  -c|t|^\alpha \widetilde h\big(\frac{1}{|t|}\big)\cdot \frac{h\big(\frac{1}{|t|}\big)}{\widetilde h\big(\frac{1}{|t|}\big)}\big(1+o(1)\big)
 =-c|t|^\alpha \widetilde h\big(\frac{1}{|t|}\big)\big(1+o(1)\big).\]
This means that $h$ is unique up to equivalence; thus, by Theorem 1.3.3. p.\ 14 of \cite{BGT} we can assume that $h$ is continuous (even $C^\infty$) on $[a, \infty)$ for some
$a>0$.

 \smallskip\par
An analogous observation is in force for $\arg\big(f(t)\big)$.
\end{remark}

\begin{remark}
 \rm Thus we deal with a subclass of strictly stable distributions.
Denoting by $\psi$ the characteristic function of $G$, we know from \cite{Z},
Theorem C.4 on p.\ 17 that  $\log \psi$ (for strictly
stable distributions) has the form
$$\log \psi(t)= - c |t|^\alpha \exp \{-i\left(\frac{\pi}{2}\right) \theta
\alpha\,{\rm sign} (t)\},$$
where $|\theta| \le  \min \{1, \frac{2}{\alpha}-1\}$ and $c >0$. For $\alpha=1$
and $|\theta|=1$ we get degenerate
distribution and in this case we say that $X_n$ is relatively stable (see e.g.
\cite{SZ}). Almost sure variant of
relative stability for dependent strictly stationary sequences will be discussed
elsewhere.\end{remark}

\medskip
Observe that Theorem 4.2.1 p.\,121 in \cite{IBLIN} implies that
\begin{equation}\label{qq}\sup_n\big\{\sup_m b_n {\mathbb P}\{S_n= m\}\big\}= C < \infty.
\end{equation}
For every $n$, let $\kappa_n$ be a number of the form $nv_0
+ mD$ and let
$$\lim_{n \to \infty}{{\kappa_n}\over{b_n}}= \kappa.$$

\subsubsection{Correlation inequality}
\vskip 1pt \noindent
We assume that $(X_n)_{n \ge  1}$ is a sequence of i.i.d. random variables verifying the following conditions:
(\ref{e2a}), $\alpha \not =1$ and
$\mu=E[X_1]=0$ when $\alpha >1$. Recall that the norming constant are $a_n=0$ and $b_n= L(n) n^{1/\alpha}$ with $L$ slowly varying. With no loss of generality, we shall
assume throughout that
$d=1$.

\begin{theorem}\label{cov} {\rm (i)} In the above setting we have
\[\Big|{\mathbb P}\{S_m=\kappa_m,S_n=\kappa_n \}- {\mathbb P}\{S_m=\kappa_m\}{\mathbb P}\{S_n=\kappa_n\}\Big|
\le  C\Big\{\Big(\frac
{n}{n-m}\Big)^{1/\alpha}\frac {L(n)}{L(n-m)}+1\Big\}.\]

{\rm (ii)} Assume that condition \eqref{princ.a2} is fulfilled.
 Then there exists $\e   >0$ such that if $M=M_\e$,     $M_\e$ being defined in \eqref{M}  we have
\begin {align}\label{secondacorrelazione}\nonumber & b_mb_n
\Big|{\mathbb P}\{S_m=\kappa_m,S_n=\kappa_n\})- {\mathbb P}\{S_m=\kappa_m
\}{\mathbb P}\{S_n=\kappa_n\}\Big|\\&\le  CL(n)\left\{n^{1/\alpha}\Big(\frac {1}{e^{(n-m)c} } + \frac {1}{e^{nc}
}\Big)+ \frac{\frac{m}{n}}{\Big(1- \frac{m}{n}\Big)^{1+ \frac{1}{\alpha}}}\Big(1+ M(n^{1+\frac{1}{\alpha}})\Big) + \frac{\big(\frac{m}{n}\big)^\frac{\eta}{\alpha}
L^\eta(m)}{(1-\frac{m}{n})^\frac{\eta +1}{\alpha}}  \right\}\end {align}
for every pair $(m,n)$ of integers, with $m \ge  1 $, $n > m+ \e   ^{-\frac{\alpha}{\alpha +1}}$, and for every $\eta \in (0,1]$, where $M=M_\e$ and   $M_\e$ is defined in
\eqref{M}.
\end{theorem}

\begin{remark}\label{nuovo}\rm
If $h$ is ultimately increasing, then condition
\eqref{princ.a2} is  satisfied.  A quick look at the proof (see {\color{orange} \cite{GS}}) shows that if $h$ is increasing and continuous, the inequality
\eqref{secondacorrelazione} holds for $1 \le  m<n$.
\end{remark}
Let $\delta > -1$ and $p>0$ two given numbers. We shall use the equality
\begin{equation}\label{formula}
\int_0^{+ \infty} t^\delta e^{-p t^\alpha} dt = \frac{\Gamma(\frac{\delta +1}{\alpha})}{\alpha}\cdot \frac{1}{p^{\frac{\delta +1}{\alpha}}}= C  \cdot\frac{1}{p^{\frac{\delta
+1}{\alpha}}}.
\end{equation}

If $L$ is non-decreasing the proof of (ii) in Theorem \ref{cov} (ii) simplifies.

\begin{corollary}\label{cor} For large $m$ and $n \geq 2m$, for every $\delta < \frac{1}{\alpha}$ and for every $\eta \in(0,1]$ we have $$b_mb_n
\Big|{\mathbb P}\{S_m=\kappa_m,S_n=\kappa_n \}- {\mathbb P}\{S_m=\kappa_m\}{\mathbb P}\{S_n=\kappa_n\}\Big| \le  C \tilde L(n)\cdot \Big(\frac{m}{n}\Big)^\rho,$$
with $\tilde L(n) =L(n)\big(1+ L(n^{1+\frac{1}{\alpha}})+L^\eta(n)\big)$ and $\rho:= \min\{\eta (\frac{1}{\alpha}- \delta),1\}$.
\end {corollary}

\begin{remark}\label{remW}
 \rm In the case $\alpha =2$ the correlation inequality of Corollary \ref{cor} furnishes $\rho= \eta\big(\frac{1}{2}-\e  \big)<  \frac{1}{2}$,
 while in Corollary \ref{cort1} the better exponent $\rho =\frac{1}{2}$ is found.
 \end{remark}
{
\begin{remark}\label{remZS}
\rm In the case $\alpha =2$ with $E(X_1^2)<\infty$ we can avoid (\ref{e2a})
(see \cite[Remark 2.1]{GS})
if we adopt the method of the proof of (ii) in Theorem \ref{cor} together with  the Esseen inequality (see \cite{Es68})
\[e^{-C_1\vert t\vert^{2}} \leq \vert f(t)\vert\leq  e^{-C_2\vert t\vert^{2}},\, \vert t\vert\leq A,
\,\,{\rm for}\,\,C_1,C_2,A>0\]
and the following elementary inequality
\[ |z^x - z^y|<|r^x-r^y| +r^{{x+y}\over{2}}|\theta||x-y|,\,\,{\rm where}\,\, z=re^{i\theta}.\]
We get  Corollary \ref{cort1} with $\rho < \frac{1}{2}$ which is still sufficient, by the proof of  Theorem \ref{princ},
to obtain Theorem \ref{t1[asllt]} in a different way.
\end{remark}
}

 \begin{remark} \rm As clearly stated at the beginning of this section, Theorem \ref{princ} holds in the case $\alpha >1$. We believe that this is due to the particular
 arguments used for the proof, and that it is possible to extend the ASLLT also to the case $\alpha<1$. The critical case $\alpha=1$ remains unexplored till now. Another not
 yet investigated situation is
for $\alpha=2$ with $x \mapsto \E[X^2 1_{\{|X|\le  x\}}]$ slowly varying and $\E[X^2]=\infty$ with $x \mapsto x^2{\mathbb P}\{|X| > x\}$ not
slowly varying.
\end{remark}

\bigskip

{\color{green}\subsection{ASLLT for Markov chains}}  \label{2.4}
{\color{red}
The question about ASLLT for dependent random variables was raised by Denker and Koch (\cite{DK}, p. 149, lines -2,-1).
Consider $\{\xi_k\}$ be $0-1$ state Markov chain generated by a $2\times 2$ matrix ${\mathbf P},$ where
\[
{\mathbf P}=
\begin{pmatrix}
 p_{11}&p_{12}\\
 p_{21}&p_{22}
\end{pmatrix},
\]
$p_{ij}>0,$ $i,j=1,2$.
By the analogy to the i.i.d. case let us call $\{\xi_k\}$ Markov trials.
Set
\[ \gamma = 1-p_{12}-p_{21},\qquad
\pi_0 = {{p_{21}}\over{p_{12}+p_{21}}},\qquad
 \pi_1 = {{p_{12}}\over{p_{12}+p_{21}}}\cdotp \]
Consider $\{f(\xi_k)\}$ where $f(0)=-\pi_1,f(1)=\pi_0$. Let $S_n=\sum_{k=1}^n f(\xi_k)$.
It is not difficult to see that the asymptotic (or spectral) variance
$\sigma^2$ of $\{f(\xi_k)\}$ satisfies
\[
\sigma^2
= { E}_{\pi}(f^2(\xi_0))+ 2\sum_{n\geq 1} { E}_{\pi}(f(\xi_0)f(\xi_n))=
\pi_0\pi_1(1+ {{2\gamma}\over{1-\gamma}}) =
\pi_0\pi_1{{1+\gamma}\over{1-\gamma}}\cdotp
\]
It turns out that for Markov trials the following corresponds to Corollary 1 in \cite{DK}
\begin{equation}
\label{admeqn}
{{1}\over{\log{n}}}\sum_{{\nu}=1}^n {{\sigma}\over{\sqrt{{\nu}}}}I_{[S_{\nu} =
\kappa_{\nu}]}
\buildrel{\rm a.s.}\over{\rightarrow}
{{1}\over{\sqrt{2\pi}}}
e^{{-\kappa^2}\over{2}}  {\qquad{\rm as}\qquad}
{{\kappa_{\nu} }\over{\sigma\sqrt{{\nu}}}}\to_{\nu} \kappa,
\end{equation}
where $\kappa_{\nu}$ are of the form $-{\nu}\pi_1+k$.
{\color{green} We refer to   Giuliano, Szewczak
\cite{GS13}}. The relation (\ref{admeqn}) is the immediate consequence of more general result in \cite{GS13}.
The proof of the latter uses ideas from \cite{GW2}. The key role {\color{green} is} play{\color{green}ed by } Edgeworth expansions in
conditional form{\color{green},} which is the special case of more general form with tied up ends (\cite{Vo58}, \cite{Sz08b})).
For example Theorem \ref{moivre} for Markov chains runs as:  {\it For $i,j=0,1$
\[
\sigma\sqrt{2\pi n}\,
{\mathbb P}\{\xi_1+\cdots+\xi_n = k\, ;\, \xi_n = j\,\vert\, \xi_0=i\}
= e^{-{{(k - n\pi_1)^2}\over{2\sigma^2n}}}
\pi_j(1+o(1)),
 \]
uniformly in $k$ such that $\vert{{ k - n\pi_1}\over{\sigma\sqrt{n}}}\vert=o(n^{{1}\over{6}})$}
(see Theorem \ref{stl} and Theorem 2 in \cite{Sz22}).
Limit theorems of this type
proved to be useful in statistics for finding asymptotic formula for the Bayes risk in discriminating between two Markov
chains (see \cite{Na01}).
}

\bigskip
{\color{green}\subsection{ASLLT for the Dickman function}}
\label{2.3}

This section is devoted to the  study the  asymptotic second order properties  of Hwang and  Tsai's  probabilistic model  for   the Dickman function in \cite{HT}.
{\color{green} We refer to   Giuliano, Szewczak and Weber
\cite{GSW}}. We prove a rather delicate correlation inequality for this model and next derive a fine  almost sure local limit theorem.  Recall briefly  that this function
arises from  Dickman's result:   the limit
\begin{equation*} \lim_{n\to \infty} \frac{1}{n}\#\big\{k; 1\le k\le n : P^+(k)\le n^{1/u}\big\}= \varrho(u)
\end{equation*}
exists, where $P^+(n)$ is largest prime factor   of a natural integer $n$.
 This limit is   called the Dickman function, and is   the continuous solution of the differential-difference equation
 $u\varrho'(u) + \varrho(u-1)=0$, $u>1   $
with the initial condition $\varrho(u)=1$ for $0\le u\le 1$.  \vskip 1 pt

Let  $(Z_k)_{k \ge  1}$ be independent random variables defined  by
\begin{equation}\label{variabiliz}
Z_k = \begin{cases}
1 & \hbox {with probability } {1}/{k}\\ 0 & \hbox {with probability }  1- {1}/{k}.
\end{cases}
\end{equation}

  Put for all integers $ m,n $ with $0\le  m <n$, $T_m^n = \sum_{k=m+1}^n k Z_k$,
and set    $T_n=T_0^n $. Then
$$\lim_{n\to \infty} P\big\{ n^{-1}T_n<x\big\}= e^{-\gamma}\int_{0}^x \varrho(v) {\rm d} v, \qq \quad (x>0)  $$
 where $\gamma$ is the Euler--Mascheroni constant ($\int_0^\infty \varrho(v){\rm d} v=e^\gamma$.).
 \smallskip\par   The {\it Dickman distribution} is denoted throughout by  $D$ and  is the distribution function with   density $  e^{-\gamma}\varrho (x) $, $x \ge  0$. It
is known  that $D$ is infinitely divisible.

\vskip 1 pt

 Apart from its  obvious link with number theory,
these type of models also appear as   borderline cases in   the local and almost sure local limit theory.    An important problem inside this theory concerns   the study of
the local and almost sure local limit theorem for weighted sums of Bernoulli  variables. The ``simple'' case when the weights are increasing   is, to say the least,  far from
being understood.
Further,  the well-known Bernoulli part extraction method used for proving  local limit theorems becomes ineffective  in this case.

  \smallskip\par
  The almost sure local limit theorem for the sequence $\{T_n, n\ge 1\}$ states as follows.

\begin{theorem} \label{ASLLT}Let   $\boldsymbol \kappa= (\kappa_n)_{n \ge  1}$ be a strictly increasing sequence such that $\lim_{n \to \infty}\frac{\kappa_n}{n}=x >0$.
Then we have, recalling that $\varrho(x)=1$ if $0\le x\le 1$,
 $$\lim_{N \to \infty} \frac{1}{\log N}\sum_{n=1}^N1_{\{T_n = \kappa_n\}}= e^{- \gamma}\varrho(x), \qquad a.s.$$
  \end{theorem}

 \begin{corollary} \label{corASLLT}We have
 $$\lim_{N \to \infty} \frac{1}{\log N}\sum_{n=1}^N1_{\{T_n = n\}}= e^{- \gamma}, \qquad a.s.$$
As a consequence, for every $x\ge  1$,
 $$\lim_{N \to \infty} \frac{\sum_{n=1}^N1_{\{T_n = [xn]\}}}{\sum_{n=1}^N1_{\{T_n = n\}}}= \varrho(x), \qquad a.s.$$
 \end{corollary}
The key tools for proving Theorem \ref{ASLLT} are the following correlation inequality (Theorem \ref{covarianzabase})
and integral convergence of characteristic functions (Proposition \ref{ZS}).

\bigskip\noindent
\subsubsection{Correlation inequality and integral convergence of characteristic functions }

\medskip

\begin{theorem} \label{covarianzabase} Let  ${\boldsymbol \kappa }=(\kappa_n)_{n \ge  1}$ be a  sequence of integers such that  $\lim_{n \to \infty} \frac{\kappa_n}{n}=x>0$.
Let
 $Y_n = n 1_{\{T_n = \kappa_n\}}$, $n\ge 1$.
 Then
 there exists a positive constant $C$, such that for any $n> m \ge  2$,
\begin{align*}
&|Cov(Y_m, Y_n)|\le  C \bigg\{\frac{n}{n-m}\chi^{({\boldsymbol \kappa },x)}_{m,n}
+\frac{m}{n-m}+\chi^{({\boldsymbol \kappa },x)}_{2,n}+\frac{1}{n}\bigg\}.
\end{align*}
  and
$$\chi_{m,n}^{({\boldsymbol \kappa },x)}=
\frac{n-m}{\kappa_n-\kappa_m}\cdot\frac{\log\frac{n}{m}}{\sqrt {n-m}}+ \frac{n-m}{\kappa_n-\kappa_m}\cdot g_{m,n}
+x\Big|\frac{n-m}{\kappa_n-\kappa_m}-\frac{1}{x}\Big|+ \frac{m+1}{\kappa_n-\kappa_m}.$$
\end{theorem}

  \smallskip\par \noindent

Recall that the characteristic function of the Dickman distribution satisfies
\begin{equation}\label{characteristicfunction}
\phi(t)=\exp \Big\{\int_0^1 \frac{e^{itu}-1}{u}\,{\color{green} du}\Big\}  .\end{equation}
See \cite{HT}. Consider also the characteristic function of $Z_k$  and $T_m^n$ defined in \eqref{variabiliz} and right after. We have
\begin{eqnarray*} \phi_{Z_k}(t)&=& 1+ \frac{e^{it}-1}{k}
\cr \phi_{T_m^n}(t)&=&\prod_{k=m+1}^n\phi_{Z_k}(tk)=\prod_{k=m+1}^n\Big(1+ \frac{e^{itk}-1}{k}\Big).
\end{eqnarray*}

The following result    is of crucial use  for  the proof of Proposition \ref{stimabase}, which together with the above correlation inequality
will lead to the proof of Theorem \ref{ASLLT}.

\begin{proposition}\label{ZS} We have
\begin{eqnarray*}
{\rm(a)} & & \int_{-\infty}^{+\infty}\big|\phi(t)\big|^2 \,dt <  \infty,
\cr {\rm(b)} & &\int_{-\pi n}^{\pi n}\big|\phi_{\frac{T_n}{n}}(u)\big|^2\,du\to \int_{-\infty}^{\infty}|\phi(u)|^2\,du, \qquad {\text as}\quad n \to \infty.\end{eqnarray*}
\end{proposition}
Another consequence of Proposition \ref{ZS} is the following local limit theorem

  \begin{corollary} {{\rm ({\rm Local limit theorem})}} \label{local}
 Let $\boldsymbol \kappa=(\kappa_n)_{n\ge  1}$ be a  sequence of integers such that  $\lim_{n \to \infty}\frac{\kappa_n}{n}=x >0$. Then
 $$\lim_{n \to \infty}n{\mathbb P}(T_n = \kappa_n)= e^{- \gamma}\varrho(x).$$
 \end{corollary}
 The  strong  form of the local limit theorem for the Dickman density states as follows.

\begin{theorem}[{\rm Strong  Local limit theorem}] \label{localIL}
 \begin{equation}
\label{sllt}
\sum_{\kappa\in{\mathbb N}} |{\mathbb P}(T_n = \kappa)- n^{-1}e^{-\gamma}\varrho(n^{-1}\kappa)|\to 0, \quad n \to \infty.
\end{equation}
 \end{theorem}
A sharper speed of convergence result was recently obtained in  La  Bret\`eche and Tenenbaum \cite{BT}. The proof  however  much appeals to   analytic number theory, details
are not given, and so seems more reserved to specialists of this area.

\vskip 15 pt
 The  following result is an extension of   the one given in \cite{HT} (Proposition 2.1) for the case $m_n \equiv 0$.

\begin{proposition}\label{law} Let $(m_n)_{n \ge  1}$ be a sequence of integers such that $\lim_{n \to \infty}(n-m_n)= +  \infty$ and $m_n=o(n)$. Then, as $n \to \infty$,
 the sequence $\frac{T_{m_n}^n}{n-m_n}$ converges in distribution to the Dickman law.
 \end{proposition}

  The next  result
 will be crucial for the proof of the correlation inequality.

\begin{proposition}\label{stimabase}
Let $\boldsymbol \kappa=(\kappa_n)_n$ be an  increasing sequence of integers. Then, for $n> m \geq 2$,
\begin{align*}\Big|(\kappa_n- \kappa_m){\mathbb P}(T_m^n = \kappa_n- \kappa_m)-&
 {\mathbb P}\big((\kappa_n- \kappa_m)-n < T_m^n \le  (\kappa_n- \kappa_m)-(m+1)\big) \Big|
 \cr &\le  C  \frac{1+\log\frac{n}{m}}{\sqrt{n-m}}.
 \end{align*}
\end{proposition}

In the Proposition below, we first  specify  Proposition 1 of \cite{HT} quantitatively in terms of the characteristic functions.
\begin{proposition}\label{W1}
There exists an absolute constant $C$ such that for all integers $n > m \ge  2$ and all real numbers $t$,
$$\bigg|\phi_{\frac{T_m^n}{n-m}}(t)-\exp\Big\{\int_0^1 \frac{e^{itu}-1}{u}\,du\Big\}\bigg|\le  f_{m,n}(t),$$ where
$$ f_{m,n}(t)=\exp\Bigg\{Ct^2\Big(\frac{\log\frac{n}{m}}{(n-m)^2}+\frac{m+2}{n-m}\Big)\Bigg\}-1 .$$
\end{proposition}

 \vskip 3 pt

The following result specifies Proposition 1 of \cite{HT} quantitatively in terms of distribution functions.

\begin{proposition}\label{W2}
There exists an absolute positive constant $C$ such that, for all positive integers $n$, $m$, with $n>m \ge  2$,
\begin{equation*}\sup_{x \in \mathbb{R}}\Big|{\mathbb P}\Big(\frac{T_m^n}{n-m}\le  x\Big)- D(x)\Big|\le  C g_{m,n},\end{equation*} where
\begin{equation*}g_{m,n}=\exp\Bigg(C\Big\{\frac{\log\frac{n}{m}}{(n-m)^2}+ \frac{m+2}{n-m}\Big\}\log^2\frac{n}{m}\Bigg)-1+ \frac{1}{\log\frac{n}{m}}.
\end{equation*}
\end{proposition}
 \vskip 15 pt

The detailed proofs of the various ASLLT presented before, as well as  new ones, are the objet of a  joint work in preparation.}


\vfill\break

$${}$$
\vskip 80 pt

{\small    \leftskip =3cm \rightskip=3cm $$\hbox{\bf Summary table}$$
 \vskip 30pt  \noi$\hbox{\ninerm 1. The local limit theorem.}$
\eightrm De Moivre-Laplace's Theorem\,$\cdot$\,Lattice-valued random variables.
 $\hbox{\sevenit The i.i.d. case}$:
Gnedenko's Theorem\,$\cdot$\,Method of   characteristic functions
\,$\cdot$\,Local limit theorems for weighted sums of i.i.d. random variables\,$\cdot$\,Local large deviations\,$\cdot$\,Diophantine measures and local limit
theorem\,$\cdot$\,Local limit theorems with arithmetical   conditions.
 $\hbox{\sevenit The independent case}$:
Prohorov's Theorem\,$\cdot$\,A general necessary condition\,$\cdot$\,Gamkrelidze's lower bound\,$\cdot$\,Stable limit distributions\,$\cdot$\,\hfill Gamkrelidze's
counter-examples\,$\cdot$\,Integral limit theorem and local limit theorem\,$\cdot$\,Mu-\hfill
khin's necessary and sufficient condition\,$\cdot$\,Characteristics of a random variable\,$\cdot$\,The Bernoulli part of a random variable\,$\cdot$\,A local limit theorem
with effective rate\,$\cdot$\,Local limit theorems and the Landau--Kolmogorov inequalities\,$\cdot$\,Mod-$\phi$ convergence and local limit theorem\,$\cdot$\,Local limit
theorem under tightness conditions\,$\cdot$\,Local limit theorems and domains of attraction\,$\cdot$\,Refinements of the Gnedenko local limit theorem\,$\cdot$\,Local limit
theorem for the number of renewals\,$\cdot$\,Local limit theorems for densities and refinements\,$\cdot$\,Local limit theorem for residues of linear transforms of sums.
 $\hbox{\sevenit The ergodic case}$:
Local limit theorem for interval expanding maps\,$\cdot$\,The continued fractions case\,$\cdot$\,Local limit theorems for Riesz-Raikov sums.
 $\hbox{\sevenit Miscellaneous results}$: Linear random fields\,$\cdot$\,Operator form\,$\cdot$\,Nonconventional sums\,$\cdot$\,A generalization of Richter's theorem.
 $\hbox{\sevenit Poisson approximation--Another coupling method}$: Poisson approximation to   binomial distribution\,$\cdot$\,Le Cam's inequality\,$\cdot$\,Measures of
 disparity.

 \noi $\hbox{\sevenit Some applications}$: Asymptotic enumeration\,$\cdot$\,Coefficients of polynomials of binomial type\,$\cdot$\,Asymptotics for Legendre
 polynomials\,$\cdot$\,Allocations problems\,$\cdot$\,Random permutations\,$\cdot$\,Random mappings\,$\cdot$\,Number theory.
\vskip 2 pt
\noi$\hbox{\ninerm 2. Almost sure versions.}$ The Chung-Erd\H{o}s ASLLT\,$\cdot$\,ASLLT for i.i.d. square integrable lattice  random variables\,$\cdot$\,  ASLLT for   i.i.d.
stable lattice random variables\,$\cdot$\,ASLLT for Markov chains\,$\cdot$\,ASLLT for the Dickman function.
\vskip 2 pt
 \hfill
\vskip 1pt \hskip 20pt
}

\end{document}